\documentclass[reqno]{amsart}
\usepackage[T1]{fontenc}
\usepackage{amsfonts}
\usepackage{amsmath}
\usepackage{amsthm}
\usepackage{amssymb}
\usepackage{graphicx}
\usepackage{subcaption}
\usepackage{url}
\usepackage{hyperref}
\usepackage{geometry}
\newtheorem{theorem}[subsection]{Theorem}
\newtheorem{definition}{Definition}[section]
\newtheorem{lemma}{Lemma}[section]
\newtheorem{proposition}{Proposition}[section]
\newtheorem{corollary}{Corollary}[section]
\newtheorem{remark}{Remark}[section]
\newtheorem{example}{Example}[section]

\title[Uniqueness of the strong positive solution]{Uniqueness of the strong positive solution for a general quasilinear elliptic problem with variable exponents and homogeneous Neumann boundary conditions using a generalization of the $p(x)$-D\' iaz-Saa inequality}
\author{Maxim Bogdan}
\email{maxim.bogdan.n6h@student.ucv.ro}
\date{\today}

\begin{document}

\maketitle

	\begin{center} \footnotesize{Department of Mathematics, University of Craiova, Al. I. Cuza Street, no. 13, 200585, Craiova, Romania}
\end{center}

	\begin{abstract}
	 In this paper, we study a generalization of the Díaz-Saa inequality and its applications to nonlinear elliptic problems. We first present the necessary hypotheses and preliminary results before introducing an improved version of the inequality, which holds in a broader functional setting and allows applications to problems with homogeneous Neumann boundary conditions. The significance of cases where the inequality becomes an equality is also analyzed, leading to uniqueness results for certain classes of partial differential equations. Furthermore, we provide a detailed proof of a uniqueness theorem for strong positive solutions and illustrate our findings with two concrete applications: a multiple-phase problem and an elliptic quasilinear equation relevant to image processing. The paper concludes with possible directions for future research.
	\end{abstract}
	
\bigskip

{\footnotesize{\textbf{Keywords}: Multiple--phase elliptic problem, variable exponents, homogeneous Neumann boundary conditions, weak solution, D\'{i}az-Saa inequality, uniqueness of the steady-state, image processing}}

\smallskip

{\footnotesize{\textbf{MSC 2020}: 35A01, 35A02, 35B09, 35D30, 35J20, 35J62, 35J92, 98A08}}

\section{Introduction}\label{sec1}

This paper is the natural continuation of the recent work \cite{max1}, so we will study the same problem introduced there:

\begin{equation}\tag{$E$}\label{eqmphase}
	\begin{cases}-\operatorname{div}\mathbf{a}\big (x,\nabla U(x)\big )=f\big (x,U(x)\big ), & x\in\Omega\\[3mm] \dfrac{\partial U}{\partial\nu}=0, & x\in \partial\Omega\\[3mm] 0\leq U(x)\leq 1, & x\in\Omega\end{cases}
\end{equation}

\noindent From \cite[Theorem 7.1]{max1} we know that problem \eqref{eqmphase} admits a minimal solution $\underline{U}$ and a maximal solution $\overline{U}$ in the sense of Definition \ref{def2weaksol} and if $\underline{U}=\overline{U}$ then there is no other solution. In this research we are interested in proving the uniqueness of a strongly positive weak solution $U$, i.e. a weak solution with $\underset{\Omega}{\operatorname{ess\ inf}}\ U>0$.

In Section 2, we provide the list of hypotheses and notations used throughout the paper. In Section 3, we present some preliminary technical results related to the Díaz-Saa inequality. Section 4 has an original character, as it introduces the generalization of the Díaz-Saa inequality. This inequality is extremely important when addressing the problem of uniqueness or nonexistence of solutions for an elliptic partial differential equation with boundary conditions.

Its history dates back to 1986, when it was used for the Laplacian in the work \cite[page 57]{brezosw}. The following year, a generalization of this inequality for the $p$-Laplacian appeared in the paper \cite[Lemme 2]{diazsaa} (from which it took its definitive name). The inequality was also employed for the $p$-Laplacian in \cite[page 145]{fleck} and \cite[page 292]{tak2}. Up to that point, no detailed proof of the result had been provided, only a few lines, omitting a multitude of technical details. There are also small errors due to the absence of a complete proof. Moreover, it had only been stated for functions in the space $W^{1,p}_0(\Omega)$ as it was applied only for homogeneous Dirichlet boundary conditions.

With the explosive growth of interest in studying problems with variable exponents after 2005, the inequality was further generalized to Sobolev spaces with variable exponents, specifically to the space $W^{1,p(x)}_0(\Omega)$, by Peter Takáč and Jacques Giacomoni in the article \cite[Theorem 2.4]{tak}. This work presents the first slightly more elaborated proof of the inequality, though many technical details (which are far from trivial) are left to the reader. In fact, the article \cite{tak} was the one that helped and inspired me the most in writing this paper.

Although we provide a significantly improved version of the D\'{i}az-Saa inequality, which holds throughout the entire space $W^{1,p(x)}(\Omega)$, making it applicable to elliptic problems with homogeneous Neumann boundary conditions, it is nonetheless surprising that the proof is immediate and elementary. The inequality holds for a wide variety of differential operators (see Section 8).

The Díaz-Saa inequality is not as important as the situations where it becomes an equality. These cases allow us to draw conclusions regarding the uniqueness of solutions for various partial differential problems. All such situations are addressed in Section 5, where we introduce a new regularity condition without which qualitative conclusions about the equality cannot be drawn, as illustrated by the examples \ref{ex51} and \ref{ex52}. Theorem \ref{thmequalitydiaz} and Lemma \ref{lemmamea} in this section represent central results of the present work.

Section 6 is quite extensive and contains all the necessary technical details required to establish a uniqueness theorem for strong positive solutions: Theorem \ref{thmuniq}.

In Sections 7 and 8, we present two illustrative examples in which the results obtained throughout the paper can be applied. The first example consists of a general multiple-phase problem, while the second is an elliptic problem with practical significance in image processing. The results from the previously cited articles do not apply to this latter example, as shown in Remark 8.1.

At the end of the paper, we outline several directions for future research. The article concludes with a very brief appendix and the relevant bibliography. All proofs are given in full detail, making the article highly suitable for beginners in the field of partial differential equations.

\section{Hypotheses and notations}

\noindent We work under the following hypotheses, which are basically the same as those from \cite[Section 2]{max1} with the exceptions \textbf{(H7)},\textbf{(H7')},\textbf{(H13)} and \textbf{(H13')}.

\begin{enumerate}
	\item[\textbf{(H1)}] $\Omega\subset\mathbb{R}^N,\ N\geq 2$ is an open, bounded and connected Lipschitz domain.
	
	\bigskip
	
	\item[$\bullet$] $\mathcal{U}:=\big\{U\in L^{\infty}(\Omega)\ |\ 0\leq U\leq 1\ \text{a.e. on}\ \Omega \big \}$.
	
	\bigskip
	
	\item[\textbf{(H2)}] $p:\overline{\Omega}\to (1,\infty)$ is a log-H\"{o}lder continuous variable exponent with the property that: $p^-\geq\dfrac{2N}{N+2}$ or $p:\overline{\Omega}\to (1,\infty)$ is a continuous exponent with $p^->\dfrac{2N}{N+2}$. 
	\bigskip
	
	\item[$\bullet$] Denote $p^{-}=\displaystyle\min_{x\in\overline{\Omega}}\ p(x) >1$ and $p^+=\displaystyle\max_{x\in\overline{\Omega}}\  p(x)<\infty$.
	
	\bigskip
	
	\item[\textbf{(H3)}] $\Psi:\overline{\Omega}\times (0,\infty)\to (0,\infty)$ with $\Psi(\cdot, s)$ measurable for each $s\in (0,\infty)$ and $\Psi(x,\cdot)\in\operatorname{AC}_{\text{loc}}\big ((0,\infty)\big )$ for a.a. $x\in\Omega$.

	\item[\textbf{(H4)}] $\lim\limits_{s\to 0^+} \Psi(x,s)s=0$ for each $x\in\Omega$.
	
	\item[\textbf{(H5)}] For a.a. $x\in\Omega$ we have that $(0,\infty)\ni s\longmapsto \Psi(x,s)s\in (0,\infty)$ is a strictly increasing function.
	
	\bigskip

	\item[$\bullet$] We take $\mathbf{a}:\overline{\Omega}\times\mathbb{R}^N\to\mathbb{R}^N,\ \mathbf{a}(x,\xi)=\Psi(x,|\xi|)\xi$ if $\xi\neq \mathbf{0}$ and $\mathbf{a}(x,\textbf{0})=\textbf{0}$.
	
	\bigskip
	
	\item[$\bullet$] Define $\Phi:\overline{\Omega}\times\mathbb{R}\to [0,\infty),\ \Phi(x,s)=\begin{cases} \Psi(x,|s|)|s|, & s\neq 0\\ 0, & s=0\end{cases}$. Notice that: $|\mathbf{a}(x,\xi)|=\Phi(x,|\xi|)$ and $\mathbf{a}(x,\xi)\cdot\xi=\Phi(x,|\xi|)|\xi|,\ \forall\ x\in\overline{\Omega},\ \forall\ \xi\in\mathbb{R}^N$.
	
	\bigskip
	
	\item[\textbf{(H6)}] There is some function $a\in L^{p'(x)}(\Omega)^+$ and a constant $b\geq 0$ such that: 
	
	\begin{equation}
		\Phi(x,s)\leq a(x)+b|s|^{p(x)-1},\ \forall\ s\in [0,\infty),\ \text{for a.a.}\ x\in\Omega.
	\end{equation}
	
	\item[\textbf{(H7)}] There is some $r\geq 1$ such that the function $(0,\infty)\ni s\mapsto\dfrac{\Phi(x,s)}{s^{r-1}}$ is increasing for a.a. $x\in\Omega$.
	
	\item[\textbf{(H7')}] There is some $r\geq 1$ such that the function $(0,\infty)\ni s\mapsto\dfrac{\Phi(x,s)}{s^{r-1}}$ is \textbf{strictly increasing} for a.a. $x\in\Omega$.
	\bigskip
	
	\item[$\bullet$] $A:\overline{\Omega}\times\mathbb{R}^N\to [0,\infty),\ A(x,\xi)=\displaystyle\int_{0}^{|\xi|}\Phi(x,s)\ ds$.

	\item[$\bullet$] $\mathcal{A}:W^{1,p(x)}(\Omega)\to [0,\infty)$ given by $\mathcal{A}(u)=\displaystyle\int_{\Omega} A(x,\nabla u(x))\ dx$.
	
	\bigskip
	
	\item[\textbf{(H8)}] There exists some constant $\delta_0>0$ (which can be very small) and some constant $\tilde{\delta}_0\geq 0$ (which can be very big) such that:
	
	\begin{equation}
		\mathcal{A}(v)=\int_{\Omega} A(x,\nabla v(x))\ dx\geq \delta_0\int_{\Omega} |\nabla v(x)|^{p(x)}\ dx-\tilde{\delta}_0,\ \forall\ v\in W^{1,p(x)}(\Omega).
	\end{equation}
	
	\item[\textbf{(H9)}] The function $f:\overline{\Omega}\times [0,1]\to\mathbb{R}$ is measurable.
	
	\item[\textbf{(H10)}] The function $\Omega\ni x\longmapsto f(x,0)$ is from $L^{\infty}(\Omega)$.
	
	\item[\textbf{(H11)}] $f(x,0)\geq 0$ and $f(x,1)\leq 0$ for a.a. $x\in\Omega$. 
	
	\item[\textbf{(H12)}] There are some constants $\lambda_0,\gamma>0$ such that for a.a. $x\in\Omega$, the function $[0,1]\ni s\mapsto f(x,s)+\lambda_0 s$ is strictly increasing and the function $[0,1]\ni s\mapsto f(x,s)$ is $\gamma$--Lipschitz.
	
	\item[$\bullet$] $F:\overline{\Omega}\times[0,1]\to\mathbb{R},\ F(x,s)=\displaystyle\int_{0}^s f(x,\tau)\ d\tau$.
	
	\item[\textbf{(H13)}] There is some constant $\min\{p^{-},2\}\geq\alpha>1$ such that the function $(0,1]\ni s\longmapsto \dfrac{f(x,\sqrt[\alpha]{s})}{\sqrt[\alpha]{s^{\alpha-1}}}$ is decreasing.
	
	\item[\textbf{(H13')}] There is some constant $\alpha\in (1,2)$ with $p^-\geq\alpha$ such that the function $(0,1]\ni s\longmapsto \dfrac{f(x,\sqrt[\alpha]{s})}{\sqrt[\alpha]{s^{\alpha-1}}}$ is \textbf{strictly decreasing}.

\end{enumerate}

\bigskip

\noindent $\blacktriangleright$ Hypotheses \textbf{(H7')} and \textbf{(H13')} will be assumed only when this is explicitly stated.

\begin{remark}
	The two hypotheses \textnormal{\textbf{(H12)}} and \textnormal{\textbf{(H13)}} when assumed together forces the source to be negative if $f(x,\cdot)$ is differentiable and $f(x,0)=0$ for a.a. $x\in\Omega$. Indeed, since from \textnormal{\textbf{(H13)}} we have for each $s\in (0,1]$ that:
	
	\begin{equation}
		\dfrac{f(x,\sqrt[\alpha]{s})}{\sqrt[\alpha]{s^{\alpha-1}}}\leq \lim\limits_{s\to 0^+}\dfrac{f(x,\sqrt[\alpha]{s})}{\sqrt[\alpha]{s^{\alpha-1}}}\stackrel{s=\tau^{\alpha}}{=}\lim\limits_{\tau\to 0^+}\dfrac{f(x,\tau)}{\tau^{\alpha-1}}=\lim\limits_{\tau\to 0^+}\dfrac{1}{\alpha-1}\underbrace{\dfrac{\partial f}{\partial\tau}(x,\tau)}_{\textnormal{bounded}}\tau^{2-\alpha}=0,\ \text{for}\ \alpha<2,
	\end{equation}
	
	\noindent we deduce that $f(x,\tau)\leq 0$ for any $\tau\in (0,1]$. We have used \textnormal{\textbf{(H12)}} to obtain that $\dfrac{\partial f}{\partial\tau}(x,\tau)$ is bounded (as $f(x,\cdot)$ is $\gamma$--Lipschitz). A better requirement in \textnormal{\textbf{(H12)}} would be to have $f(x,\cdot):[0,1]\to\mathbb{R}$ to be locally-Lipschitz, but this will be a topic for a further research (with this assumption we can't use the existence theorem proved in \cite[page 23]{max1}).
\end{remark}

\begin{remark} Hypotheses \textnormal{\textbf{(H9)}}, \textnormal{\textbf{(H10)}}, \textnormal{\textbf{(H11)}}, \textnormal{\textbf{(H12)}} and \textnormal{\textbf{(H13)}} are satisfied for a negative source $f:\overline{\Omega}\times [0,1]\to\mathbb{R},\ f(x,s)=-r_1(x)s^{q_1(x)}-r_2(x)s^{q_2(x)}$, where $r_1,r_2\in L^{\infty}(\Omega),\ r_1,r_2\geq 0$ a.e. on $\Omega$ and $q_1,q_2:\Omega\to [1,\infty)$ are two measurable and bounded variable exponents.
	
\bigskip
	
\noindent Indeed, $f$ is clearly measurable, $f(x,0)=0\in L^{\infty}(\Omega)$ and $f(x,0)=0\geq 0$ a.e. on $\Omega$. Note that $f(x,1)=-r_1(x)-r_2(x)\leq 0$ a.e. on $\Omega$, $\left |\dfrac{\partial f}{\partial s}(x,s)\right |=r_1(x)q_1(x)s^{q_1(x)-1|}+r_2(x)q_2(x)s^{q_2(x)-1}\leq r_1(x)q_1(x)+r_2(x)q_2(x)\leq \Vert r_1\Vert_{L^{\infty}(\Omega)}\Vert q_1\Vert_{L^{\infty}(\Omega)}+\Vert r_2\Vert_{L^{\infty}(\Omega)}\Vert q_2\Vert_{L^{\infty}(\Omega)}:=\gamma$. So $[0,1]\ni s\mapsto f(x,s)$ is $\gamma$--Lipschitz for a.a. $x\in\Omega$. We may choose any $\lambda_0>\gamma$ to obtain that $[0,1]\ni s\mapsto f(x,s)+\lambda_0 s$ is a strictly increasing function for a.a. $x\in\Omega$, as pointed in \cite[Remark 2.2]{max1}.
\end{remark}

\begin{remark}
	In the most important results of the paper, like Theorems \ref{thmbeta} and \ref{thmuniq}, we will need \textnormal{\textbf{(H7)}} or \textnormal{\textbf{(H7')}} to hold for $r=\alpha$.
\end{remark}

\bigskip

\noindent In order to use a variational approach for problem \eqref{eqmphase} we need first to extend the domain of the function $f$ in such a way that it will not lose its properties. The idea of the below definition is taken from \cite[Theorem 3.1, page 102]{gariepy}.

\begin{definition} We introduce the following functions $\overline{f},\overline{F}:\overline{\Omega}\times\mathbb{R}\to\mathbb{R}$ given by:
	
	\begin{equation}
		\overline{f}(x,s)=\sup_{\tau\in [0,1]}f(x,\tau)-\gamma |s-\tau|=\begin{cases} f(x,0)+\gamma s, & s\in (-\infty,0)\\ f(x,s), & s\in [0,1] \\ f(x,1)-\gamma(s-1), & s\in (1,\infty) \end{cases}
	\end{equation}
	
	\noindent and:
	
	\begin{equation}
		\overline{F}(x,s)=\int_{0}^s\overline{f}(x,\tau)\ d\tau=\begin{cases} f(x,0)+\frac{\gamma}{2}s^2& s<0 \\[3mm] F(x,s)=\displaystyle\int_{0}^s f(x,\tau)\ d\tau, &  s\in [0,1] \\[3mm] F(x,1)+f(x,1)(s-1)-\frac{\gamma}{2}(s-1)^2, & s>1 \end{cases}.
	\end{equation}
\end{definition}

\begin{proposition}\label{propf} The following properties of $\overline{f}$ and $\overline{F}$ hold for a.a. $x\in\Omega$:
	\begin{enumerate}
		\item[\textnormal{\textbf{(1)}}] $\overline{f}(x,\cdot):\mathbb{R}\to\mathbb{R}$ is $\gamma-$Lipschitz;
		
		\item[\textnormal{\textbf{(2)}}] For any $\overline{\lambda}>\gamma$ the function $\mathbb{R}\ni s\longmapsto\overline{f}(x,s)+\overline{\lambda}s$ is strictly increasing for a.a. $x\in\Omega$;
		
		\item[\textnormal{\textbf{(3)}}]  The function $[0,\infty)\ni s\mapsto -\overline{F}(x,\sqrt[\alpha]{s})$ is convex and the function $(0,\infty)\ni s\mapsto\dfrac{\overline{f}(x,\sqrt[\alpha]{s})}{\sqrt[\alpha]{s^{\alpha-1}}}$ is decreasing for a.a. $x\in\Omega$.
		
		\item[\textnormal{\textbf{(4)}}]  If \textnormal{\textbf{(H13')}} holds then the function $[0,\infty)\ni s\mapsto -\overline{F}(x,\sqrt[\alpha]{s})$ is \textbf{strictly convex} and the function $(0,\infty)\ni s\mapsto\dfrac{\overline{f}(x,\sqrt[\alpha]{s})}{\sqrt[\alpha]{s^{\alpha-1}}}$ is \textbf{strictly decreasing} for a.a. $x\in\Omega$.
	\end{enumerate}
	
\end{proposition}

\begin{proof}\textbf{(1)} Consider some real numbers $s_1,s_2$. We have that: $\overline{f}(x,s_1)=\displaystyle\sup_{\tau\in [0,1]} f(x,\tau)-\gamma|s_1-\tau|=\displaystyle\sup_{\tau\in [0,1]} f(x,\tau)-\gamma|(s_2-\tau)-(s_2-s_1)|\leq\displaystyle\sup_{\tau\in [0,1]} f(x,\tau)-\gamma|s_2-\tau|+\gamma |s_2-s_1|=\overline{f}(x,s_2)+\gamma |s_2-s_1|$. So $\overline{f}(x,s_1)-\overline{f}(x,s_2)\leq \gamma |s_2-s_1|$. By swapping $s_1$ and $s_2$ we get that $\overline{f}(x,s_2)-\overline{f}(x,s_1)\leq \gamma |s_1-s_2|$. Thus, $|\overline{f}(x,s_2)-\overline{f}(x,s_1)|\leq\gamma |s_2-s_1|$, as desired.
	
	\noindent\textbf{(2)} Let $s_1>s_2$ be two real numbers. Then $\overline{f}(x,s_1)+\overline{\lambda}s_1-\overline{f}(x,s_2)-\overline{\lambda}s_2=\overline{f}(x,s_1)-\overline{f}(x,s_2)+\overline{\lambda}(s_1-s_2)\geq -\gamma (s_1-s_2)+\overline{\lambda}(s_1-s_2)=(\overline{\lambda}-\gamma)(s_1-s_2)>0$.
	
	\noindent\textbf{(3)} Since $\overline{f}(x,\cdot)$ is continuous (being $\gamma$--Lipschitz) we deduce that its antiderivative $\overline{F}(x,\cdot)$ is continuously differentiable on $\mathbb{R}$.
	Thus, from the chain rule, for any $s>0$, we get that: $\dfrac{\partial }{\partial s}\overline{F}(x,\sqrt[\alpha]{s})=\overline{f}(x,\sqrt[\alpha]{s})\dfrac{1}{\alpha\sqrt[\alpha]{s^{\alpha-1}}}=\dfrac{1}{\alpha}\begin{cases}\dfrac{f(x,\sqrt[\alpha]{s})}{\sqrt[\alpha]{s^{\alpha-1}}}, & s\in (0,1]\\[3mm] \dfrac{f(x,1)-\gamma(\sqrt[\alpha]{s}-1)}{\sqrt[\alpha]{s^{\alpha-1}}}, & s>1 \end{cases}$. From \textbf{(H13)} we get that this function is decreasing on $(0,1]$. Being continuous at $s=1$ we only need to check if the second branch defines a decreasing function too. Note that for any $s>1$ we have:
	
	\begin{equation}
		\dfrac{\partial}{\partial s}\left [ \dfrac{f(x,1)-\gamma(\sqrt[\alpha]{s}-1)}{\sqrt[\alpha]{s^{\alpha-1}}}\right ]=\dfrac{s^{\frac{1}{\alpha}-2}}{\alpha}\big [\underbrace{(1-\alpha)}_{<0}\underbrace{\big (f(x,1)+\gamma \big )}_{\geq 0} -\underbrace{(2-\alpha)}_{\geq 0}\sqrt[\alpha]{s}\big ]\leq 0.
	\end{equation}
	
	\noindent This is true, since from \textbf{(H11)} and \textbf{(H12)}, we get that: $f(x,0)-f(x,1)=|f(x,1)-f(x,0)|\leq \gamma |1-0|=\gamma$. So $f(x,1)+\gamma\geq f(x,0)\geq 0$. We conclude that $-\dfrac{\partial }{\partial s}\overline{F}(x,\sqrt[\alpha]{s})\geq 0$ for any $s>0$. Therefore $-\overline{F}(x,\sqrt[\alpha]{s})$ is convex with respect to $s\in [0,\infty)$.
	
	\noindent\textbf{(4)} In a similar manner as above we obtain $-\dfrac{\partial }{\partial s}\overline{F}(x,\sqrt[\alpha]{s})>0$ for any $s>0$. It is essential that $\alpha<2$. Therefore $-\overline{F}(x,\sqrt[\alpha]{s})$ is strictly convex with respect to $s\in [0,\infty)$.
\end{proof}

\begin{definition}\label{def2weaksol}
	We say that $U\in W^{1,p(x)}(\Omega)$ is a \textbf{weak solution} for \eqref{eqmphase} if $1\geq U\geq 0$ a.e. on $\Omega$ (i.e. $U\in\mathcal{U})$ and for any test function $\phi\in W^{1,p(x)}(\Omega)$ we have that:
	\begin{equation}
		\int_{\Omega} \mathbf{a}(x,\nabla U(x))\cdot \nabla\phi(x)\ dx=\int_{\Omega} f(x,U(x))\phi(x)\ dx.
	\end{equation}
\end{definition}

\begin{definition} We consider the energy functional $\mathcal{J}:W^{1,p(x)}(\Omega)\to\mathbb{R}$ associated to problem \eqref{eqmphase}, that is defined by:
	
	\begin{equation}
		\mathcal{J}(U)=\int_{\Omega} A(x,\nabla U(x))\ dx-\int_{\Omega} \overline{F}(x,U(x))\ dx,\ \forall\ U\in W^{1,p(x)}(\Omega).
	\end{equation}	
\end{definition}

\begin{proposition}\label{propc1j} We have that $\mathcal{J}\in C^1\big (W^{1,p(x)}(\Omega) \big )$ and:
	
	\begin{equation}
		\langle \mathcal{J}'(U),\phi\rangle=\int_{\Omega} \mathbf{a}(x,\nabla U(x))\cdot\nabla\phi(x)\ dx-\int_{\Omega} \overline{f}(x,U(x))\phi(x)\ dx,\ \forall\ \phi\in W^{1,p(x)}(\Omega).
	\end{equation}
	
	\noindent In particular $U\in W^{1,p(x)}(\Omega)$ is a \textbf{weak solution} of \eqref{eqmphase} iff $U\in\mathcal{U}$ and $\mathcal{J}'(U)\equiv 0$.
	
\end{proposition}

\begin{proof} The proof is standard and it is made in full details in \cite[Proposition 3.2]{max1}.
	
\end{proof}

\begin{definition} Consider the following two sets: $W=\left \{w\in L^{\frac{p(x)}{\alpha}}(\Omega)\ \big |\ \sqrt[\alpha]{|w|}\in W^{1,p(x)}(\Omega) \right\}$ and $\overset{\bullet}{W}=\left\{w\in W\ \big |\ w>0\ \text{a.e. on}\ \Omega\right\}$.
\end{definition}

\section{Preliminary results}

\begin{proposition} The subset $\overset{\bullet}{W}\subset W$ has the following properties:
	
	\begin{enumerate}
		\item[\textnormal{\textbf{(1)}}] If $w\in\overset{\bullet}{W}$ then $\theta w\in\overset{\bullet}{W}$ for any $\theta\in (0,\infty)$.
		
		\item[\textnormal{\textbf{(2)}}] If $w_1,w_2\in \overset{\bullet}{W}$ then $w_1+w_2\in \overset{\bullet}{W}$.
		
		\item[\textnormal{\textbf{(3)}}] $\overset{\bullet}{W}$ is a convex subset of $W$.
	\end{enumerate}
	
\end{proposition}

\begin{proof} \textbf{(1)} From $\alpha\leq p^{-}$ we have that $\dfrac{p(x)}{\alpha}\geq\dfrac{p(x)}{p^-}\geq 1$ for a.a. $x\in\Omega$. Thus $L^{\frac{p(x)}{\alpha}}(\Omega)$ is a Banach space. So knowing that $w\in L^{\frac{p(x)}{\alpha}}(\Omega)$ implies that $\theta w$ is also in the same space. Moreover from $\sqrt[\alpha]{w}\in W^{1,p(x)}(\Omega)$ we get that $\sqrt[\alpha]{\theta w}=\sqrt[\alpha]{\theta}\sqrt[\alpha]{w}\in W^{1,p(x)}(\Omega)$. Finally $w>0$ a.e. on $\Omega$ implies that $\theta w>0$ a.e. on $\Omega$ for any $\theta>0$. In conclusion $\theta w\in\overset{\bullet}{W}$.
	
	\noindent\textbf{(2)} A first remark is that from $w_1,w_2>0$ a.e. on $\Omega$ we easily get that $w:=w_1+w_2>0$ a.e. on $\Omega$. From $w_1,w_2\in L^{\frac{p(x)}{\alpha}}(\Omega)$ we also get that $w\in L^{\frac{p(x)}{\alpha}}(\Omega)$. From this we immediately get that $\sqrt[\alpha]{w}\in L^{p(x)}(\Omega)$. It remains to show the existence of $\nabla\sqrt[\alpha]{w}$ and the fact that $\nabla\sqrt[\alpha]{w}\in L^{p(x)}(\Omega)^N$. 
	
	\noindent\textbf{(3)} This is a straightforward consequence of \textbf{(1)} and \textbf{(2)}.
	
\end{proof}

\begin{lemma}\label{lemmatrunc} Let $w\in\overset{\bullet}{W}$ and define for each $\epsilon\in (0,1)$ the truncated function $w_{\epsilon}:\Omega\to\mathbb{R},\ w_{\epsilon}(x)=\begin{cases} \epsilon, & w(x)\leq \epsilon\\w(x), & w(x)\in I_{\epsilon}:=\left (\epsilon,\dfrac{1}{\epsilon} \right )\\ \dfrac{1}{\epsilon}, & w(x)\geq \dfrac{1}{\epsilon}\end{cases}$. Then:
	
	\begin{enumerate}
		\item[\textnormal{\textbf{(1)}}] $\lim\limits_{\epsilon\searrow 0} w_{\epsilon}(x)=w(x)$ for a.a. $x\in\Omega$.
		
		\item[\textnormal{\textbf{(2)}}] $\nabla\sqrt[\alpha]{w_{\epsilon}}=\chi_{I_{\epsilon}}(w)\nabla\sqrt[\alpha]{w}$.
		
		\item[\textnormal{\textbf{(3)}}] $w_{\epsilon}\in W^{1,p(x)}(\Omega)$ and $\nabla w_{\epsilon}=\alpha\sqrt[\alpha]{w^{\alpha-1}}\nabla \sqrt[\alpha]{w} \chi_{I_{\epsilon}}(w)$. 
		
		\noindent In particular $\lim\limits_{\epsilon\searrow 0} \nabla w_{\epsilon}=\alpha\sqrt[\alpha]{w^{\alpha-1}}\nabla \sqrt[\alpha]{w}$ a.e. on $\Omega$.
		
		\item[\textnormal{\textbf{(4)}}] $w_{\epsilon}^r\in W^{1,p(x)}(\Omega)$ and $\nabla w_{\epsilon}^r=rw_{\epsilon}^{r-1}\nabla w_{\epsilon}$ for any $r\in\mathbb{R}$.
		
		\item[\textnormal{\textbf{(5)}}] $w_{\epsilon}\stackrel{\epsilon\searrow 0}{\longrightarrow} w$ in $L^1(\Omega)$.
		
		\item[\textnormal{\textbf{(6)}}] For any $r>0$ we have that $\lim\limits_{\epsilon\searrow 0} \displaystyle\int_{\Omega} |w_\epsilon-w|^r\ dx=0$.
		
		\item[\textnormal{\textbf{(7)}}] $w\in W^{1,1}(\Omega)$, i.e. $\overset{\bullet}{W}\subset W^{1,1}(\Omega)$. In particular $\nabla w=\alpha\sqrt[\alpha]{w^{\alpha-1}}\nabla \sqrt[\alpha]{w}$ a.e. on $\Omega$.
		
		\item[\textnormal{\textbf{(8)}}] $\lim\limits_{\epsilon\searrow 0}\nabla w_{\epsilon}(x)=\nabla w(x)$, for a.a. $x\in\Omega$.
		
		\item[\textnormal{\textbf{(9)}}] $|w_{\epsilon}-w|\leq |w-1|$ a.e. on $\Omega$.
	\end{enumerate}
	
\end{lemma}

\begin{proof} \textbf{(1)} For a.a. $x\in\Omega$ we have that $w(x)>0$. So for $\epsilon>0$ sufficiently small we have that $w(x)\in I_{\epsilon}$. Thus $\lim\limits_{\epsilon\searrow 0} w_{\epsilon}(x)=\lim\limits_{\epsilon\searrow 0} w(x)=w(x)$ as needed.

\noindent \textbf{(2)} The function $\mathbb{R}\ni s\mapsto\begin{cases}\epsilon^{\frac{1}{\alpha}}, & s\leq\epsilon^{\frac{1}{\alpha}}\\ s, & s\in I_{\epsilon^{\frac{1}{\alpha}}}\\ \left (\dfrac{1}{\epsilon}\right )^{\frac{1}{\alpha}}, & s\geq \left (\frac{1}{\epsilon} \right )^{\frac{1}{\alpha}} \end{cases}$ is 1--Lipschitz. Composing it with $w^{\frac{1}{\alpha}}$ gives:
	
	$\begin{cases}\epsilon^{\frac{1}{\alpha}}, & w^{\frac{1}{\alpha}}\leq\epsilon^{\frac{1}{\alpha}}\\ w^{\frac{1}{\alpha}}, & w^{\frac{1}{\alpha}}\in I_{\epsilon^{\frac{1}{\alpha}}}\\ \left (\dfrac{1}{\epsilon}\right )^{\frac{1}{\alpha}}, & w^{\frac{1}{\alpha}}\geq \left (\frac{1}{\epsilon} \right )^{\frac{1}{\alpha}} \end{cases}=\begin{cases}\epsilon^{\frac{1}{\alpha}}, & w\leq\epsilon\\ w^{\frac{1}{\alpha}}, & w\in I_{\epsilon}\\ \left (\dfrac{1}{\epsilon}\right )^{\frac{1}{\alpha}}, & w\geq \frac{1}{\epsilon} \end{cases}=w_{\epsilon}^{\frac{1}{\alpha}}$.

\noindent  So, from the \textit{chain rule}\footnote{See \cite[Theorem 11.1.2, page 345]{Hasto}.}, we get that $\sqrt[\alpha]{w_{\epsilon}}\in W^{1,p(x)}(\Omega)$ and:
	
\begin{equation}
	\nabla\sqrt[\alpha]{w_{\epsilon}}=\chi_{I_{\epsilon^{\frac{1}{\alpha}}}}(\sqrt[\alpha]{w})\nabla \sqrt[\alpha]{w}=\chi_{I_{\epsilon}}(w)\nabla \sqrt[\alpha]{w}.
\end{equation}
	
\noindent \textbf{(3)} Also, the function $\mathbb{R}\ni s\mapsto\begin{cases}\epsilon, & s\leq\epsilon^{\frac{1}{\alpha}}\\ s^{\alpha}, & s\in I_{\epsilon^{\frac{1}{\alpha}}}\\ \dfrac{1}{\epsilon}, & s\geq \left (\dfrac{1}{\epsilon}\right )^{\frac{1}{\alpha}} \end{cases}$ is $\alpha\left (\dfrac{1}{\epsilon} \right )^{\frac{\alpha-1}{\alpha}}$-- Lipschitz. When we compose it with $\sqrt[\alpha]{w_{\epsilon}}$ we get exactly $w_{\epsilon}$. So, again from the \textit{chain rule}, we obtain that $w_{\epsilon}\in W^{1,p(x)}(\Omega)$ and:

\begin{equation}\label{eqderweak}\nabla w_{\epsilon}=\alpha \sqrt[\alpha]{w_{\epsilon}^{\alpha-1}}\nabla \sqrt[\alpha]{w_{\epsilon}}\chi_{I_{\epsilon^{\frac{1}{\alpha}}}}(\sqrt[\alpha]{w_{\epsilon}})=\alpha \sqrt[\alpha]{w^{\alpha-1}}\nabla \sqrt[\alpha]{w}\chi_{I_{\epsilon}}(w).
\end{equation}

\noindent\textbf{(4)} If $r>0$ take the function $\mathbb{R}\ni s\mapsto\begin{cases} \epsilon^r, & s\leq \epsilon^r\\ s^r, & s\in \left (\epsilon^r,\frac{1}{\epsilon^r}\right )\\ \left (\dfrac{1}{\epsilon}\right )^r, & s\geq\frac{1}{\epsilon^r}\end{cases}$ and if $r<0$ consider $\mathbb{R}\ni s\mapsto\begin{cases}\left (\dfrac{1}{\epsilon}\right )^r, & s\leq\frac{1}{\epsilon^r} \\ s^r, & s\in \left (\epsilon^r,\frac{1}{\epsilon^r}\right )\\ \epsilon^r, & s\geq \epsilon^r \end{cases}$. Both functions are $|r|\max\left \{\epsilon^{r-1},\frac{1}{\epsilon^{r-1}} \right \}$--Lipschitz. Composing this function with $w_{\epsilon}\in W^{1,p(x)}(\Omega)$ will give as result $w_{\epsilon}^r$. So, from the \textit{chain rule}, we deduce that $w_{\epsilon}^r\in W^{1,p(x)}(\Omega)$ and the formula given in the statement holds.

\noindent\textbf{(5)} and \textbf{(6)} Take any decreasing sequence $(\epsilon_n)_{n\geq 1}$ that converges to $0$. Notice that:

\begin{align*}
	\int_{\Omega} |w_{\epsilon_n}-w|^r\ dx&=\int_{[w\leq\epsilon_n]} (\epsilon_n-w)^r\ dx+\int_{\big [w\geq\frac{1}{\epsilon_n}\big ]} \left (w-\dfrac{1}{\epsilon_n}\right )^r\\
	&=\int_{\Omega} (\epsilon_n-w)^r\chi_{[w\leq\epsilon_n]}(x)\ dx+\int_{\Omega}\left( w-\dfrac{1}{\epsilon_n}\right)^r \chi_{\big [w\geq\frac{1}{\epsilon_n}\big ]}(x)\ dx.
\end{align*}

\noindent If we set $g_n=(\epsilon_n-w)^r\chi_{[w\leq\epsilon_n]}$ and $h_n=\left( w-\dfrac{1}{\epsilon_n}\right)^r \chi_{\big [w\geq\frac{1}{\epsilon_n}\big ]}$ we can easily remark that: $g_{n}(x)\geq g_{n+1}(x)$ and $h_n(x)\geq h_{n+1}(x)$ for a.a. $x\in\Omega$ and any $n\geq 1$. Moreover: $\lim\limits_{n\to\infty} g_n(x)=\lim\limits_{n\to\infty} h_n(x)=0$ for a.a. $x\in\Omega$ because $w(x)>0$ for a.a. $x\in\Omega$. Then from \textit{Beppo-Levi monotone convergence theorem} we obtain that $g_n,h_n\longrightarrow 0$ in $L^1(\Omega)$. In particular, for $r=1$ this shows that $w_{\epsilon_n}\to w$ in $L^1(\Omega)$ as required.

\noindent\textbf{(7)} \noindent From \eqref{eqderweak} we can write for any $\phi\in C^{\infty}_c(\Omega)$ and any $i\in\overline{1,N}$ that:

\begin{equation}\label{eqtolimit}
	\int_{\Omega} w_{\epsilon}\dfrac{\partial\phi}{\partial x_i}\ dx=-\int_{\Omega} \alpha w^{\frac{\alpha-1}{\alpha}}\dfrac{\partial\sqrt[\alpha]{w}}{\partial x_i}\chi_{I_{\epsilon}}(w)\phi\ dx.
\end{equation}

\noindent Note that:

\begin{align*}
	\int_{\Omega} \left |w_{\epsilon}\dfrac{\partial\phi}{\partial x_i}-w\dfrac{\partial\phi}{\partial x_i} \right |\ dx\leq \left\Vert\dfrac{\partial\phi}{\partial x_i}\right \Vert_{L^{\infty}(\Omega)}\int_{\Omega} |w_{\epsilon}-w|\ dx\ \stackrel{\epsilon\searrow 0}{\longrightarrow} 0.
\end{align*}

\noindent Here I have used that $w_{\epsilon}\to w$ in $L^1(\Omega)$. For the right-hand side we have that $w\in L^{\frac{p(x)}{\alpha}}(\Omega)\ \Rightarrow\ w^{\frac{\alpha-1}{\alpha}}\in L^{\frac{p(x)}{\alpha-1}}(\Omega)\hookrightarrow L^{p'(x)}(\Omega)$, because $\dfrac{p(x)}{\alpha-1}\geq\dfrac{p(x)}{p(x)-1}=p'(x)\ \Leftrightarrow\ \alpha\leq p(x)$ which is true since $\alpha\leq p^-\leq p(x)$ (from \textbf{(H13)} or \textbf{(H13')}) for all $x\in\overline{\Omega}$. Note also that $\dfrac{\partial \sqrt[\alpha]{w}}{\partial x_i}\in L^{p(x)}(\Omega)$ (because $\sqrt[\alpha]{w}\in W^{1,p(x)}(\Omega)$), so from \textit{H\"{o}lder inequality} we obtain that $\alpha w^{\frac{\alpha-1}{\alpha}}\dfrac{\partial \sqrt[\alpha]{w}}{\partial x_i}\in L^{1}(\Omega)$. Since:

\begin{enumerate}
	\item[$\bullet$] $\lim\limits_{\epsilon\searrow 0} \alpha w^{\frac{\alpha-1}{\alpha}}\dfrac{\partial \sqrt[\alpha]{w}}{\partial x_i}\chi_{I_{\epsilon}}(w)\phi=\alpha w^{\frac{\alpha-1}{\alpha}}\dfrac{\partial \sqrt[\alpha]{w}}{\partial x_i}\phi$, and
	
	\item[$\bullet$] $\left |\alpha w^{\frac{\alpha-1}{\alpha}}\dfrac{\partial \sqrt[\alpha]{w}}{\partial x_i}\chi_{I_{\epsilon}}(w)\phi \right |\leq \alpha w^{\frac{\alpha-1}{\alpha}}\dfrac{\partial \sqrt[\alpha]{w}}{\partial x_i}\phi\in L^1(\Omega)$ for every $\epsilon\in (0,1)$,
\end{enumerate}

\noindent we deduce from \textit{Lebesgue dominated convergence theorem} that $\displaystyle\int_{\Omega} \alpha w^{\frac{\alpha-1}{\alpha}}\dfrac{\partial\sqrt[\alpha]{w}}{\partial x_i}\chi_{I_{\epsilon}}(w)\phi\ dx\stackrel{\epsilon\searrow 0}{\longrightarrow} \displaystyle\int_{\Omega} \alpha w^{\frac{\alpha-1}{\alpha}}\dfrac{\partial\sqrt[\alpha]{w}}{\partial x_i}\phi\ dx$. 

\noindent So, making $\epsilon\searrow 0$ in \eqref{eqtolimit}, we'll get that:

\begin{equation}
	\int_{\Omega} w\dfrac{\partial\phi}{\partial x_i}\ dx=-\int_{\Omega} \alpha w^{\frac{\alpha-1}{\alpha}}\dfrac{\partial\sqrt[\alpha]{w}}{\partial x_i}\phi\ dx,
\end{equation}

\noindent for each $i\in\overline{1,N}$ and for any $\phi\in C^{\infty}_c(\Omega)$. This shows that the weak gradient $\nabla w$ exists and $\nabla w=\alpha w^{\frac{\alpha-1}{\alpha}}\nabla\sqrt[\alpha]{w}\in L^1(\Omega)$. We conclude that $w\in W^{1,1}(\Omega)$.

\noindent\textbf{(8)} From \textbf{(3)} we already know that $\lim\limits_{\epsilon\searrow 0} \nabla w_{\epsilon}(x)=\alpha\sqrt[\alpha]{w(x)^{\alpha-1}}\nabla\sqrt[\alpha]{w(x)}=\nabla w(x)$, from the proof of \textbf{(7)}.

\bigskip

\noindent\textbf{(9)} $|w_{\epsilon}-w|=\begin{cases}\epsilon-w\leq 1-w=|w-1|, & w\leq \epsilon\\[3mm] 0\leq |w-1|, & w\in I_{\epsilon}\\[3mm] w-\frac{1}{\epsilon}\leq w-1=|w-1|, & w\geq\frac{1}{\epsilon} \end{cases}\leq |w-1|$, a.e. on $\Omega$.

\end{proof}

\begin{lemma}\label{lemsaaintro} For any $w_1,w_2\in W^{1,p(x)}(\Omega)$ with $w_1,w_2>0$ a.e. on $\Omega$ and $\dfrac{w_2}{w_1},\dfrac{w_1}{w_2}\in L^{\infty}(\Omega)$ we have that $\dfrac{w_2^r}{w_1^{r-1}},\ \dfrac{w_1^r}{w_2^{r-1}}\in W^{1,p(x)}(\Omega)$ for any $r\in\mathbb{R}$ and moreover:
	
	\begin{align}
		&\nabla\left (\dfrac{w_2^r}{w_1^{r-1}} \right )=r\left (\dfrac{w_2}{w_1} \right )^{r-1}\nabla w_2-(r-1)\left (\dfrac{w_2}{w_1}\right )^r\nabla w_1\\
		&\nabla\left (\dfrac{w_1^r}{w_2^{r-1}} \right )=r\left (\dfrac{w_1}{w_2} \right )^{r-1}\nabla w_1-(r-1)\left (\dfrac{w_1}{w_2}\right )^r\nabla w_2.
	\end{align}
	
\end{lemma} 

\begin{proof} For each $0<\epsilon<1$ consider the following function:
	
	\begin{equation}
		\eta_{\epsilon}:\mathbb{R}\to (0,\infty),\ \eta_{\epsilon}(w)=\begin{cases} \epsilon, & w\leq\epsilon\\ w, & w\in I_{\epsilon}:=\left(\epsilon,\dfrac{1}{\epsilon}\right)\\ \dfrac{1}{\epsilon}, & w\geq\dfrac{1}{\epsilon}\end{cases}.
	\end{equation}
	
	\noindent It is easy to check that $\eta_{\epsilon}$ is a $1-$Lipschitz function, i.e. $|\eta_{\epsilon}(\tilde{w})-\eta_{\epsilon}(w)|\leq |\tilde{w}-w|$ for any $\tilde{w},w\in (0,\infty)$.

	\noindent Therefore, using the chain rule for weak derivatives\footnote{See Exercise 11.51 (i) from \cite[Page 340]{leonibook}.} which works fine even for variable exponents as pointed in \cite[Theorem 11.1.2, page 345]{Hasto}, we deduce that\footnote{Here we have used that $\Omega$ is bounded so that there is no need to have $\eta_{\epsilon}(0)=0$.} $w_{1,\epsilon}:=\eta_{\epsilon}\circ w_1,\ w_{2,\epsilon}:=\eta_{\epsilon}\circ w_2\in W^{1,p(x)}(\Omega)\cap L^{\infty}(\Omega)$ and moreover for a.a. $x\in\Omega$:\footnote{Since $w_1,w_2\in W^{1,p(x)}(\Omega)\hookrightarrow W^{1,p^-}(\Omega)$ and $\eta_\epsilon$ is a Lipschitz function we get that $w_{1,\epsilon},w_{2,\epsilon}\in W^{1,p^-}(\Omega)$ and the give formulas hold. Noticing that $|\nabla w_{1,\epsilon}|\leq |\nabla w_1|$ and $|\nabla w_{2,\epsilon}|\leq |\nabla w_2|$ we get that $|\nabla w_{1,\epsilon}|,|\nabla w_{2,\epsilon}|\in L^{p(x)}(\Omega)$. Similarly it is easy to observe that $w_{1,\epsilon},w_{2,\epsilon}\in L^{p(x)}(\Omega)$ as $w_1,w_2\in L^{p(x)}(\Omega)$ and constant functions are in $L^{p(x)}(\Omega)$. Combining all these facts we get that $w_{1,\epsilon},w_{2,\epsilon}\in W^{1,p(x)}(\Omega)$ and the given chai rule formulas hold.}
	
	\begin{equation}\label{apeqdiazimp1}
		\begin{cases} \nabla w_{1,\epsilon}=\nabla \eta_{\epsilon}\circ w_1=\chi_{I_{\epsilon}}(w_1)\nabla w_1\\ \nabla w_{2,\epsilon} \nabla \eta_{\epsilon}\circ w_2=\chi_{I_{\epsilon}}(w_2)\nabla w_2\end{cases}\ \Longrightarrow\ \begin{cases} |\nabla w_{1,\epsilon}|\leq|\nabla w_1|\\ |\nabla w_{2,\epsilon}|\leq |\nabla w_2|\end{cases},\ \forall\ \epsilon>0.
	\end{equation}

	\noindent The most important fact here is to notice that:
	
	\begin{equation}\label{apeqdiazimp2}
		\begin{cases}\dfrac{w_{1,\epsilon}}{w_{2,\epsilon}}\leq\max\left\{\dfrac{w_1}{w_2},1\right\}, \text{a.e. on}\ \Omega\\[3mm] \dfrac{w_{2,\epsilon}}{w_{1,\epsilon}}\leq\max\left\{\dfrac{w_2}{w_1},1\right\}, \text{a.e. on}\ \Omega\end{cases}, \forall\ \epsilon>0.
	\end{equation}
	
	\noindent We have pointwise that for a.a. $x\in\Omega$:
	
	\begin{equation}
		\lim\limits_{\epsilon\searrow 0} \dfrac{w_{2,\epsilon}(x)^r}{w_{1,\epsilon}(x)^{r-1}}=\dfrac{w_2(x)^r}{w_1(x)^{r-1}},\quad\ \lim\limits_{\epsilon\searrow 0} \dfrac{w_{1,\epsilon}(x)^r}{w_{2,\epsilon}(x)^{r-1}}=\dfrac{w_1(x)^r}{w_2(x)^{r-1}},
	\end{equation}
	
	\noindent and:
	
	\begin{equation}\label{apeqpointdiaz1}
		\lim\limits_{\epsilon\searrow 0} r\left(\dfrac{w_{2,\epsilon}(x)}{w_{1,\epsilon}(x)}\right )^{r-1}\nabla w_{2,\epsilon}(x)-(r-1)\left (\dfrac{w_{2,\epsilon}(x)}{w_{1,\epsilon}(x)}\right )^{r}\nabla w_{1,\epsilon}(x)=r\left (\dfrac{w_{2}(x)}{w_{1}(x)}\right )^{r-1}\nabla w_{2}(x)-(r-1)\left (\dfrac{w_{2}(x)}{w_{1}(x)}\right )^r\nabla w_{1}(x).
	\end{equation}
	
	\begin{equation}\label{apeqpointdiaz2}
		\lim\limits_{\epsilon\searrow 0} r\left(\dfrac{w_{1,\epsilon}(x)}{w_{2,\epsilon}(x)}\right )^{r-1}\nabla w_{1,\epsilon}(x)-(r-1)\left (\dfrac{w_{1,\epsilon}(x)}{w_{2,\epsilon}(x)}\right )^{r}\nabla w_{2,\epsilon}(x)=r\left (\dfrac{w_{1}(x)}{w_{2}(x)}\right )^{r-1}\nabla w_{1}(x)-(r-1)\left (\dfrac{w_{1}(x)}{w_{2}(x)}\right )^r\nabla w_{2}(x).
	\end{equation}
	
	\noindent Now we argue that the following formulas hold:
	
	\begin{align}
		&\nabla\left (\dfrac{w_{2,\epsilon}^r}{w_{1,\epsilon}^{r-1}} \right )=r\left (\dfrac{w_{2,\epsilon}}{w_{1,\epsilon}} \right )^{r-1}\nabla w_{2,\epsilon}-(r-1)\left (\dfrac{w_{2,\epsilon}}{w_{1,\epsilon}}\right )^r\nabla w_{1,\epsilon}\\
		&\nabla\left (\dfrac{w_{1,\epsilon}^r}{w_{2,\epsilon}^{r-1}} \right )=r\left (\dfrac{w_{1,\epsilon}}{w_{2,\epsilon}} \right )^{r-1}\nabla w_{1,\epsilon}-(r-1)\left (\dfrac{w_{1,\epsilon}}{w_{2,\epsilon}}\right )^r\nabla w_{2,\epsilon}.
	\end{align}
	
	\noindent Indeed since $w_{1,\epsilon},w_{2,\epsilon}\in W^{1,p(x)}(\Omega)\hookrightarrow W^{1,p^-}(\Omega)$ and $w_{1,\epsilon},w_{2,\epsilon}\in L^{\infty}(\Omega)$ we deduce that the \textit{chain rule} and \textit{product rule} formulas work. We will explain the first formula. For the second one we can proceed in the same way.
	
	\noindent Fix some arbitrary $i\in\overline{1,N}$. The function $\mathbb{R}\ni s\longmapsto\begin{cases} \dfrac{1}{\epsilon}, & s\leq\epsilon\\[3mm] \dfrac{1}{s}, & s\in \left (\epsilon,\dfrac{1}{\epsilon}\right )\\[3mm]\epsilon, & s\geq \dfrac{1}{\epsilon}\end{cases}$ is $\dfrac{1}{\epsilon^2}-$Lipschitz. Therefore from the \textit{chain rule for weak derivatives} we have that $\dfrac{1}{w_{1,\varepsilon}}\in W^{1,p(x)}(\Omega)$ and:
	
	\begin{equation}\label{apeqdiaz10} \dfrac{\partial}{\partial x_i}\left(\dfrac{1}{w_{1,\epsilon}}\right)=-\dfrac{1}{w_{1,\epsilon}^2}\chi_{I_{\epsilon}}(w_{1,\epsilon})\dfrac{\partial w_{1,\epsilon}}{\partial x_i}=-\dfrac{1}{w_{1,\epsilon}^2}\dfrac{\partial w_{1,\epsilon}}{\partial x_i},
	\end{equation} 
	
	\noindent since $\dfrac{1}{\epsilon}\geq w_{1,\epsilon}\geq \epsilon$. This also shows that $\dfrac{1}{w_{1,\epsilon}}\in L^{\infty}(\Omega)$.
	
	\noindent Now, since $w_{2,\epsilon},\dfrac{1}{w_{1,\epsilon}}\in W^{1,p(x)}(\Omega)\cap L^{\infty}(\Omega)$, we get from the \textit{product rule for weak derivatives}\footnote{See \cite[Exercise 11.51 (ii), page 340]{leonibook}.} that $\dfrac{w_{2,\epsilon}}{w_{1,\epsilon}}\in W^{1,p(x)}(\Omega)\cap L^{\infty}(\Omega)$ and using \eqref{apeqdiaz10}:
	
	\begin{equation}\label{apeqdiaz20}
		\dfrac{\partial}{\partial x_i}\left(\dfrac{w_{2,\epsilon}}{w_{1,\epsilon}}\right)=\dfrac{1}{w_{1,\epsilon}}\dfrac{\partial w_{2,\epsilon}}{\partial x_i}+w_{2,\epsilon}\dfrac{\partial}{\partial x_i}\left(\dfrac{1}{w_{1,\epsilon}}\right)=\dfrac{1}{w_{1,\epsilon}}\dfrac{\partial w_{2,\epsilon}}{\partial x_i}-\dfrac{w_{2,\epsilon}}{w_{1,\epsilon}^2}\dfrac{\partial w_{1,\epsilon}}{\partial x_i}.
	\end{equation}
	
	\noindent The function $\mathbb{R}\ni s\longmapsto\begin{cases}\epsilon^{2r-2}, & s\leq \epsilon^2\\ s^{r-1}, & s\in \left(\epsilon^2,\dfrac{1}{\epsilon^2}\right)\\[3mm] \dfrac{1}{\epsilon^{2r-2}}, & s\geq \dfrac{1}{\epsilon^2}\end{cases}$ is $(r-1)\max\{\epsilon^{2r-2},\epsilon^{-(2r-2)}\}-$Lipschitz. Therefore using one more time the \textit{chain rule for weak derivatives} we obtain that $\left(\dfrac{w_{2,\epsilon}}{w_{1,\epsilon}}\right)^{r-1}\in W^{1,p(x)}(\Omega)$ and with the help of \eqref{apeqdiaz20}:
	
	\begin{align}\label{apeqdiaz30}
		\dfrac{\partial}{\partial x_i}\left(\dfrac{w_{2,\epsilon}}{w_{1,\epsilon}}\right)^{r-1}&=(r-1)\left(\dfrac{w_{2,\epsilon}}{w_{1,\epsilon}}\right)^{r-2}\chi_{I_{\epsilon^2}}\left (\dfrac{w_{2,\epsilon}}{w_{1,\epsilon}} \right)	\dfrac{\partial}{\partial x_i}\left(\dfrac{w_{2,\epsilon}}{w_{1,\epsilon}}\right)=(r-1)\left(\dfrac{w_{2,\epsilon}}{w_{1,\epsilon}}\right)^{r-2}\dfrac{\partial}{\partial x_i}\left(\dfrac{w_{2,\epsilon}}{w_{1,\epsilon}}\right) \nonumber \\ 
		&=(r-1)\left(\dfrac{w_{2,\epsilon}}{w_{1,\epsilon}}\right)^{r-2}\left(\dfrac{1}{w_{1,\epsilon}}\dfrac{\partial w_{2,\epsilon}}{\partial x_i}-\dfrac{w_{2,\epsilon}}{w_{1,\epsilon}^2}\dfrac{\partial w_{1,\epsilon}}{\partial x_i} \right),
	\end{align}
	
	\noindent as $\dfrac{1}{\epsilon^2}\geq \dfrac{w_{2,\epsilon}}{w_{1,\epsilon}}\geq\epsilon^2$. From here we deduce one more time that $\left (\dfrac{w_{2,\epsilon}}{w_{1,\epsilon}}\right )^{r-1}\in L^{\infty}(\Omega)$.

	\noindent At last, because $w_{2,\epsilon},\ \left (\dfrac{w_{2,\epsilon}}{w_{1,\epsilon}}\right )^{r-1}\in W^{1,p(x)}(\Omega)\cap L^{\infty}(\Omega)$, we deduce from the \textit{product rule for the weak derivatives} that $\dfrac{w_{2,\epsilon}^r}{w_{1,\epsilon}^{r-1}}=w_{2,\epsilon}\left (\dfrac{w_{2,\epsilon}}{w_{1,\epsilon}}\right )^{r-1}\in W^{1,p(x)}(\Omega)\cap L^{\infty}(\Omega)$ and using \eqref{apeqdiaz30}:
	
	\begin{align*}
		\dfrac{\partial}{\partial x_i}\left(\dfrac{w_{2,\epsilon}^r}{w_{1,\epsilon}^{r-1}}\right)&=\left (\dfrac{w_{2,\epsilon}}{w_{1,\epsilon}}\right )^{r-1}\dfrac{\partial w_{2,\epsilon}}{\partial x_i}+w_{2,\epsilon}	\dfrac{\partial}{\partial x_i}\left(\dfrac{w_{2,\epsilon}}{w_{1,\epsilon}}\right )^{r-1}\\
		&= \left (\dfrac{w_{2,\epsilon}}{w_{1,\epsilon}}\right )^{r-1}\dfrac{\partial w_{2,\epsilon}}{\partial x_i}+(r-1)w_{2,\epsilon}\left(\dfrac{w_{2,\epsilon}}{w_{1,\epsilon}}\right)^{r-2}\left(\dfrac{1}{w_{1,\epsilon}}\dfrac{\partial w_{2,\epsilon}}{\partial x_i}-\dfrac{w_{2,\epsilon}}{w_{1,\epsilon}^2}\dfrac{\partial w_{1,\epsilon}}{\partial x_i} \right)\\
		&=r \left (\dfrac{w_{2,\epsilon}}{w_{1,\epsilon}}\right )^{r-1}\dfrac{\partial w_{2,\epsilon}}{\partial x_i}-(r-1)\left(\dfrac{w_{2,\epsilon}}{w_{1,\epsilon}}\right)^{r}\dfrac{\partial w_{1,\epsilon}}{\partial x_i}
	\end{align*}
	
	\noindent Writing this relation for any $i\in\overline{1,N}$ we may write: $\nabla\left (\dfrac{w_{2,\epsilon}^r}{w_{1,\epsilon}^{r-1}} \right )=r\left (\dfrac{w_{2,\epsilon}}{w_{1,\epsilon}} \right )^{r-1}\nabla w_{2,\epsilon}-(r-1)\left (\dfrac{w_{2,\epsilon}}{w_{1,\epsilon}}\right )^r\nabla w_{1,\epsilon}$.
	
	\noindent We are allowed now, from the definition of the weak derivative, to write for all $\phi\in C^{\infty}_c(\Omega)$ that:
	
	\begin{equation}\label{apeqdiazlim}
		\int_{\Omega} \dfrac{w_{2,\epsilon}^r}{w_{1,\epsilon}^{r-1}}\nabla\phi(x)\ dx=-\int_{\Omega} \left [r\left (\dfrac{w_{2,\epsilon}}{w_{1,\epsilon}} \right )^{r-1}\nabla w_{2,\epsilon}-(r-1)\left (\dfrac{w_{2,\epsilon}}{w_{1,\epsilon}}\right )^r\nabla w_{1,\epsilon} \right ]\phi(x)\ dx.
	\end{equation}
	
	\noindent Next, we will show that we can pass this relation to the limit as $\epsilon\searrow 0$ and obtain the formulas from the statement of this lemma.
	
	\noindent We begin with the right hand side: from \eqref{apeqpointdiaz1} and \eqref{apeqpointdiaz2} we have that:
	
	\begin{equation}
		r\left (\dfrac{w_{2,\epsilon}}{w_{1,\epsilon}} \right )^{r-1}\nabla w_{2,\epsilon}-(r-1)\left (\dfrac{w_{2,\epsilon}}{w_{1,\epsilon}}\right )^r\nabla w_{1,\epsilon}\longrightarrow r\left (\dfrac{w_{2}}{w_{1}} \right )^{r-1}\nabla w_{2}-(r-1)\left (\dfrac{w_{2}}{w_{1}}\right )^r\nabla w_{1}\ \text{pointwise a.e. on}\ \Omega.
	\end{equation}
	
	\noindent On the other hand, from \eqref{apeqdiazimp1} and \eqref{apeqdiazimp2} we have a.e. on $\Omega$ that:
	
	\begin{align}
		&\left |	r\left (\dfrac{w_{2,\epsilon}}{w_{1,\epsilon}} \right )^{r-1}\nabla w_{2,\epsilon}-(r-1)\left (\dfrac{w_{2,\epsilon}}{w_{1,\epsilon}}\right )^r\nabla w_{1,\epsilon} \right |\leq\\
		\leq & r\max\left\{\left \Vert\dfrac{w_2}{w_1}\right \Vert_{L^{\infty}(\Omega)}^{r-1},1\right\}|\nabla w_2|+(r-1)\max\left\{\left \Vert\dfrac{w_2}{w_1}\right \Vert_{L^{\infty}(\Omega)}^{r},1\right\}|\nabla w_1|\in L^{p(x)}(\Omega)\hookrightarrow L^1(\Omega).
	\end{align}
	
	\noindent Applying \textit{Lebesgue dominated convergence theorem} we obtain that:
	
	\begin{equation}
		r\left (\dfrac{w_{2,\epsilon}}{w_{1,\epsilon}} \right )^{r-1}\nabla w_{2,\epsilon}-(r-1)\left (\dfrac{w_{2,\epsilon}}{w_{1,\epsilon}}\right )^r\nabla w_{1,\epsilon}\stackrel{\epsilon\to 0}{\longrightarrow} r\left (\dfrac{w_{2}}{w_{1}} \right )^{r-1}\nabla w_{2}-(r-1)\left (\dfrac{w_{2}}{w_{1}}\right )^r\nabla w_{1}\ \text{in}\ L^1(\Omega).
	\end{equation}

	\noindent Taking now into account that $\phi\in C^{\infty}_c(\Omega)\subset L^{\infty}(\Omega)$ we deduce that:
	
	\begin{equation}\label{apeqdiazgata}
		\int_{\Omega} \left [r\left (\dfrac{w_{2,\epsilon}}{w_{1,\epsilon}} \right )^{r-1}\nabla w_{2,\epsilon}-(r-1)\left (\dfrac{w_{2,\epsilon}}{w_{1,\epsilon}}\right )^r\nabla w_{1,\epsilon} \right ]\phi\ dx\stackrel{\epsilon\to 0}{\longrightarrow}\int_{\Omega}\left [r\left (\dfrac{w_{2}}{w_{1}} \right )^{r-1}\nabla w_{2}-(r-1)\left (\dfrac{w_{2}}{w_{1}}\right )^r\nabla w_{1} \right ]\phi\ dx.
	\end{equation}
	
	\noindent For the left hand side of the equality \eqref{apeqdiazlim} we have that:
	
	\begin{align}\label{apeqdiazfinal}
		&\left |\int_{\Omega} \left (\dfrac{w_{2,\epsilon}^r}{w_{1,\epsilon}^{r-1}}- \dfrac{w_{2}^r}{w_{1}^{r-1}}\right )\nabla\phi(x) \ dx\right | \leq \Vert|\nabla\phi|\Vert_{L^{\infty}(\Omega)}\int_{\Omega} |w_{2,\epsilon}-w_{2}|\cdot \left (\dfrac{w_{2,\epsilon}}{w_{1,\epsilon}}\right )^{r-1}\ +w_2\left | \left (\dfrac{w_{2,\epsilon}}{w_{1,\epsilon}}\right )^{r-1}- \left (\dfrac{w_{2}}{w_{1}}\right )^{r-1} \right |\ dx\nonumber \\
		&\stackrel{\eqref{apeqdiazimp2}}{\leq} \Vert|\nabla\phi|\Vert_{L^{\infty}(\Omega)}\max\left\{\left \Vert\dfrac{w_2}{w_1}\right \Vert_{L^{\infty}(\Omega)}^{r-1},1\right\}\int_{\Omega} |w_{2,\epsilon}-w_{2}|\ dx+\int_{\Omega} w_2\left | \left (\dfrac{w_{2,\epsilon}}{w_{1,\epsilon}}\right )^{r-1}- \left (\dfrac{w_{2}}{w_{1}}\right )^{r-1} \right |\ dx.
	\end{align}
	
	\noindent Remark that: $\lim\limits_{\epsilon\searrow 0} w_2\left (\dfrac{w_{2,\epsilon}}{w_{1,\epsilon}}\right )^{r-1}- w_2\left (\dfrac{w_{2}}{w_{1}}\right )^{r-1}=0$ pointwise a.e. on $\Omega$, and from \eqref{apeqdiazimp2} we also get that $w_2\left | \left (\dfrac{w_{2,\epsilon}}{w_{1,\epsilon}}\right )^{r-1}- \left (\dfrac{w_{2}}{w_{1}}\right )^{r-1} \right |\leq w_2\max\left\{\left\Vert\dfrac{w_{2}}{w_1}\right \Vert_{L^{\infty}(\Omega)}^{r-1},1\right\}\in L^{p(x)}(\Omega)\hookrightarrow L^{1}(\Omega)$.
	
	\noindent Thus from the \textit{Lebesgue dominated convergence theorem} we obtain that: 
	
	\begin{equation}
		\lim\limits_{\epsilon\searrow 0} \int_{\Omega} w_2\left | \left (\dfrac{w_{2,\epsilon}}{w_{1,\epsilon}}\right )^{r-1}- \left (\dfrac{w_{2}}{w_{1}}\right )^{r-1} \right |\ dx=0.
	\end{equation}
	
	\noindent It remains to show that $\lim\limits_{\epsilon\searrow 0} w_{2,\epsilon}=w_2$ in $L^1(\Omega)$, but this follows from Lemma \ref{lemmatrunc} \textbf{(4)}. Hence making $\epsilon\to 0$ in \eqref{apeqdiazfinal} we get that $\displaystyle\lim\limits_{\epsilon\to 0}\int_{\Omega}  \dfrac{w_{2,\epsilon}^r}{w_{1,\epsilon}^{r-1}}\nabla\phi(x)\ dx=\int_{\Omega} \dfrac{w_{2}^r}{w_{1}^{r-1}}\nabla\phi(x) \ dx$. Henceforth, using also \eqref{apeqdiazgata}, we arrive at:
	
	\begin{equation}
		\int_{\Omega} \dfrac{w_{2}^r}{w_{1}^{r-1}}\nabla\phi(x)\ dx=-\int_{\Omega} \left [r\left (\dfrac{w_{2}}{w_{1}} \right )^{r-1}\nabla w_{2}-(r-1)\left (\dfrac{w_{2}}{w_{1}}\right )^r\nabla w_{1} \right ]\phi(x)\ dx,\ \forall\ \phi\in C^{\infty}_c(\Omega).
	\end{equation}
	
	\noindent This shows that $\dfrac{w_{2}^r}{w_{1}^{r-1}}$ has a weak gradient and the formula from the statement holds. It is easy to see that $\dfrac{w_{2}^r}{w_{1}^{r-1}}=\underbrace{w_2}_{\in L^{p(x)}(\Omega)}\underbrace{\left (\dfrac{w_2}{w_1}\right )^{r-1}}_{\in L^{\infty}(\Omega)}\in L^{p(x)}(\Omega)$ and $\nabla\left ( \dfrac{w_{2}^r}{w_{1}^{r-1}}\right )=r\underbrace{\left (\dfrac{w_{2}}{w_{1}} \right )^{r-1}}_{\in L^{\infty}(\Omega)}\underbrace{\nabla w_{2}}_{\in L^{p(x)}(\Omega)^N}-(r-1)\underbrace{\left (\dfrac{w_{2}}{w_{1}}\right )^r}_{\in L^{\infty}(\Omega)}\underbrace{\nabla w_{1}}_{\in L^{p(x)}(\Omega)^N}\in L^{p(x)}(\Omega)^N$. In conclusion $\dfrac{w_{2}^r}{w_{1}^{r-1}}\in W^{1,p(x)}(\Omega)$ and similarly $\dfrac{w_{1}^r}{w_{2}^{r-1}}\in W^{1,p(x)}(\Omega)$.
	
	\noindent 
\end{proof}

\section{An elementary proof for a general version of D\'{i}az-Saa inequality}

\noindent We start this section with a simple inequality lemma based on which we will prove the general version of D\'{i}az-Saa inequality (Theorem \ref{thmdiazsaa}). In the following section we will detect all the situations in which equality can occur.

\begin{lemma}\label{lemineqdiazsaa} Let $\phi:[0,\infty)\to [0,\infty),\ \phi(0)=0$ and $r\geq 1$ such that the function $(0,\infty)\ni\theta\mapsto\dfrac{\phi(\theta)}{\theta^{r-1}}$ is increasing. Then for any $a,b>0$ and $c,d\geq 0$ the following inequality holds:
	
	\begin{equation} \left [1+(r-1)\left (\frac{a}{b} \right )^r \right ]\phi(c)c+\left [1+(r-1)\left (\frac{b}{a} \right )^r \right ]\phi(d)d\geq r\left (\dfrac{a}{b} \right )^{r-1}\phi(c)d+r\left (\frac{b}{a} \right )^{r-1}\phi(d)c.
	\end{equation}
	
\noindent Suppose that equality holds. The following statements are true:

\begin{enumerate}
	\item[(i)] If $c=0$ then $d=0$ and equality holds.
	
	\item[(i)] If $r=1$, then $c=d$ or the function $\phi$ is constant between $c$ and $d$, and equality holds.
	
	\item[(ii)] If $r>1$ then $ac=bd$ and if moreover one of the following two requirements are fulfilled, equality holds:
	
	\begin{enumerate}
		\item[$\bullet$] $c=d$
		
		\item[$\bullet$] The function $\theta\mapsto\dfrac{\phi(\theta)}{\theta^{r-1}}$ is constant between $c$ and $d$.
	\end{enumerate}
	
\end{enumerate}

\noindent There are no other situations in which equality can occur.

\noindent In case we have that the function $[0,\infty)\ni\theta\mapsto\dfrac{\phi(\theta)}{\theta^{r-1}}$ is \textbf{strictly increasing} and if equality holds, the following statements are true:

\begin{enumerate}
	\item[(I)] If $c=0$ then $d=0$ and equality holds.
	
	\item[(II)] If $r=1$ then $c=d$ and equality holds.
	
	\item[(III)] If $r>1$ then $c=d$, $a=b$ and equality holds.
\end{enumerate}

\noindent There are no other situations in which equality can occur.
\end{lemma}

\begin{proof} If $c=0$ or $d=0$ there is nothing to prove. Also, if $r=1$, the inequality becomes: $\phi(c)c+\phi(d)d\geq \phi(c)d+\phi(d)c\ \Leftrightarrow\ \big (\phi(c)-\phi(d) \big )(c-d)\geq 0$ which is true, since $\phi$ is an increasing function.
	
\noindent Assume now that $c,d\neq 0$ and $r>1$. First we write:
	
	\begin{equation}
		r\left(\dfrac{a}{b} \right)^{r-1}\phi(c)d+r\left(\frac{b}{a} \right)^{r-1}\phi(d)c = r \left(\dfrac{ac}{b} \right)^{r-1}d\dfrac{\phi(c)}{c^{r-1}}+r\left(\frac{bd}{a} \right)^{r-1}c\dfrac{\phi(d)}{d^{r-1}},
	\end{equation}
	
	\noindent then we use Young's inequality with $r>1$ and its H\"{o}lder conjugate $\dfrac{r}{r-1}>0$:
	
	\begin{align}
		\left (\dfrac{ac}{b} \right )^{r-1}d &\leq \frac{r-1}{r}\left (\dfrac{ac}{b} \right )^{r} + \frac{d^r}{r},\\
		\left (\dfrac{bd}{a} \right )^{r-1}c &\leq \frac{r-1}{r}\left (\dfrac{bd}{a} \right )^{r} + \frac{c^r}{r}.
	\end{align}
	
	\noindent We should mention here that the equality holds for $r>1$ only when $ac=bd$. This gives us:
	
	\begin{equation}
		r \left(\dfrac{ac}{b} \right)^{r-1}d\dfrac{\phi(c)}{c^{r-1}}+r\left(\frac{bd}{a} \right)^{r-1}c\dfrac{\phi(d)}{d^{r-1}} \leq (r-1)\left(\dfrac{a}{b} \right)^{r}\phi(c)c + \dfrac{\phi(c)}{c^{r-1}}d^r + (r-1)\left(\dfrac{b}{a} \right)^{r}\phi(d)d + \dfrac{\phi(d)}{d^{r-1}}c^r,
	\end{equation}
	
	\noindent so all that remains to show is that: $\dfrac{\phi(c)}{c^{r-1}}d^r+\dfrac{\phi(d)}{d^{r-1}}c^r \leq \phi(c)c +\phi(d)d$ which is easy to see since the difference factors nicely:
	\begin{equation}
\phi(c)c+\phi(d)d-\dfrac{\phi(c)}{c^{r-1}}d^r-\dfrac{\phi(d)}{d^{r-1}}c^r=\left (\dfrac{\phi(c)}{c^{r-1}}-\dfrac{\phi(d)}{d^{r-1}} \right )(c^r-d^r)\geq 0.
	\end{equation}
	
	\noindent Here we have used that $(0,\infty)\ni \theta\mapsto\dfrac{\phi(\theta)}{\theta^{r-1}}$ is an increasing function. It is straightforward to check the equality cases. The proof is now completed.
	
\end{proof}

\begin{theorem}[\textbf{D\'{i}az-Saa inequality}]\label{thmdiazsaa}
	Let $w_1,w_2\in W^{1,p(x)}(\Omega)$ with $w_1,w_2>0$ a.e. on $\Omega$ and $\dfrac{w_2}{w_1},\dfrac{w_1}{w_2}\in L^{\infty}(\Omega)$. Then the following inequality holds:
	
	\begin{equation}\label{apeqdiaz}
		\int_{\Omega} \mathbf{a}(x,\nabla w_1)\cdot\nabla\left (w_1-\dfrac{w_2^r}{w_1^{r-1}}\right )\ dx\geq \int_{\Omega} \mathbf{a}(x,\nabla w_2)\cdot\nabla\left (\dfrac{w_1^r}{w_2^{r-1}}-w_2\right )\ dx.
	\end{equation}

\end{theorem}

\begin{proof} From Lemma \ref{lemsaaintro} we have that the expressions that appear in \eqref{apeqdiaz} are well-defined and finite. The inequality becomes:
	
	\begin{align*}  \int_{\Omega} \Phi(x,|\nabla w_1|)|\nabla w_1|&- \Psi(x,|\nabla w_1|)\nabla w_1\cdot\nabla\left (\dfrac{w_2^r}{w_1^{r-1}} \right )\ dx\geq \int_{\Omega}   \Psi(x,|\nabla w_2|)\nabla w_2\nabla\left (\dfrac{w_1^r}{w_2^{r-1}}\right )-\Phi(x,|\nabla w_2|)|\nabla w_2|\ dx \\
		\Longleftrightarrow \int_{\Omega} &\left [1+(r-1)\left(\dfrac{w_2}{w_1}\right)^r\right ]\Phi(x,|\nabla w_1|)|\nabla w_1|+\left [1+(r-1)\left(\dfrac{w_1}{w_2}\right)^r\right ]\Phi(x,|\nabla w_2|)|\nabla w_2|\ dx\geq\\
		\geq	& \int_{\Omega} r\left [\left(\dfrac{w_2}{w_1}\right)^{r-1}\Psi(x,|\nabla w_1|)+\left(\dfrac{w_1}{w_2}\right)^{r-1}\Psi(x,|\nabla w_2|)\right ]\nabla w_1\cdot\nabla w_2\ dx.\ \ \ \ \ \ \ \ \ \ 
	\end{align*}
	
	\noindent Fix some $x\in\Omega$. Since $[0,\infty)\ni s\mapsto\dfrac{\Phi(x,s)}{s^{r-1}}$ is increasing (see \textbf{(H7)}) from Lemma \ref{lemineqdiazsaa} applied for $a=w_2(x),\ b=w_1(x)>0$ and $c=|\nabla w_1(x)|,\ d=|\nabla w_2(x)|\geq 0$ we have that:
	
	\begin{align}\label{ineqeq}
		\left [1+(r-1)\left(\dfrac{w_2(x)}{w_1(x)}\right)^r\right ]&\Phi(x,|\nabla w_1|)|\nabla w_1(x)|+\left [1+(r-1)\left(\dfrac{w_1(x)}{w_2(x)}\right)^r\right ]\Phi(x,|\nabla w_2|)|\nabla w_2(x)|\geq \nonumber\\ 
		\geq	&\ r\left(\dfrac{w_2(x)}{w_1(x)}\right)^{r-1}\Phi(x,|\nabla w_1(x)|)|\nabla w_2(x)|+r\left(\dfrac{w_1(x)}{w_2(x)}\right)^{r-1}\Phi(x,|\nabla w_2(x)|)|\nabla w_1(x)|\nonumber\\
\text{(Cauchy ineq.)}\ \ \ \ \ \ \ 	\geq	&\ r\left [\left(\dfrac{w_2}{w_1}\right)^{r-1}\Psi(x,|\nabla w_1|)+\left(\dfrac{w_1}{w_2}\right)^{r-1}\Psi(x,|\nabla w_2|)\right ]\nabla w_1\cdot\nabla w_2.
	\end{align}

\noindent Integrating on $\Omega$ leads to the desired inequality.
\end{proof}

\begin{corollary}\label{cordiaz} In fact we have proved that pointwise:

\begin{equation}
	\mathbf{a}(x,\nabla w_1)\cdot\nabla\left (w_1-\dfrac{w_2^r}{w_1^{r-1}}\right )\geq \mathbf{a}(x,\nabla w_2)\cdot\nabla\left (\dfrac{w_1^r}{w_2^{r-1}}-w_2\right ),\ \text{a.e. in}\ \Omega.
\end{equation}

\noindent So we may integrate this inequality on any measurable subset $\tilde{\Omega}\subset\Omega$ to obtain that:

\begin{equation}\label{apeqdiaztilde}
	\int_{\tilde{\Omega}} \mathbf{a}(x,\nabla w_1)\cdot\nabla\left (w_1-\dfrac{w_2^r}{w_1^{r-1}}\right )\ dx\geq \int_{\tilde{\Omega}} \mathbf{a}(x,\nabla w_2)\cdot\nabla\left (\dfrac{w_1^r}{w_2^{r-1}}-w_2\right )\ dx.
\end{equation}
	
\end{corollary}

\begin{remark}
	Theorem \ref{thmdiazsaa} generalizes the form of the D\'{i}az-Saa inequality proved in \cite[Theorem 4.2]{giaco1} from two important points of view: 
	
	\begin{enumerate}
		\item[$\bullet$] The restrictive space $W_{0}^{p(x)}(\Omega)$ is replaced with the entire space $W^{1,p(x)}(\Omega)$ so that the inequality may be used also for homogeneous Neumann boundary conditions and not only homogeneous Dirichlet conditions.
		
		\item[$\bullet$] Hypothesis \textnormal{\textbf{(A7)}} from \cite[page 7]{giaco1} is a particular case in which hypothesis \textnormal{\textbf{(H7)}} from the present paper holds. For a better understanding of this fact we will provide in what follows a simple lemma.
	\end{enumerate}
\end{remark}

\begin{definition}
	We say that $A:\overline{\Omega}\times\mathbb{R}^N\to\mathbb{R}$ is $p(x)$--homogeneous if for a.a. $x\in\Omega$: 
	
	\begin{equation}
		A(x,t\xi)=|t|^{p(x)}A\left (x,\xi\right )\text{, for all}\ t\in\mathbb{R}\ \text{and}\ \xi\in\mathbb{R}^N.
	\end{equation} 
	
	\noindent This is hypothesis \textnormal{\textbf{(A7)}} from \cite[page 7]{giaco1} that can be also found in \cite[page 5]{tak}.
\end{definition}

\begin{lemma}\label{lemmasimplehomo} $A:\overline{\Omega}\times\mathbb{R}^N\to\mathbb{R}$ is $p(x)$--homogeneous $\Longleftrightarrow$ $\Phi:\overline{\Omega}\times\mathbb{R}\to\mathbb{R}$ is $p(x)-1$--homogeneous.
\end{lemma}

\begin{proof} If $A$ is $p(x)$--homogeneous then for some fixed $x\in\Omega$, any $t>0$ we may write for any $\xi\in\mathbb{R}^N$ that:
	
\begin{equation}
	t^{p(x)}\int_{0}^{|\xi|}\Phi(x,\tau)\ d\tau=t^{p(x)}A(x,\xi)=A(x,t\xi)=\int_{0}^{t|\xi|}\Phi(x,s)\ ds\stackrel{s=t\tau }{=}\int_{0}^{|\xi|}t\Phi(x,t\tau)\ d\tau.
\end{equation}

\noindent Dividing this relation by $t>0$ will give us that:

\begin{equation}
	\int_{0}^{|\xi|}\Phi(x,t\tau)-t^{p(x)-1}\Phi(x,\tau)\ d\tau=0,\ \forall\ \xi\in\mathbb{R}^N\ \Rightarrow\ \int_{0}^{s}\Phi(x,t\tau)-t^{p(x)-1}\Phi(x,\tau)\ d\tau=0,\ \forall\ s\in [0,\infty).
\end{equation}

\noindent Since the function $[0,s]\ni\tau\mapsto \Phi(x,t\tau)-t^{p(x)-1}\Phi(x,\tau)$ is continuous, from the \textit{Fundamental theorem of calculus}, we may differentiate with respect to $s$, and obtain: $\Phi(x,ts)-t^{p(x)-1}\Phi(x,s)=0$ for any $s\geq 0$. Since $\Phi(x,\cdot):\mathbb{R}\to\mathbb{R}$ is an even function we may also write that:

\begin{equation}
	\Phi(x,ts)=\Phi(x,|t|\cdot |s|)=|t|^{p(x)-1}\Phi(x,|s|)=|t|^{p(x)-1}\Phi(x,s),\ \forall\ t\in\mathbb{R}^*,\ \forall\ s\in\mathbb{R}.
\end{equation}
\noindent For $t=0$ the same relation holds, because $\Phi(x,0)=0$.

\noindent In the reverse sense, if $\Phi$ is $p(x)-1$--homogeneous then we easily obtain for any $t\in\mathbb{R}$ and $\xi\in\mathbb{R}^N$ that:

\begin{equation}
	A(x,t\xi)=\int_{0}^{|t|\cdot |\xi|}\Phi(x,s)\ ds\stackrel{s=|t|\tau}{=}\int_{0}^{|\xi|}\Phi(x,|t|\tau) |t|\ d\tau=|t|^{p(x)}\int_{0}^{|\xi|}\Phi(x,\tau)\ d\tau=|t|^{p(x)}A(x,\xi),
\end{equation}

\noindent which completes the proof.
\end{proof}

\begin{remark} We have proved that hypothesis \textnormal{\textbf{(A7)}} from \cite{giaco1} holds iff $\Phi$ is $p(x)-1$--homogeneous. Note that in that case hypothesis \textnormal{\textbf{(H7')}} is satisfied for any $r\in [1,p^-)$ since for a.a. $x\in\Omega$ we have that:
	
\begin{equation}
	(0,\infty)\ni s\mapsto\dfrac{\Phi(x,s)}{s^{r-1}}=\dfrac{s^{p(x)-1}\Phi(x,1)}{s^{r-1}}=s^{p(x)-r}\underbrace{\Phi(x,1)}_{>0}
\end{equation}

\noindent is a \textbf{strictly increasing} function on $(0,\infty)$. This is true because $p(x)\geq p^->r$. So, we have presented here a significantly more general version of the \textbf{D\'{i}az-Saa} inequality than the version that can be found in \cite[Theorem 4.2]{giaco1} and \cite[Theorem 2.4]{tak}. The present version applies in practical situations in which the old versions cannot be used (take a look at the problem \eqref{eqimage} from Section 8 that is important in image processing). Not to mention that in those versions there is also a convexity assumption of the mapping $\xi\mapsto A(x,\xi)^{r}{p(x)}:\mathbb{R}^N\to\mathbb{R}$, that in the present paper is not necessary. This omission will not affect the equality cases, that we will analyze in very deep details in the next section (and say more than in the the two cited papers).
\end{remark}

\section{Equality in D\'{i}az-Saa inequality for homogeneous Neumann conditions}

\begin{theorem}[\textnormal{\textbf{Equality in D\'{i}az-Saa inequality}}]\label{thmequalitydiaz} Let $w_1,w_2\in W^{1,p(x)}(\Omega)$ with $w_1,w_2>0$ a.e. on $\Omega$, $\dfrac{w_2}{w_1},\dfrac{w_1}{w_2}\in L^{\infty}(\Omega)$ and $\dfrac{\nabla w_1}{w_1},\ \dfrac{\nabla w_2}{w_2}\in L^1_{\textnormal{loc}}(\Omega)^N$ such that the following equality holds:
	
	\begin{equation}
		\int_{\Omega} \mathbf{a}(x,\nabla w_1)\cdot\nabla\left (w_1-\dfrac{w_2^r}{w_1^{r-1}}\right )\ dx= \int_{\Omega} \mathbf{a}(x,\nabla w_2)\cdot\nabla\left (\dfrac{w_1^r}{w_2^{r-1}}-w_2\right )\ dx.
	\end{equation}
\noindent The following statements are true:

\begin{enumerate}
	\item[(i)] If $r=1$ then there is some real number $c$ such that $w_1-w_2=c$ a.e. on $\Omega$ and equality holds.
	
	\item[(ii)] If $r>1$ then there is some $\lambda>0$ such that $w_2=\lambda w_1$ a.e. on $\Omega$. If $\lambda=1$, i.e. $w_1\equiv w_2$ then equality holds, but if $\lambda\neq 1$ we have equality if and only if:
	
	\begin{equation}
		\Phi(x,\lambda |\nabla w_1(x)|)=\lambda^{r-1}\Phi (x,|\nabla w_1(x)|),\ \text{for a.a.}\ x\in\Omega.
	\end{equation}
	
	\noindent In this situation, if in addition we know that the function $(0,\infty)\ni s\mapsto\dfrac{\Phi(x,s)}{s^{r-1}}$ is \textbf{strictly increasing} for a.a. $x\in\Omega$ (i.e. \textnormal{\textbf{(H7')}} holds), then equality holds iff $w_1\equiv w_2$ or both $w_1,w_2$ are different constant functions on $\Omega$.
\end{enumerate}	

\noindent There are no other situations in which equality can occur.
	
\end{theorem}

\begin{proof} In order to have equality in the \textbf{D\'{i}az-Saa inequality}, using Corollary \ref{cordiaz}, we need to have:
	
\begin{equation}\label{eqpointdiaz}
	\mathbf{a}(x,\nabla w_1)\cdot\nabla\left (w_1-\dfrac{w_2^r}{w_1^{r-1}}\right )= \mathbf{a}(x,\nabla w_2)\cdot\nabla\left (\dfrac{w_1^r}{w_2^{r-1}}-w_2\right ),\ \text{a.e. in}\ \Omega.
\end{equation}

\noindent Thus for a.a. $x\in\Omega$ both inequalities from \eqref{ineqeq} will become equalities.

\bigskip

\noindent\textbf{Case I:} If $r=1$ then we have that $\big (\mathbf{a}(x,\nabla w_1(x))-\mathbf{a}(x,\nabla w_2(x))\big )\cdot (\nabla w_1(x)-\nabla w_2(x))=0$ for a.a. $x\in\Omega$. Thus, using Proposition 3.2 \textbf{(1)} from \cite{max1}, we get that $\nabla w_1=\nabla w_2$ a.e. on $\Omega$. Then $\nabla (w_1-w_2)=0$ a.e. on $\Omega$ (which is a connected domain) i.e. $w_1-w_2$ is a constant function on $\Omega$. So (i) is proved.

\bigskip

\noindent\textbf{Case II:} $r>1$. We will prove step by step some facts:

\medskip

\noindent$\blacktriangleright$ \textbf{Fact I:} Fix some $x\in\Omega$. Then $\nabla w_1(x)=0$ iff $\nabla w_2(x)=0$. 

\medskip

\noindent Indeed if $\nabla w_1(x)=0$ then: $\left [1+(r-1)\left(\dfrac{w_1(x)}{w_2(x)}\right)^r\right ]\Phi(x,|\nabla w_2|)|\nabla w_2(x)|=0\ \Longrightarrow\ |\nabla w_2(x)|=0\ \Longrightarrow\ \nabla w_2(x)=0$. Similarly if $\nabla w_2(x)=0$ then: $\left [1+(r-1)\left(\dfrac{w_2(x)}{w_1(x)}\right)^r\right ]\Phi(x,|\nabla w_1|)|\nabla w_1(x)|=0\ \Longrightarrow\ |\nabla w_1(x)|=0\ \Longrightarrow\ \nabla w_1(x)=0$.

\noindent Now we can introduce the measurable set $\tilde{\Omega}=\{x\in\Omega\ |\ \nabla w_1\neq 0\}=\{x\in\Omega\ |\ \nabla w_2\neq 0\}=\{x\in\Omega\ |\ \nabla w_1,\nabla w_2\neq 0\}$. The equality between these three sets has to be understood in the measure sense (i.e. any two of them can differ on a set of null measure). So we have that $\nabla w_1(x)=\nabla w_2(x)=0$ for a.a. $x\in\Omega\setminus\tilde{\Omega}$. We may also write that $\dfrac{\nabla w_1(x)}{w_1(x)}=\dfrac{\nabla w_2(x)}{w_2(x)}=0$ for a.a. $x\in\Omega\setminus\tilde{\Omega}$.

\bigskip

\noindent\textbf{Fact II:} For a.a. $x\in\tilde{\Omega}$ we have that $\dfrac{\nabla w_1(x)}{w_1(x)}=\dfrac{\nabla w_2(x)}{w_2(x)}$.

\medskip

\noindent Since $x\in\tilde{\Omega}$ we have that $\nabla w_1(x),\nabla w_2(x)\neq 0$. But to have equalities in \eqref{ineqeq} we deduce that $|\nabla w_1(x)|\cdot |\nabla w_2(x)|=\nabla w_1(x)\cdot\nabla w_2(x)$. This means that there is some $\lambda(x)>0$ such that $\nabla w_2(x)=\lambda(x)\nabla w_1(x)$. Following the cases in which equality can occur in Lemma \ref{lemineqdiazsaa}, having $r>1$ and $|\nabla w_1(x)|,|\nabla w_2(x)|>0$, we are in the case (iii) from which we can write that: $\dfrac{|\nabla w_1(x)|}{w_1(x)}=\dfrac{|\nabla w_2(x)|}{w_2(x)}\ \Longrightarrow\ \dfrac{|\nabla w_1(x)|}{w_1(x)}=\dfrac{\lambda(x)|\nabla w_1(x)|}{w_2(x)}\ \Longrightarrow\ \lambda(x)=\dfrac{w_2(x)}{w_1(x)}$. So the claim is proved, since $\nabla w_2(x)=\dfrac{w_2(x)}{w_1(x)}\nabla w_1(x)$.

\bigskip

\noindent\textbf{Fact III:} The function $\dfrac{w_2}{w_1}$ is constant.

\medskip

\noindent At this moment we have that $\dfrac{\nabla w_1(x)}{w_1(x)}=\dfrac{\nabla w_2(x)}{w_2(x)}$ a.e. on the entire domain $\Omega$. Using now Lemma \ref{lemmamea} proved at the end of this proof, we deduce that $\dfrac{w_2}{w_1}\in W^{1,1}_{\text{loc}}(\Omega)$ and: $\nabla\left (\dfrac{w_2}{w_1} \right )=\dfrac{w_2}{w_1}\left (\dfrac{\nabla w_2}{w_2}-\dfrac{\nabla w_1}{w_1} \right )=0$ a.e. on $\Omega$. Since $\Omega$ is a connected domain we infer that $\dfrac{w_2}{w_1}$ is a constant function.

\noindent Let $w_2=\lambda w_1$ a.e. on $\Omega$ for some $\lambda>0$. Substituting this in \eqref{eqpointdiaz}  we obtain a.e. on $\Omega$ that:

\begin{align*}
	\mathbf{a}(x,\nabla w_1)\cdot\nabla\big (w_1-\lambda^r w_1\big )&=\mathbf{a}(x,\lambda\nabla w_1)\cdot\nabla\left (\dfrac{w_1}{\lambda^{r-1}}-\lambda w_1\right )\\
	\Longleftrightarrow\ \big (1-\lambda^r\big )\mathbf{a}(x,\nabla w_1)\cdot\nabla w_1&=\dfrac{1-\lambda^r}{\lambda^{r-1}}\mathbf{a}(x,\lambda\nabla w_1)\cdot\nabla w_1\\
	\Longleftrightarrow\ \big (1-\lambda^r\big )\lambda^{r-1}\Phi(x,|\nabla w_1|)|\nabla w_1|&=\big(1-\lambda^r\big )\Psi(x,\lambda|\nabla w_1|)\lambda\nabla w_1\cdot\nabla w_1\\
	\Longleftrightarrow\ \big (1-\lambda^r\big )\lambda^{r-1}\Phi(x,|\nabla w_1|)|\nabla w_1|&=\big(1-\lambda^r\big )\Phi(x,\lambda|\nabla w_1|)|\nabla w_1|\\
	\Longleftrightarrow\ \big (1-\lambda^r\big )|\nabla w_1|\big [\Phi(x,\lambda|\nabla w_1|)&-\lambda^{r-1}\Phi(x,|\nabla w_1|)\big ]=0.
\end{align*}

\noindent Therefore we may have $\lambda=1$ and the equality holds. If $\lambda\neq 1$ then $\nabla w_1(x)=0$ or $\Phi(x,\lambda|\nabla w_1|)=\lambda^{r-1}\Phi(x,|\nabla w_1|)$ for a.a. $x\in\Omega$. This disjunction is equivalent to the second relation since $\Phi(x,0)=0$.

\noindent If we know that the function $(0,\infty)\ni s\mapsto\dfrac{\Phi(x,s)}{s^{r-1}}$ is \textbf{strictly increasing} for a.a. $x\in\Omega$ (i.e. \textbf{(H7')} holds), then equality can occur when $\lambda=1$ (i.e. $w_1\equiv w_2$). If $\lambda\neq 1$ then for $x\in\tilde{\Omega}$ (i.e. $\nabla w_1(x)\neq 0$) we have that $\Phi(x,\lambda|\nabla w_1(x)|)=\lambda^{r-1}\Phi(x,|\nabla w_1(x)|)$ which rewrites as $\dfrac{\Phi(x,|\nabla w_1(x)|)}{|w_1(x)|^{r-1}}=\dfrac{\Phi(x,\lambda|\nabla w_1(x)|)}{\big (\lambda|w_1(x)|\big )^{r-1}}$. This equality is impossible to hold from the strict monotony of the function involved there. So $|\tilde{\Omega}|=0$ which means that $\nabla w_1=0$ a.e. on $\Omega$. Therefore $\nabla w_2=\lambda\nabla w_1=0$ a.e. on $\Omega$. Finally, since $\Omega$ is connected we have that $w_1$ and $w_2$ are different constant functions. The theorem is now completely proved.
\end{proof}

\noindent Now we shall state and prove the lemma we have used above.

\begin{lemma}\label{lemmamea}
	Let $\Omega\subset\mathbb{R}^N$ be an open, bounded and connected domain. If $w_1,w_2\in W^{1,p(x)}(\Omega)$ with:
	
	\begin{enumerate} 
		
		\item[(i)] $w_1,w_2>0$ a.e. on $\Omega$;
		\item[(ii)] $\dfrac{w_1}{w_2},\ \dfrac{w_2}{w_1}\in L^{\infty}(\Omega)$;
		\item[(iii)] $\dfrac{\nabla w_1}{w_1},\ \dfrac{\nabla w_2}{w_2}\in L^1_{\textnormal{loc}}(\Omega)^N$.
		
	\end{enumerate}
	
	Then $\dfrac{w_2}{w_1},\dfrac{w_1}{w_2}\in W^{1,1}_{\textnormal{loc}}(\Omega)$ and moreover:
	
	\begin{equation}
		\nabla\left (\dfrac{w_2}{w_1} \right )=\dfrac{w_1\nabla w_2-w_2\nabla w_1}{w_1^2}=\dfrac{w_2}{w_1}\left (\dfrac{\nabla w_2}{w_2}-\dfrac{\nabla w_1}{w_1} \right ).
	\end{equation} 
\end{lemma}

\begin{proof}
	For each $0<\epsilon<1$ consider the following function:
	
	\begin{equation}
		\eta_{\epsilon}:\mathbb{R}\to (0,\infty),\ \eta_{\epsilon}(w)=\begin{cases} \epsilon, & w\leq\epsilon\\ w, & w\in I_{\epsilon}:=\left(\epsilon,\dfrac{1}{\epsilon}\right)\\ \dfrac{1}{\epsilon}, & w\geq\dfrac{1}{\epsilon}\end{cases}.
	\end{equation}
	
	\noindent It is easy to check that $\eta_{\epsilon}$ is a $1-$Lipschitz function, i.e. $|\eta_{\epsilon}(\tilde{w})-\eta_{\epsilon}(w)|\leq |\tilde{w}-w|$ for any $\tilde{w},w\in \mathbb{R}$.

	\noindent Therefore, using the chain rule for weak derivatives, we deduce that $w_{1,\epsilon}:=\eta_{\epsilon}\circ w_1,\ w_{2,\epsilon}:=\eta_{\epsilon}\circ w_2\in W^{1,p(x)}(\Omega)\cap L^{\infty}(\Omega)$ and moreover for a.a. $x\in\Omega$:
	
	\begin{equation}\label{apeqdiazimp10}
		\begin{cases} \nabla w_{1,\epsilon}=\nabla \eta_{\epsilon}\circ w_1=\chi_{I_{\epsilon}}(w_1)\nabla w_1\\ \nabla w_{2,\epsilon} \nabla \eta_{\epsilon}\circ w_2=\chi_{I_{\epsilon}}(w_2)\nabla w_2\end{cases}\ \Longrightarrow\ \begin{cases} |\nabla w_{1,\epsilon}|\leq|\nabla w_1|\\ |\nabla w_{2,\epsilon}|\leq |\nabla w_2|\end{cases},\ \forall\ \epsilon>0.
	\end{equation}

	\noindent The most important fact here is to notice that:
	
	\begin{equation}\label{apeqdiazimp20}
		\begin{cases}\dfrac{w_{1,\epsilon}}{w_{2,\epsilon}}\leq\max\left\{\dfrac{w_1}{w_2},1\right\}, \text{a.e. on}\ \Omega\\[3mm] \dfrac{w_{2,\epsilon}}{w_{1,\epsilon}}\leq\max\left\{\dfrac{w_2}{w_1},1\right\}, \text{a.e. on}\ \Omega\end{cases}, \forall\ \epsilon>0.
	\end{equation}
	
	\noindent Setting $M:=\max\left\{\left\Vert\dfrac{w_1}{w_2} \right\Vert_{L^{\infty}(\Omega)},\left\Vert\dfrac{w_2}{w_1} \right\Vert_{L^{\infty}(\Omega)}\right\}\geq 1$ we have that $\dfrac{w_{1,\epsilon}}{w_{2,\epsilon}},\dfrac{w_{2,\epsilon}}{w_{1,\epsilon}}\leq M$ a.e. on $\Omega$ for any $\epsilon>0$.

	\noindent Fix some arbitrary $i\in\overline{1,N}$. The function $\mathbb{R}\ni s\longmapsto\begin{cases} \dfrac{1}{\epsilon}, & s\leq\epsilon\\[3mm] \dfrac{1}{s}, & s\in \left (\epsilon,\dfrac{1}{\epsilon}\right )\\[3mm]\epsilon, & s\geq \dfrac{1}{\epsilon}\end{cases}$ is $\dfrac{1}{\epsilon^2}-$Lipschitz. Therefore from the \textit{chain rule for weak derivatives} we have that $\dfrac{1}{w_{1,\varepsilon}}\in W^{1,p(x)}(\Omega)$ and:
	
	\begin{equation}\label{apeqdiaz1} \dfrac{\partial}{\partial x_i}\left(\dfrac{1}{w_{1,\epsilon}}\right)=-\dfrac{1}{w_{1,\epsilon}^2}\chi_{I_{\epsilon}}(w_{1,\epsilon})\dfrac{\partial w_{1,\epsilon}}{\partial x_i}=-\dfrac{1}{w_{1,\epsilon}^2}\dfrac{\partial w_{1,\epsilon}}{\partial x_i},
	\end{equation} 
	
	\noindent since $\dfrac{1}{\epsilon}\geq w_{1,\epsilon}\geq \epsilon$. This also shows that $\dfrac{1}{w_{1,\epsilon}}\in L^{\infty}(\Omega)$.
	
	\noindent Now, since $w_{2,\epsilon},\dfrac{1}{w_{1,\epsilon}}\in W^{1,p(x)}(\Omega)\cap L^{\infty}(\Omega)$, we get from the \textit{product rule for weak derivatives} that $\dfrac{w_{2,\epsilon}}{w_{1,\epsilon}}\in W^{1,p(x)}(\Omega)\cap L^{\infty}(\Omega)$ and using \eqref{apeqdiaz1}:
	
	\begin{equation}\label{apeqdiaz2}
		\dfrac{\partial}{\partial x_i}\left(\dfrac{w_{2,\epsilon}}{w_{1,\epsilon}}\right)=\dfrac{1}{w_{1,\epsilon}}\dfrac{\partial w_{2,\epsilon}}{\partial x_i}+w_{2,\epsilon}\dfrac{\partial}{\partial x_i}\left(\dfrac{1}{w_{1,\epsilon}}\right)=\dfrac{1}{w_{1,\epsilon}}\dfrac{\partial w_{2,\epsilon}}{\partial x_i}-\dfrac{w_{2,\epsilon}}{w_{1,\epsilon}^2}\dfrac{\partial w_{1,\epsilon}}{\partial x_i}.
	\end{equation}
	
	\noindent So, for any $\phi\in C^{\infty}_c(\Omega)$ we may write:
	
	\begin{equation}\label{eqlimit}
		\int_{\Omega} \dfrac{w_{2,\epsilon}}{w_{1,\epsilon}}\dfrac{\partial\phi}{\partial x_i}\ dx=-\int_{\Omega} \left [\dfrac{1}{w_{1,\epsilon}}\dfrac{\partial w_{2,\epsilon}}{\partial x_i}-\dfrac{w_{2,\epsilon}}{w_{1,\epsilon}^2}\dfrac{\partial w_{1,\epsilon}}{\partial x_i} \right ]\phi\ dx.
	\end{equation}
	
	\noindent We want to pass $\epsilon\searrow 0$ in \eqref{eqlimit}. In the left hand side there is no problem since:
	
	\begin{equation}
		\left |\int_{\Omega}\left (\dfrac{w_{2,\epsilon}}{w_{1,\epsilon}}-\dfrac{w_2}{w_1}\right )\dfrac{\partial\phi}{\partial x_i} \ dx \right |\leq\left\Vert\dfrac{\partial\phi}{\partial x_i}\right\Vert_{L^{\infty}(\Omega)}\int_{\Omega}\left |\dfrac{w_{2,\epsilon}}{w_{1,\epsilon}}-\dfrac{w_2}{w_1} \right |\ dx\stackrel{\epsilon\searrow 0}{\longrightarrow} 0,
	\end{equation}
	
	\noindent from \textit{Lebesgue dominated convergence theorem}, since:
	
	\begin{enumerate}
		\item[$\bullet$] $\lim\limits_{\epsilon\searrow 0} \dfrac{w_{2,\epsilon}}{w_{1,\epsilon}}-\dfrac{w_2}{w_1}=0$ pointwise, a.e. on $\Omega$,
		
		\item[$\bullet$] $\left |\dfrac{w_{2,\epsilon}}{w_{1,\epsilon}}-\dfrac{w_2}{w_1} \right |\leq M$ for each $1>\epsilon>0$.
	\end{enumerate}
	
	\noindent For the right-hand side of \eqref{eqlimit} we denote $\tilde{\Omega}=\operatorname{supp}(\phi)\Subset\Omega$ and remark first that:
	
	\begin{align}
		\int_{\Omega} \left [\dfrac{1}{w_{1,\epsilon}}\dfrac{\partial w_{2,\epsilon}}{\partial x_i}-\dfrac{w_{2,\epsilon}}{w_{1,\epsilon}^2}\dfrac{\partial w_{1,\epsilon}}{\partial x_i} \right ]\phi\ dx&=\int_{\tilde{\Omega}} \left [\dfrac{1}{w_{1,\epsilon}}\dfrac{\partial w_{2,\epsilon}}{\partial x_i}-\dfrac{w_{2,\epsilon}}{w_{1,\epsilon}^2}\dfrac{\partial w_{1,\epsilon}}{\partial x_i} \right ]\phi\ dx\\
		&=\int_{\tilde{\Omega}} \left [\dfrac{w_{2,\epsilon}}{w_{1,\epsilon}}\dfrac{\partial_{x_i}w_{2}}{w_{2,\epsilon}}\chi_{I_{\epsilon}}(w_2)-\dfrac{w_{2,\epsilon}}{w_{1,\epsilon}}\dfrac{\partial_{x_i}w_{1}}{w_{1,\epsilon}}\chi_{I_{\epsilon}}(w_1) \right ]\phi\ dx\\
		&=\int_{\tilde{\Omega}} \dfrac{w_{2,\epsilon}}{w_{1,\epsilon}}\left [\dfrac{\partial_{x_i}w_{2}}{w_{2}}\chi_{I_{\epsilon}}(w_2)-\dfrac{\partial_{x_i}w_{1}}{w_{1}}\chi_{I_{\epsilon}}(w_1) \right ]\phi\ dx\\
		&\stackrel{\epsilon\searrow 0}{\longrightarrow} \int_{\tilde{\Omega}} \dfrac{w_{2}}{w_{1}}\left [\dfrac{\partial_{x_i}w_{2}}{w_{2}}-\dfrac{\partial_{x_i}w_{1}}{w_{1}}\right ]\phi\ dx,
	\end{align}
	
	\noindent from \textit{Lebesgue dominated convergence theorem}, as:
	
	\begin{enumerate}
		\item[$\bullet$] $\lim\limits_{\epsilon\searrow 0}\dfrac{w_{2,\epsilon}}{w_{1,\epsilon}}\left [\dfrac{\partial_{x_i}w_{2}}{w_{2}}\chi_{I_{\epsilon}}(w_2)-\dfrac{\partial_{x_i}w_{1}}{w_{1}}\chi_{I_{\epsilon}}(w_1) \right ]\phi=\dfrac{w_{2}}{w_{1}}\left [\dfrac{\partial_{x_i}w_{2}}{w_{2}}-\dfrac{\partial_{x_i}w_{1}}{w_{1}}\right ]\phi$ pointwise, a.e. on $\Omega$,
		
		\item[$\bullet$] $\left |\dfrac{w_{2,\epsilon}}{w_{1,\epsilon}}\left [\dfrac{\partial_{x_i}w_{2}}{w_{2}}\chi_{I_{\epsilon}}(w_2)-\dfrac{\partial_{x_i}w_{1}}{w_{1}}\chi_{I_{\epsilon}}(w_1) \right ]\phi \right |\leq M\Vert\phi\Vert_{L^{\infty}(\tilde{\Omega})}\left (\left |\dfrac{\partial_{x_i}w_{2}}{w_{2}} \right |+\left | \dfrac{\partial_{x_i}w_{1}}{w_{1}}\right | \right )\in L^{1}(\tilde{\Omega})$ from (iii) for any $\epsilon\in (0,1)$.
	\end{enumerate}
	
	\noindent Thus, making $\epsilon\searrow 0$ in \eqref{eqlimit} we get for each $i\in\overline{1,N}$ that:
	
	\begin{equation}
		\int_{\Omega} \dfrac{w_{2}}{w_{1}}\dfrac{\partial\phi}{\partial x_i}\ dx=-\int_{\Omega} \left [\dfrac{1}{w_{1}}\dfrac{\partial w_{2}}{\partial x_i}-\dfrac{w_{2}}{w_{1}^2}\dfrac{\partial w_{1}}{\partial x_i} \right ]\phi\ dx,\ \forall\ \phi\in C^{\infty}_c(\Omega).
	\end{equation}
	
	\noindent This shows that $\dfrac{w_2}{w_1}$ has a weak gradient and it is given by:
	
	\begin{equation}
		\nabla\left (\dfrac{w_2}{w_1} \right )=\dfrac{w_1\nabla w_2-w_2\nabla w_1}{w_1^2}=\dfrac{w_2}{w_1}\left (\dfrac{\nabla w_2}{w_2}-\dfrac{\nabla w_1}{w_1} \right )\in L^{1}_{\text{loc}}(\Omega).
	\end{equation}
	
	\noindent Taking again into account that $\dfrac{w_2}{w_1}\in L^{\infty}(\Omega)$ we can conclude that $\dfrac{w_2}{w_1}\in W^{1,1}_{\text{loc}}(\Omega)$.
\end{proof}

\bigskip

\begin{remark} The condition $\dfrac{\nabla w_1}{w_1},\dfrac{\nabla w_2}{w_2}\in L^1_{\textnormal{loc}}(\Omega)^N$ that is introduced in the present paper is essential. Here we provide two counterexamples in which $\dfrac{\nabla w_1}{w_1}=\dfrac{\nabla w_2}{w_2}\Rightarrow \nabla\left (\dfrac{w_1}{w_2}\right )=0$ but $\dfrac{w_2}{w_1}$ fails to be a constant function.
\end{remark}

\begin{example}\label{ex51} Set $\Omega=(-1,1)$ and define $w_1:\Omega\to \mathbb{R},\ w_1(x)=|x|$, $w_2:\Omega\to\mathbb{R},\ w_2(x)=\begin{cases} x, & x\geq 0\\ -2x, & x<0\end{cases}$. Then $w_1,w_2\in W^{1,p}(\Omega)$ for any $p>1$ (because they are Lipschitz functions), $w_1,w_2>0$ a.e. on $\Omega$ (only if $x=0$ we have $w_1(x)=w_2(x)=0$), $\dfrac{w_1}{w_2}\in\left\{1,\dfrac{1}{2}\right\},\dfrac{w_2}{w_1}\in\{1,2\}$ are both in $L^{\infty}(\Omega)$ and $\dfrac{w'_1(x)}{w_1(x)}=\dfrac{w'_2(x)}{w_2(x)}=\dfrac{1}{x}$ for any $x\in (-1,1)\setminus\{0\}$ (so for a.a. $x\in\Omega$). However: $\dfrac{w_2(x)}{w_1(x)}=\begin{cases} 1, & x> 0\\ 2, & x<0\end{cases}$ which is not a constant function.
\end{example}

\begin{example}\label{ex52} We can adjust the previous example to obtain another one with function from $W^{1,p}_0(\Omega)$. For this purpose consider the cutt-off function $\eta:(-1,1)\to (0,\infty), \eta(x)=\exp\left(\dfrac{1}{x^2-1}\right )$ ($\eta\in C^{\infty}_c([-1,1]),\ \operatorname{supp}(\eta)=[-1,1]$), and define $\tilde{w}_1:\Omega\to\mathbb{R},\ \tilde{w}_1(x)=w_1(x)\eta(x)$ and $\tilde{w}_2:\Omega\to\mathbb{R},\  \tilde{w}_2(x)=w_2(x)\eta(x)$. Clearly $\tilde{w}_1,\tilde{w}_2\in W^{1,p}_0(\Omega)$ for any $p>1$, because $\tilde{w}_1,\tilde{w}_2\in C(\overline{\Omega})$ (even Lipschitz) and $\tilde{w}_1(1)=\tilde{w}_2(1)=\tilde{w}_1(-1)=\tilde{w}_2(-1)=0$ (see \cite[Theorem 9.17]{brezis}). Again $\tilde{w}_1,\tilde{w}_2>0$ a.e. on $\Omega$, $\dfrac{\tilde{w}_1}{\tilde{w}_2}\in\left\{1,\frac{1}{2}\right\},\ \dfrac{\tilde{w}_2}{\tilde{w}_1}\in\{1,2\}$, so both functions belong to $L^{\infty}(\Omega)$. Note that $\dfrac{ \tilde{w}'_1(x)}{\tilde{w}_1(x)}=\dfrac{w'_1(x)}{w_1(x)}+\dfrac{\eta'(x)}{\eta(x)}=\dfrac{1}{x}+\dfrac{\eta'(x)}{\eta(x)}=\dfrac{w'_2(x)}{w_2(x)}+\dfrac{\eta'(x)}{\eta(x)}=\dfrac{\tilde{w}'_2(x)}{\tilde{w}_2(x)}$ for any $x\in (-1,1)\setminus\{0\}$. So, all the hypotheses are fulfilled but $\dfrac{\tilde{w}_1}{\tilde{w}_2}=\dfrac{w_2}{w_1}$ fails to be a constant. This example is important because it offers a refinement of Theorem 2.2 from \cite{tak} which does not assume in its statement the condition $\dfrac{\nabla w_1}{w_1},\dfrac{\nabla w_2}{w_2}\in L^1_{\textnormal{loc}}(\Omega)^N$ but in equation (2.7) from page 7 the authors used a result like Lemma \ref{lemmamea} from this paper.
\end{example}

\section{Uniqueness of the steady-state}

\begin{proposition}\label{propwtheta} Fix some $w_1,w_2\in\overset{\bullet}{W}$ with $\dfrac{w_1}{w_2},\ \dfrac{w_2}{w_1}\in L^{\infty}(\Omega)$. We denote $w_{\theta}=\theta w_1+(1-\theta)w_2$ for each $\theta\in\mathbb{R}$. Then for $M:=\max\left\{\left\Vert\dfrac{w_1}{w_2} \right\Vert_{L^{\infty}(\Omega)},\left\Vert\dfrac{w_2}{w_1} \right\Vert_{L^{\infty}(\Omega)}\right\}\geq 1$ and $\theta_0:=\dfrac{1}{2(M-1)}\in (0,\infty]$ the following properties hold:
	
	\begin{enumerate}
		\item[\textnormal{\textbf{(1)}}] $\dfrac{1}{2}\min\{w_1,w_2\}\leq w_\theta\leq\dfrac{3}{2}\max\{w_1,w_2\}$ for any $\theta\in (-\theta_0,1+\theta_0)$.
		
		\item[\textnormal{\textbf{(2)}}] $\dfrac{1}{M}\leq\dfrac{w_{2,\epsilon}}{w_{1,\epsilon}},\dfrac{w_{1,\epsilon}}{w_{2,\epsilon}}\leq M$ for any $\epsilon\in (0,1)$.
		
		\item[\textnormal{\textbf{(3)}}] $\dfrac{\epsilon}{2}\leq\dfrac{1}{2}\min\{w_{1,\epsilon};w_{2,\epsilon}\}\leq w_{\theta,\epsilon}\leq\dfrac{3}{2}\max\{w_{1,\epsilon};w_{2,\epsilon}\}\leq\dfrac{3}{2\epsilon}$ for any $\theta\in (-\theta_0,1+\theta_0)$, where $w_{\theta,\epsilon}:=\theta w_{1,\epsilon}+(1-\theta)w_{2,\epsilon}$.
		
		\item[\textnormal{\textbf{(4)}}] $\dfrac{2}{3M}\leq\dfrac{w_1}{w_{\theta}},\dfrac{w_2}{w_{\theta}}\leq 2M$ for any $\theta\in (-\theta_0,1+\theta_0)$.
		
		\item[\textnormal{\textbf{(5)}}] $\dfrac{2}{3M}\leq\dfrac{w_{1,\epsilon}}{w_{\theta,\epsilon}},\dfrac{w_{2,\epsilon}}{w_{\theta,\epsilon}}\leq 2M$ for any $\theta\in (-\theta_0,1+\theta_0)$.

		\item[\textnormal{\textbf{(6)}}] $w_{\theta}\in\overset{\bullet}{W}$ and $\nabla w_{\theta}^{\frac{1}{\alpha}}=\dfrac{1}{\alpha}w_{\theta}^{\frac{1}{\alpha}-1}\nabla w_{\theta}$ for any $\theta\in (-\theta_0,1+\theta_0)$.
		
		\item[\textnormal{\textbf{(7)}}] $\gamma_{\theta}:=\dfrac{1}{\alpha}\cdot\dfrac{w_1-w_2}{\sqrt[\alpha]{w_{\theta}^{\alpha-1}}}\in W^{1,p(x)}(\Omega)$ and $\nabla\gamma_{\theta}=\dfrac{1}{\alpha}\dfrac{\nabla w_{1}-\nabla w_{2}}{w_{\theta}^{1-\frac{1}{\alpha}}}+\dfrac{1}{\alpha}\left(\dfrac{1}{\alpha}-1\right)\dfrac{w_{1}-w_{2}}{w_{\theta}^{2-\frac{1}{\alpha}}}\nabla w_{\theta}$ for any $\theta\in (-\theta_0,1+\theta_0)$.
		
		\item[\textnormal{\textbf{(8)}}] Let $\Gamma: (-\theta_0,1+\theta_0)\to W^{1,p(x)}(\Omega),\ \Gamma(\theta)=\sqrt[\alpha]{w_{\theta}}\in W^{1,p(x)}(\Omega)$.\footnote{From \textbf{(6)} the function $\Gamma$ is well-defined.} Then $\Gamma$ is differentiable on $(-\theta_0,1+\theta_0)$ and:
		
		\begin{equation}
			\Gamma'(\theta)=\gamma_{\theta}=\dfrac{1}{\alpha}\cdot\dfrac{w_1-w_2}{\sqrt[\alpha]{w_{\theta}^{\alpha-1}}}\in W^{1,p(x)}(\Omega),\ \forall\ \theta\in (-\theta_0,1+\theta_0).
		\end{equation}
	\end{enumerate}	
\end{proposition}

\begin{proof} \textbf{(1)}  Indeed, if $\theta\in (-\theta_0,0)$:
	
	 $w_{\theta}=\begin{cases} w_1\left [1+\theta\left (\dfrac{w_2}{w_1}-1 \right ) \right]\geq w_1 [1+\theta(M-1)]\geq w_1 [1-\theta_0(M-1)]=\dfrac{w_1}{2}\\[3mm] w_1\left [1-(-\theta)\left (\dfrac{w_2}{w_1}-1 \right ) \right ]\leq w_1\left [1-(-\theta)\left (\dfrac{1}{M}-1\right ) \right ]\leq w_1\left (1+\theta_0\dfrac{M-1}{M} \right )\leq\dfrac{3w_1}{2}\end{cases}$.
	
	\medskip
	
	\noindent If $\theta\in (1,1+\theta_0)$: 
	
	$w_{\theta}=\begin{cases} w_2\left [(\theta-1)\left (1-\dfrac{w_1}{w_2} \right )+1 \right ]\geq w_2 [1-(\theta-1)(M-1)]\geq w_2[1-\theta_0(M-1)]=\dfrac{w_2}{2}\\[3mm] w_2\left [(\theta-1)\left (1-\dfrac{w_1}{w_2} \right )+1 \right ]\leq w_2\left [(\theta-1)\left (1-\dfrac{1}{M} \right )+1 \right ]\leq w_2\left [\theta_0\left (1-\dfrac{1}{M} \right )+1 \right ]\leq\dfrac{3w_2}{2}\end{cases}$.
	
	\medskip
	
	\noindent If $\theta\in [0,1]|$ there is nothing to prove since $\min\{w_1,w_2\}\leq w_{\theta}\leq \max\{w_1,w_2\}$.
	
	\noindent From this fact it follows that $w_{\theta}>0$ a.e. on $\Omega$ for any $\theta\in (-\theta_0,1+\theta_0)$.
	
	\bigskip
	
	\noindent\textbf{(2)} If $w_{1}\leq \epsilon$ then $\dfrac{w_{1,\epsilon}}{w_{2,\epsilon}}=\begin{cases}\frac{1}{M}\leq\frac{\epsilon}{\epsilon}=1\leq M, & w_2\leq \epsilon\\ \frac{1}{M}\leq\frac{w_1}{w_2}\leq\frac{\epsilon}{w_2}\leq 1\leq M, & w_2\in I_{\epsilon}\\ \frac{1}{M}\leq\frac{w_1}{w_2}\leq\frac{\epsilon}{\frac{1}{\epsilon}}=\epsilon^2< 1\leq M, & w_2\geq\frac{1}{\epsilon} \end{cases}$. If $w_1\in I_{\epsilon}$ then:

	$\dfrac{w_{1,\epsilon}}{w_{2,\epsilon}}=\begin{cases}\frac{1}{M}\leq 1=\frac{\epsilon}{\epsilon}\leq\frac{w_1}{\epsilon}\leq\frac{w_1}{w_2}\leq M, & w_2\leq\epsilon \\ \frac{1}{M}\leq\frac{w_1}{w_2}\leq M, & w_2\in I_{\epsilon}\\ \frac{1}{M}\leq\frac{w_1}{w_2}\leq\frac{w_1}{\frac{1}{\epsilon}}\leq 1\leq M, & w_2\geq\frac{1}{\epsilon} \end{cases}$. Also, if $w_1\geq\frac{1}{\epsilon}$, then $\dfrac{w_{1,\epsilon}}{w_{2,\epsilon}}=\begin{cases}\frac{1}{M}\leq 1\leq\frac{\frac{1}{\epsilon}}{\epsilon}\leq\frac{w_1}{w_2}\leq M, & w_{2}\leq\epsilon \\ \frac{1}{M}\leq 1\leq\frac{\frac{1}{\epsilon}}{w_2}\leq\frac{w_1}{w_2}\leq M, & w_2\in I_{\epsilon}\\ \frac{1}{M}\leq\frac{\frac{1}{\epsilon}}{\frac{1}{\epsilon}}=1\leq M, & w_2\geq\frac{1}{\epsilon} \end{cases}$. Thus $\dfrac{1}{M}\leq\dfrac{w_{1,\epsilon}}{w_{2,\epsilon}}\leq M$, which also means that $\dfrac{1}{M}\leq\dfrac{w_{2,\epsilon}}{w_{1,\epsilon}}\leq M$ a.e. on $\Omega$.
	
	\bigskip
	
	\noindent\textbf{(3)} The inequalities follow in the same manner as in \textbf{(1)}, taking into account \textbf{(2)}.

	\bigskip
	
	\noindent\textbf{(4)} From \textbf{(1)} we have that $\dfrac{2}{3M}\leq \dfrac{2}{3\max\{1,\frac{w_2}{w_1}\}}=\dfrac{2w_1}{3\max\{w_1,w_2\}}\leq\dfrac{w_1}{w_\theta}\leq \dfrac{2w_1}{\min\{w_1,w_2\}}=\dfrac{2}{\min\{1,\frac{w_2}{w_1}\}}\leq 2M$. Similarly $\dfrac{2}{3M}\leq \dfrac{2}{3\max\{1,\frac{w_1}{w_2}\}}=\dfrac{2w_2}{3\max\{w_1,w_2\}}\leq\dfrac{w_2}{w_\theta}\leq \dfrac{2w_2}{\min\{w_1,w_2\}}=\dfrac{2}{\min\{1,\frac{w_1}{w_2}\}}\leq 2M$.

	\bigskip
	
	\noindent\textbf{(5)} Using \textbf{(3)} we proceed exactly as in \textbf{(4)}.
	
	\bigskip

	\noindent\textbf{(6)} Since $w_1,w_2\in L^{\frac{p(x)}{\alpha}}(\Omega)$ it follows that $w_{\theta}=\theta w_1+(1-\theta)w_2\in L^{\frac{p(x)}{\alpha}}(\Omega)$. 
	
	\noindent Knowing that $w_1,w_2\in\overset{\bullet}{W}$, we can use Lemma \ref{lemmatrunc} \textbf{(3)} to deduce that $w_{1,\epsilon},w_{2,\epsilon}\in W^{1,p(x)}(\Omega)$ and $\nabla w_{1,\epsilon}=\alpha\sqrt[\alpha]{w_1^{\alpha-1}}\nabla\sqrt[\alpha]{w_1}\chi_{I_{\epsilon}}(w_1), \ \nabla w_{2,\epsilon}=\alpha\sqrt[\alpha]{w_2^{\alpha-1}}\nabla\sqrt[\alpha]{w_2}\chi_{I_{\epsilon}}(w_2)$. Thus $w_{\theta,\epsilon}=\theta w_{1,\epsilon}+(1-\theta)w_{2,\epsilon}\in W^{1,p(x)}$.
	
	\noindent The function $\mathbb{R}\ni s\mapsto\begin{cases}\left ( \dfrac{\epsilon}{2}\right )^{\frac{1}{\alpha}}, & s<\frac{\epsilon}{2}\\[3mm] s^{\frac{1}{\alpha}}, & s\in\left [\frac{\epsilon}{2},\frac{3}{2\epsilon} \right ]\\[3mm] \left (\dfrac{3}{2\epsilon} \right )^{\frac{1}{\alpha}}, & s> \frac{3}{2\epsilon}\end{cases}$ is $\dfrac{1}{\alpha}\left (\dfrac{2}{\epsilon} \right )^{1-\frac{1}{\alpha}}$--Lipschitz. When we compose it with $w_{\theta,\epsilon}\in \left [\frac{\epsilon}{2},\frac{3}{2\epsilon} \right ]$ we obtain precisely $\sqrt[\alpha]{w_{\theta,\epsilon}}$. Thus from the \textit{chain rule} we get that $\sqrt[\alpha]{w_{\theta,\epsilon}}\in W^{1,p(x)}(\Omega)$ and $\nabla\sqrt[\alpha]{w_{\theta,\epsilon}}=\dfrac{1}{\alpha} w_{\theta,\epsilon}^{\frac{1}{\alpha}-1}\nabla w_{\theta,\epsilon}$. This allows us to write for each $i\in\overline{1,N}$ that:
	
	\begin{equation}\label{eqlimit2}
		\int_{\Omega} w_{\theta,\epsilon}^{\frac{1}{\alpha}}\dfrac{\partial\phi}{\partial x_i}\ dx=-\int_{\Omega}\dfrac{1}{\alpha}w_{\theta,\epsilon}^{\frac{1}{\alpha}-1}\dfrac{\partial w_{\theta,\epsilon}}{\partial x_i}\phi\ dx,\ \forall\ \phi\in C^{\infty}_c(\Omega).
	\end{equation}	
	
	\noindent For the left-hand side we have that:
	
	\begin{align*}
		\left |\int_{\Omega}w_{\theta,\epsilon}^{\frac{1}{\alpha}}\dfrac{\partial\phi}{\partial x_i}\ dx-\int_{\Omega}w_{\theta}^{\frac{1}{\alpha}}\dfrac{\partial\phi}{\partial x_i}\ dx \right |&\leq \left\Vert\dfrac{\partial\phi}{\partial x_i}\right\Vert_{L^{\infty}(\Omega)}\int_{\Omega} \left |w_{\theta,\epsilon}^{\frac{1}{\alpha}}-
		w_{\theta}^{\frac{1}{\alpha}} \right |\ dx\\
	\text{(Lemma \ref{lemrunu} \textbf{(1)})}\ \ \ \ \ \ \ \ 	&\leq \left\Vert\dfrac{\partial\phi}{\partial x_i}\right\Vert_{L^{\infty}(\Omega)}\int_{\Omega} \left |w_{\theta,\epsilon}-
	w_{\theta} \right |^{\frac{1}{\alpha}}\ dx\\
	&=\left\Vert\dfrac{\partial\phi}{\partial x_i}\right\Vert_{L^{\infty}(\Omega)}\int_{\Omega} \left |\theta (w_{1,\epsilon}-w_1)+(1-\theta)(w_{2,\epsilon}-w_2)
	 \right |^{\frac{1}{\alpha}}\ dx \\
	 &\leq \left\Vert\dfrac{\partial\phi}{\partial x_i}\right\Vert_{L^{\infty}(\Omega)}\left [|\theta|^{\frac{1}{\alpha}}\int_{\Omega} |w_{1,\epsilon}-w_1|^{\frac{1}{\alpha}}\ dx+|1-\theta|^{\frac{1}{\alpha}}\int_{\Omega}|w_{2,\epsilon}-w_2|^{\frac{1}{\alpha}}\ dx  \right ]\\
	 &\leq \left\Vert\dfrac{\partial\phi}{\partial x_i}\right\Vert_{L^{\infty}(\Omega)} |1+\theta_0|^{\frac{1}{\alpha}}\left (\int_{\Omega} |w_{1,\epsilon}-w_1|^{\frac{1}{\alpha}}\ dx+\int_{\Omega}|w_{2,\epsilon}-w_2|^{\frac{1}{\alpha}}\ dx  \right )\\
\text{(Lemma \ref{lemmatrunc} \textbf{(6)} with $r=\frac{1}{\alpha}$)}\ \ \ \ \ \ \ 	&\stackrel{\epsilon\searrow 0}{\longrightarrow} 0.
	\end{align*}
	
	\noindent For the right-hand side of \eqref{eqlimit2} remark that:
	
 \begin{align*}
 	\lim\limits_{\epsilon\searrow 0}\dfrac{1}{\alpha}w_{\theta,\epsilon}^{\frac{1}{\alpha}-1}\dfrac{\partial w_{\theta,\epsilon}}{\partial x_i}\phi&=\lim\limits_{\epsilon\searrow 0}\frac{1}{\alpha}\big [\theta w_{1,\epsilon}+(1-\theta)w_{2,\epsilon}\big ]^{\frac{1}{\alpha}-1}\left [\theta\dfrac{\partial w_{1,\epsilon}}{\partial x_i}+(1-\theta)\dfrac{\partial w_{2,\epsilon}}{\partial x_i} \right ]\phi \\
\text{(Lemma \ref{lemmatrunc} \textbf{(1)} and \textbf{(8)})}\ \ \ \ \ \ \ \  	&=\frac{1}{\alpha}\big [\theta w_{1}+(1-\theta)w_{2}\big ]^{\frac{1}{\alpha}-1}\left [\theta\dfrac{\partial w_{1}}{\partial x_i}+(1-\theta)\dfrac{\partial w_{2}}{\partial x_i} \right ]\phi \\
 	&=\dfrac{1}{\alpha}w_{\theta}^{\frac{1}{\alpha}-1}\dfrac{\partial w_{\theta}}{\partial x_i}\phi\ \text{a.e. in }\Omega.
 \end{align*}
	
\noindent Now, for each $\epsilon\in (0,1)$ we have that:		
\begin{align*}
\left |\dfrac{1}{\alpha}w_{\theta,\epsilon}^{\frac{1}{\alpha}-1}\dfrac{\partial w_{\theta,\epsilon}}{\partial x_i}\phi \right |&\leq\Vert\phi\Vert_{L^{\infty}(\Omega)}\frac{1}{\alpha}\left |\dfrac{\theta\frac{\partial w_{1,\epsilon}}{\partial x_i}+(1-\theta)\frac{\partial w_{2,\epsilon}}{\partial x_i}}{w_{\theta,\epsilon}^{1-\frac{1}{\alpha}}}\right |\\
\text{(Lemma \ref{lemmatrunc} \textbf{(3)})}\ \ \ \ \ \ \ &=\Vert\phi\Vert_{L^{\infty}(\Omega)}\left |\dfrac{\theta w_{1}^{1-\frac{1}{\alpha}}\frac{\partial w_1^{\frac{1}{\alpha}}}{\partial x_i}\chi_{I_{\epsilon}}(w_1)+(1-\theta) w_{2}^{1-\frac{1}{\alpha}}\frac{\partial w_2^{\frac{1}{\alpha}}}{\partial x_i}\chi_{I_{\epsilon}}(w_2)}{w_{\theta,\epsilon}^{1-\frac{1}{\alpha}}}\right |\\
&\leq \Vert\phi\Vert_{L^{\infty}(\Omega)}\left [|\theta|\left (\dfrac{w_{1,\epsilon}}{w_{\theta,\epsilon}} \right )^{1-\frac{1}{\alpha}}\left |\frac{\partial w_1^{\frac{1}{\alpha}}}{\partial x_i} \right |+|1-\theta| \left (\dfrac{w_{2,\epsilon}}{w_{\theta,\epsilon}} \right )^{1-\frac{1}{\alpha}}\left |\frac{\partial w_2^{\frac{1}{\alpha}}}{\partial x_i} \right |\right ]\\
\text{\textbf{(5)}}\ \ \ \ \ \ \ \ &\leq \Vert\phi\Vert_{L^{\infty}(\Omega)}(2M)^{1-\frac{1}{\alpha}}|1+\theta_0|\left ( \left |\frac{\partial w_1^{\frac{1}{\alpha}}}{\partial x_i} \right |+\left |\frac{\partial w_2^{\frac{1}{\alpha}}}{\partial x_i} \right | \right )\in L^{p(x)}(\Omega),
\end{align*}

\noindent because $w_1^{\frac{1}{\alpha}},w_2^{\frac{1}{\alpha}}\in W^{1,p(x)}(\Omega)$. From \textit{Lebesgue dominated convergence theorem} we deduce that:

\[
\lim\limits_{\epsilon\searrow 0} \displaystyle\int_{\Omega}\dfrac{1}{\alpha}w_{\theta,\epsilon}^{\frac{1}{\alpha}-1}\dfrac{\partial w_{\theta,\epsilon}}{\partial x_i}\phi\ dx=\displaystyle\int_{\Omega} \dfrac{1}{\alpha}w_{\theta}^{\frac{1}{\alpha}-1}\dfrac{\partial w_{\theta}}{\partial x_i}\phi\ dx.
\]

\noindent Therefore, making $\epsilon\searrow 0$ in \eqref{eqlimit2}, we get for any $\phi\in C^{\infty}_c(\Omega)$ and any $i\in\overline{1,N}$ that:

\begin{equation}
		\int_{\Omega} w_{\theta}^{\frac{1}{\alpha}}\dfrac{\partial\phi}{\partial x_i}\ dx=-\int_{\Omega}\dfrac{1}{\alpha}w_{\theta}^{\frac{1}{\alpha}-1}\dfrac{\partial w_{\theta}}{\partial x_i}\phi\ dx.
\end{equation} 

\noindent This proves that the weak gradient of $w_{\theta}^{\frac{1}{\alpha}}$ exists and $\nabla w_{\theta}^{\frac{1}{\alpha}}=\dfrac{1}{\alpha}w_{\theta}^{\frac{1}{\alpha}-1}\nabla w_{\theta}$. Note also that from Lemma \ref{lemmatrunc} \textbf{(7)} we may write for any $i\in\overline{1,N}$:

\begin{equation}
	\dfrac{1}{\alpha}w_{\theta}^{\frac{1}{\alpha}-1}\dfrac{\partial w_{\theta}}{\partial x_i}=\theta\underbrace{\left (\dfrac{w_1}{w_{\theta}} \right )^{1-\frac{1}{\alpha}}}_{\in L^{\infty}(\Omega)}\underbrace{\dfrac{w_1^{\frac{1}{\alpha}}}{\partial x_i}}_{\in L^{p(x)}(\Omega)}+(1-\theta)\underbrace{\left (\dfrac{w_2}{w_{\theta}} \right )^{1-\frac{1}{\alpha}}}_{\in L^{\infty}(\Omega)}\underbrace{\dfrac{w_2^{\frac{1}{\alpha}}}{\partial x_i}}_{\in L^{p(x)}(\Omega)}\in L^{p(x)}(\Omega).
\end{equation}

\noindent We have used here \textbf{(4)} and the fact that $w_1^{\frac{1}{\alpha}},w_2^{\frac{1}{\alpha}}\in W^{1,p(x)}(\Omega)$. Thus $\exists\ \nabla w_{\theta}^{\frac{1}{\alpha}}=\dfrac{1}{\alpha}w_{\theta}^{\frac{1}{\alpha}-1}\nabla w_{\theta}\in L^{p(x)}(\Omega)^N$.
	
	\noindent Hence $w_{\theta}>0$ a.e. on $\Omega$, $w_{\theta}\in L^{\frac{p(x)}{\alpha}}(\Omega)$ and $\sqrt[\alpha]{w_{\theta}}\in W^{1,p(x)}(\Omega)$. This shows that for each $\theta\in (-\theta_0,1+\theta_0)$ we have $w_{\theta}\in\overset{\bullet}{W}$. 
	\bigskip
	
	\noindent\textbf{(7)} Introduce for each $\epsilon\in (0,1)$: $\gamma_{\theta,\epsilon}=\dfrac{1}{\alpha}\cdot\dfrac{w_{1,\epsilon}-w_{2,\epsilon}}{\sqrt[\alpha]{w_{\theta,\epsilon}^{\alpha-1}}}$. First, we prove that $\gamma_{\theta,\epsilon}\in W^{1,p(x)}(\Omega)$ for any $\epsilon\in (0,1)$. 
	
	\noindent From \textbf{(3)} we have that $\dfrac{\epsilon}{2}\leq w_{\theta,\epsilon}\leq \dfrac{3}{2\epsilon}$. Since $\dfrac{1}{\alpha}-1<0$ we get that $\left(\dfrac{2\epsilon}{3}\right )^{1-\frac{1}{\alpha}}\leq w_{\theta,\epsilon}^{\frac{1}{\alpha}-1}\leq \left (\dfrac{2}{\epsilon}\right )^{1-\frac{1}{\alpha}}$. Thus $w_{\theta,\epsilon}^{\frac{1}{\alpha}-1}\in L^{\infty}(\Omega)$. We also know, from Lemma \ref{lemmatrunc} \textbf{(3)}, that $w_{\theta,\epsilon}=\theta w_{1,\epsilon}+(1-\theta)w_{2,\epsilon}\in W^{1,p(x)}(\Omega)$. The function $\mathbb{R}\ni s\mapsto\begin{cases}\left (\frac{\epsilon}{2} \right )^{\frac{1}{\alpha}-1}, & s<\frac{\epsilon}{2}\\[3mm] s^{\frac{1}{\alpha}-1}, & s\in\left[\frac{\epsilon}{2},\frac{3}{2\epsilon}\right]\\[3mm] \left (\frac{3}{2\epsilon} \right )^{\frac{1}{\alpha}-1}, & s>\frac{3}{2\epsilon} \end{cases}$ is $\left(1-\dfrac{1}{\alpha}\right)\left(\dfrac{2}{\epsilon}\right)^{2-\frac{1}{\alpha}}$--Lipschitz. Composing it with $w_{\theta,\epsilon}\in\left[\frac{\epsilon}{2},\frac{3}{2\epsilon} \right ]$ gives $w_{\theta,\epsilon}^{\frac{1}{\alpha}-1}$. So, from the \textit{chain rule}, we have that $w_{\theta,\epsilon}^{\frac{1}{\alpha}-1}\in W^{1,p(x)}(\Omega)$ and $\nabla w_{\theta,\epsilon}^{\frac{1}{\alpha}-1}=\left (\frac{1}{\alpha}-1 \right )w_{\theta,\epsilon}^{\frac{1}{\alpha}-2}\nabla w_{\theta,\epsilon}$. Thus $\dfrac{1}{\alpha}(w_{1,\epsilon}-w_{2,\epsilon}),w_{\theta,\epsilon}^{\frac{1}{\alpha}-1}\in W^{1,p(x)}(\Omega)\cap L^{\infty}(\Omega)$. From the \textit{product rule} we obtain that $\gamma_{\theta,\epsilon}\in W^{1,p(x)}(\Omega)$ and:
	
	\begin{equation}\label{eqnablagamma}
		\nabla\gamma_{\theta,\epsilon}=\dfrac{1}{\alpha}\dfrac{\nabla w_{1,\epsilon}-\nabla w_{2,\epsilon}}{w_{\theta,\epsilon}^{1-\frac{1}{\alpha}}}+\dfrac{1}{\alpha}\left(\dfrac{1}{\alpha}-1\right)\dfrac{w_{1,\epsilon}-w_{2,\epsilon}}{w_{\theta,\epsilon}^{2-\frac{1}{\alpha}}}\nabla w_{\theta,\epsilon}.
	\end{equation}
	
	\noindent We have from the \textit{Lebesgue dominated convergence theorem} that $\lim\limits_{\epsilon\searrow 0} \gamma_{\theta,\epsilon}=\gamma_\theta$ in $L^{p(x)}(\Omega)$ because\footnote{See the \textit{Lebesgue dominated convergence theorem} for variable exponents in \cite[Lemma 3.2.8, page 77]{Hasto}.}:
	
	\begin{enumerate}
		\item[$\bullet$] $\lim\limits_{\epsilon\searrow 0} \gamma_{\theta,\epsilon}(x)=\lim\limits_{\epsilon\searrow 0} \dfrac{1}{\alpha}\cdot\dfrac{w_{1,\epsilon}(x)-w_{2,\epsilon}(x)}{\sqrt[\alpha]{w_{\theta,\epsilon}(x)^{\alpha-1}}}=\dfrac{1}{\alpha}\cdot\dfrac{w_{1}(x)-w_{2}(x)}{\sqrt[\alpha]{w_{\theta}(x)^{\alpha-1}}}=\gamma_{\theta}(x)$ for a.a. $x\in\Omega$. This follows from Lemma \ref{lemmatrunc} \textbf{(1)}.
		
		\item[$\bullet$] $|\gamma_{\theta,\epsilon}|=\dfrac{1}{\alpha}|w_{1,\epsilon}-w_{2,\epsilon}|w_{\theta,\epsilon}^{\frac{1}{\alpha}-1}=\dfrac{1}{\alpha}\left |\dfrac{w_{1,\epsilon}}{w_{\theta,\epsilon}}-\dfrac{w_{2,\epsilon}}{w_{\theta,\epsilon}} \right |w_{\theta,\epsilon}^{\frac{1}{\alpha}}\stackrel{\textbf{(5)}}{\leq}\dfrac{1}{\alpha}\left (2M-\dfrac{2}{3M}\right )w_{\theta,\epsilon}^{\frac{1}{\alpha}}$
		
		\begin{align*}
	\text{(triangle ineq.)}\ \ \ \ \ \ \ \ 	&\leq \dfrac{1}{\alpha}\left (2M-\dfrac{2}{3M}\right )\left[w_{\theta}^{\frac{1}{\alpha}}+\left | w_{\theta,\epsilon}^{\frac{1}{\alpha}}-w_{\theta}^{\frac{1}{\alpha}} \right | \right ]\\
\text{(Lemma \ref{lemrunu} \textbf{(1)})}\ \ \ \ \ \ \ 	&\leq \dfrac{1}{\alpha}\left (2M-\dfrac{2}{3M}\right )\left[w_{\theta}^{\frac{1}{\alpha}}+\left | w_{\theta,\epsilon}-w_{\theta} \right |^{\frac{1}{\alpha}} \right ]\\
&=\dfrac{1}{\alpha}\left (2M-\dfrac{2}{3M}\right )\left[w_{\theta}^{\frac{1}{\alpha}}+\left |\theta(w_{1,\epsilon}-w_1)+(1-\theta)(w_{2,\epsilon}-w_2)\right |^{\frac{1}{\alpha}}\right ]\\
\text{(triangle ineq.)}\ \ \ \ \ \ \ \ &\dfrac{1}{\alpha}\left (2M-\dfrac{2}{3M}\right )\left[w_{\theta}^{\frac{1}{\alpha}}+\left (|\theta|\cdot \big |w_{1,\epsilon}-w_1\big |+|1-\theta|\cdot\big |w_{2,\epsilon}-w_2\big |\right )^{\frac{1}{\alpha}}\right ]\\
\text{(Lemma \ref{lemrunu} \textbf{(2)})}\ \ \ \ \ \ \ &\leq \dfrac{1}{\alpha}\left (2M-\dfrac{2}{3M}\right )\left[w_{\theta}^{\frac{1}{\alpha}}+|\theta|^{\frac{1}{\alpha}}\cdot \big |w_{1,\epsilon}-w_1\big |^{\frac{1}{\alpha}}+|1-\theta|^{\frac{1}{\alpha}}\cdot\big |w_{2,\epsilon}-w_2\big |^{\frac{1}{\alpha}}\right ]\\
\text{(Lemma \ref{lemmatrunc} \textbf{(9)})}\ \ \ \ \ \ \ & \leq\dfrac{1}{\alpha}\left (2M-\dfrac{2}{3M}\right )\left[w_{\theta}^{\frac{1}{\alpha}}+|\theta|^{\frac{1}{\alpha}}\cdot \big |w_1-1\big |^{\frac{1}{\alpha}}+|1-\theta|^{\frac{1}{\alpha}}\cdot\big |w_2-1\big |^{\frac{1}{\alpha}}\right ]\\
\text{(triangle ineq.)}\ \ \ \ \ \ \ \ &\leq \dfrac{1}{\alpha}\left (2M-\dfrac{2}{3M}\right )\left[w_{\theta}^{\frac{1}{\alpha}}+|\theta|^{\frac{1}{\alpha}}\cdot \big (w_1+1\big )^{\frac{1}{\alpha}}+|1-\theta|^{\frac{1}{\alpha}}\cdot\big (w_2+1\big )^{\frac{1}{\alpha}}\right ]\\
\text{(Lemma \ref{lemrunu} \textbf{(2)})}\ \ \ \ \ \ \ &\leq \dfrac{1}{\alpha}\left (2M-\dfrac{2}{3M}\right )\left[w_{\theta}^{\frac{1}{\alpha}}+|\theta|^{\frac{1}{\alpha}}\cdot \big (w_1^{\frac{1}{\alpha}}+1\big )+|1-\theta|^{\frac{1}{\alpha}}\cdot\big (w_2^{\frac{1}{\alpha}}+1\big )\right ]\in L^{p(x)}(\Omega).
		\end{align*}
	\end{enumerate}
	
	\noindent As a consequence, $\gamma_{\theta}\in L^{p(x)}(\Omega)$.
	
	\noindent Next, we show in the same manner that $\lim\limits_{\epsilon\searrow 0} \nabla \gamma_{\theta,\epsilon}=\dfrac{1}{\alpha}\dfrac{\nabla w_{1}-w_{2}}{w_{\theta}^{1-\frac{1}{\alpha}}}+\dfrac{1}{\alpha}\left(\dfrac{1}{\alpha}-1\right)\dfrac{w_{1}-w_{2}}{w_{\theta}^{2-\frac{1}{\alpha}}}\nabla w_{\theta}$ in $L^{p(x)}(\Omega)$.
	
	\noindent Indeed, for each $i\in\overline{1,N}$, taking into account \eqref{eqnablagamma}, we have that:
	
	\begin{enumerate}
		\item[$\bullet$] $\lim\limits_{\epsilon\searrow 0} \dfrac{\partial \gamma_{\theta,\epsilon}}{\partial x_i}(x)=\dfrac{1}{\alpha}\dfrac{\frac{\partial w_{1}}{\partial x_i}(x)-\frac{\partial w_{2}}{\partial x_i}(x)}{w_{\theta}(x)^{1-\frac{1}{\alpha}}}+\dfrac{1}{\alpha}\left(\dfrac{1}{\alpha}-1\right)\dfrac{w_{1}(x)-w_{2}(x)}{w_{\theta}(x)^{2-\frac{1}{\alpha}}}\dfrac{\partial w_{\theta}}{\partial x_i}(x)$, from Lemma \ref{lemmatrunc} \textbf{(1)},\textbf{(8)}.
		
		\item[$\bullet$] $\left |\dfrac{\partial\gamma_{\theta,\epsilon}}{\partial x_i} \right |\leq \dfrac{1}{\alpha}\left |\dfrac{\frac{\partial w_{1,\epsilon}}{\partial x_i}-\frac{\partial w_{2,\epsilon}}{\partial x_i}}{w_{\theta,\epsilon}^{1-\frac{1}{\alpha}}} \right |+\dfrac{1}{\alpha}\left (1-\dfrac{1}{\alpha}\right )\cdot \left |\dfrac{w_{1,\epsilon}-w_{2,\epsilon}}{w_{\theta,\epsilon}^{2-\frac{1}{\alpha}}} \right |\cdot\left |\dfrac{\partial w_{\theta,\epsilon}}{\partial x_i} \right |$.
		
		\noindent For the first term, using Lemma \ref{lemmatrunc} \textbf{(3)} and \textbf{(7)} we have that:
		
		\begin{align*}
			\dfrac{1}{\alpha}\left |\dfrac{\frac{\partial w_{1,\epsilon}}{\partial x_i}-\frac{\partial w_{2,\epsilon}}{\partial x_i}}{w_{\theta,\epsilon}^{1-\frac{1}{\alpha}}} \right |&=\dfrac{1}{\alpha}\left |\dfrac{\frac{\partial w_{1}}{\partial x_i}\chi_{I_{\epsilon}}(w_1)-\frac{\partial w_{2}}{\partial x_i}\chi_{I_{\epsilon}}(w_2)}{w_{\theta,\epsilon}^{1-\frac{1}{\alpha}}} \right |=\left |\dfrac{w_1^{1-{\frac{1}{\alpha}}}\frac{\partial w_{1}^{\frac{1}{\alpha}}}{\partial x_i}\chi_{I_{\epsilon}}(w_1)-w_2^{1-{\frac{1}{\alpha}}}\frac{\partial w_{2}^{\frac{1}{\alpha}}}{\partial x_i}\chi_{I_{\epsilon}}(w_2)}{w_{\theta,\epsilon}^{1-\frac{1}{\alpha}}} \right | \\
			&=\left |\left (\dfrac{w_{1,\epsilon}}{w_{\theta,\epsilon}} \right )^{1-{\frac{1}{\alpha}}}\frac{\partial w_{1}^{\frac{1}{\alpha}}}{\partial x_i}-\left (\dfrac{w_{2,\epsilon}}{w_{\theta,\epsilon}} \right )^{1-{\frac{1}{\alpha}}}\frac{\partial w_{2}^{\frac{1}{\alpha}}}{\partial x_i}\right |\leq \left (\dfrac{w_{1,\epsilon}}{w_{\theta,\epsilon}} \right )^{1-{\frac{1}{\alpha}}}\left |\frac{\partial w_{1}^{\frac{1}{\alpha}}}{\partial x_i}\right |+\left (\dfrac{w_{2,\epsilon}}{w_{\theta,\epsilon}} \right )^{1-{\frac{1}{\alpha}}}\left |\frac{\partial w_{2}^{\frac{1}{\alpha}}}{\partial x_i}\right |\\
\textbf{(5)}\ \ \ \ \ \ \ &\leq (2M)^{1-\frac{1}{\alpha}}\left (\left |\frac{\partial w_{1}^{\frac{1}{\alpha}}}{\partial x_i}\right | +\left |\frac{\partial w_{2}^{\frac{1}{\alpha}}}{\partial x_i}\right |\right )\in L^{p(x)}(\Omega).
		\end{align*}
		
		\noindent For the second term we have that:
		
		\begin{align*}
			\dfrac{1}{\alpha}\left (1-\dfrac{1}{\alpha}\right )\cdot& \left |\dfrac{w_{1,\epsilon}-w_{2,\epsilon}}{w_{\theta,\epsilon}^{2-\frac{1}{\alpha}}} \right |\cdot\left |\dfrac{\partial w_{\theta,\epsilon}}{\partial x_i} \right |=\dfrac{\alpha-1}{\alpha^2}\left |\dfrac{w_{1,\epsilon}}{w_{\theta,\epsilon}}-\dfrac{w_{2,\epsilon}}{w_{\theta,\epsilon}} \right | \dfrac{1}{w_{\theta,\epsilon}^{1-\frac{1}{\alpha}}}\left |\theta\dfrac{\partial w_{1,\epsilon}}{\partial x_i}+(1-\theta)\dfrac{\partial w_{2,\epsilon}}{\partial x_i} \right | \\
\textbf{(5)}\ \ \ \  &\leq\dfrac{\alpha-1}{\alpha^2}\left (2M-\dfrac{2}{3M}\right )\dfrac{1}{w_{\theta,\epsilon}^{1-\frac{1}{\alpha}}}\left |\theta\dfrac{\partial w_{1,\epsilon}}{\partial x_i}+(1-\theta)\dfrac{\partial w_{2,\epsilon}}{\partial x_i} \right |  \\
\text{(Lemma \ref{lemmatrunc} \textbf{(3)})}	\ \ \ \  	&=\dfrac{\alpha-1}{\alpha}\left (2M-\dfrac{2}{3M}\right )\dfrac{1}{w_{\theta,\epsilon}^{1-\frac{1}{\alpha}}}\left |\theta w_1^{1-\frac{1}{\alpha}}\chi_{I_{\epsilon}}(w_1)\dfrac{\partial w_1^{\frac{1}{\alpha}}}{\partial x_i}+(1-\theta)w_2^{1-\frac{1}{\alpha}}\chi_{I_{\epsilon}}(w_2)\dfrac{\partial w_2^{\frac{1}{\alpha}}}{\partial x_i} \right |\\
&=\dfrac{\alpha-1}{\alpha}\left (2M-\dfrac{2}{3M}\right )\left |\theta\left(\dfrac{w_{1,\epsilon}}{w_{\theta,\epsilon}} \right )^{1-\frac{1}{\alpha}}\dfrac{\partial w_1^{\frac{1}{\alpha}}}{\partial x_i}+(1-\theta)\left(\dfrac{w_{2,\epsilon}}{w_{\theta,\epsilon}} \right )^{1-\frac{1}{\alpha}}\dfrac{\partial w_2^{\frac{1}{\alpha}}}{\partial x_i} \right |\\
&\leq \dfrac{\alpha-1}{\alpha}\left (2M-\dfrac{2}{3M}\right )\left [|\theta|\left(\dfrac{w_{1,\epsilon}}{w_{\theta,\epsilon}} \right )^{1-\frac{1}{\alpha}}\left |\dfrac{\partial w_1^{\frac{1}{\alpha}}}{\partial x_i} \right |+|1-\theta|\left(\dfrac{w_{2,\epsilon}}{w_{\theta,\epsilon}} \right )^{1-\frac{1}{\alpha}}\left |\dfrac{\partial w_2^{\frac{1}{\alpha}}}{\partial x_i} \right | \right ]\\
\textbf{(5)}\ \ \ \ &\leq \dfrac{\alpha-1}{\alpha}\left (2M-\dfrac{2}{3M}\right )\cdot (2M)^{1-\frac{1}{\alpha}}\left (|\theta|\cdot \left |\dfrac{\partial w_1^{\frac{1}{\alpha}}}{\partial x_i} \right | +|1-\theta|\cdot \left |\dfrac{\partial w_2^{\frac{1}{\alpha}}}{\partial x_i} \right | \right )\in L^{p(x)}(\Omega).
		\end{align*}

		\noindent Combining the last two relations we get that:
		
		\[
		\left |\dfrac{\partial\gamma_{\theta,\epsilon}}{\partial x_i} \right |\leq(2M)^{1-\frac{1}{\alpha}}\left \{\left |\frac{\partial w_{1}^{\frac{1}{\alpha}}}{\partial x_i}\right |\left[1+\dfrac{|\theta|(\alpha-1)\left (2M-\frac{2}{3M}\right )}{\alpha}\right] +\left |\frac{\partial w_{2}^{\frac{1}{\alpha}}}{\partial x_i}\right |\left[1+\dfrac{|1-\theta|(\alpha-1)\left (2M-\frac{2}{3M}\right )}{\alpha}\right]\right \}\in L^{p(x)}(\Omega). 
		\]
		
	\end{enumerate}
	
	\noindent From \textit{Lebesgue dominated convergence theorem} it follows that: $\lim\limits_{\epsilon\searrow 0}\dfrac{\partial \gamma_{\theta,\epsilon}}{\partial x_i}=\dfrac{1}{\alpha}\dfrac{\frac{\partial w_{1}}{\partial x_i}-\frac{\partial w_{2}}{\partial x_i}}{w_{\theta}^{1-\frac{1}{\alpha}}}+\dfrac{1}{\alpha}\left(\dfrac{1}{\alpha}-1\right)\dfrac{w_{1}-w_{2}}{w_{\theta}^{2-\frac{1}{\alpha}}}\dfrac{\partial w_{\theta}}{\partial x_i}:=\tilde{\gamma}_i$ in $L^{p(x)}(\Omega)$, for each $i\in\overline{1,N}$. In particular $\tilde{\gamma}_i\in L^{p(x)}(\Omega)$.
	
	\noindent Now, equation \eqref{eqnablagamma} can be equivalently rewritten as:
	
	\begin{equation}\label{eqnablagammalimit}
		\int_{\Omega} \gamma_{\theta,\epsilon}\dfrac{\partial\phi}{\partial x_i}\ dx=-\int_{\Omega}\dfrac{\partial\gamma_{\theta,\epsilon}}{\partial x_i}\phi\ dx,\ \forall\ \phi\in C^{\infty}_c(\Omega)\ \text{and any}\ i\in\overline{1,N}.
	\end{equation}
	
	\noindent Taking into consideration that:
	
	\begin{enumerate}
		\item[$\bullet$] $\left |\displaystyle\int_{\Omega}\gamma_{\theta,\epsilon}\dfrac{\partial\phi}{\partial x_i}\ dx-\displaystyle\int_{\Omega}\gamma_{\theta}\dfrac{\partial\phi}{\partial x_i}\ dx\right |\leq\left\Vert\dfrac{\partial\phi}{\partial x_i}\right\Vert_{L^{\infty}(\Omega)}\displaystyle\int_{\Omega}\big |\gamma_{\theta,\epsilon}-\gamma_{\theta}\big |\ dx=\left\Vert\dfrac{\partial\phi}{\partial x_i}\right\Vert_{L^{\infty}(\Omega)}\Vert\gamma_{\theta,\epsilon}-\gamma_{\theta}\Vert_{L^1(\Omega)}\stackrel{\epsilon\searrow 0}{\longrightarrow} 0$.
		
		\medskip
		
		\item[$\bullet$] $\left |\displaystyle\int_{\Omega}\dfrac{\partial\gamma_{\theta,\epsilon}}{\partial x_i}\phi\ dx-\displaystyle\int_{\Omega}\tilde{\gamma}_i\phi\ dx \right |\leq \Vert\phi\Vert_{L^{\infty}(\Omega)}\displaystyle\int_{\Omega}\left |\dfrac{\partial\gamma_{\theta,\epsilon}}{\partial x_i}-\tilde{\gamma}_i \right |\ dx=\Vert\phi\Vert_{L^{\infty}(\Omega)}\left\Vert\dfrac{\partial\gamma_{\theta,\epsilon}}{\partial x_i}-\tilde{\gamma}_i\right\Vert_{L^1(\Omega)}\stackrel{\epsilon\searrow 0}{\longrightarrow} 0$,
	\end{enumerate}
	
	\noindent we can make $\epsilon\searrow 0$ in \eqref{eqnablagammalimit} to obtain that: $\displaystyle\int_{\Omega} \gamma_{\theta}\dfrac{\partial\phi}{\partial x_i}\ dx=-\displaystyle\int_{\Omega}\tilde{\gamma}_i\phi\ dx$ for any $\phi\in C^{\infty}_c(\Omega)$ and each $i\in\overline{1,N}$.
	
	\noindent This proves that the weak gradient of $\gamma_{\theta}\in L^{p(x)}(\Omega)$ exists and:

\[
\nabla\gamma_{\theta}=\tilde{\gamma}:=(\tilde{\gamma}_1,\tilde{\gamma}_2,\hdots,\tilde{\gamma}_N)=\dfrac{1}{\alpha}\dfrac{\nabla w_{1}-\nabla w_{2}}{w_{\theta}^{1-\frac{1}{\alpha}}}+\dfrac{1}{\alpha}\left(\dfrac{1}{\alpha}-1\right)\dfrac{w_{1}-w_{2}}{w_{\theta}^{2-\frac{1}{\alpha}}}\nabla w_{\theta}\in L^{p(x)}(\Omega)^N.
\]
	
	\bigskip
	
	\noindent\textbf{(8)} Fix some $\theta\in (-\theta_0,1+\theta_0)$ and take any sequence of real numbers $(\epsilon_n)_{n\geq 1}$ that converges to $0$ such that $\theta+\epsilon_n\in (-\theta_0,1+\theta_0)$ for each $n\geq 1$. Denote $g_n:=\dfrac{\Gamma(\theta+\epsilon_n)-\Gamma(\theta)}{\epsilon_n}\in W^{1,p(x)}(\Omega)$ for each whole number $n\geq 1$. We need to show that: $\lim\limits_{n\to\infty} \Vert g_n-\gamma_{\theta}\Vert_{W^{1,p(x)}(\Omega)}=0\ \Longleftrightarrow\ \lim\limits_{n\to\infty} \Vert g_n-\gamma_{\theta}\Vert_{L^{p(x)}(\Omega)}=0$ and for each $i\in\overline{1,N}$: $\lim\limits_{n\to\infty} \left\Vert \dfrac{\partial g_n}{\partial x_i}-\dfrac{\partial\gamma_{\theta}}{\partial x_i}\right \Vert_{L^{p(x)}(\Omega)}=0$.
	
	\noindent For an arbitrarily fixed $x\in\Omega$ set $h: (-\theta_0-\theta,1+\theta_0-\theta )\to(0,\infty),\ h(\epsilon)=w_{\theta+\epsilon}^{\frac{1}{\alpha}}=\big [(\theta+\epsilon)w_1(x)+(1-\theta-\epsilon)w_2(x) \big ]^{\frac{1}{\alpha}}$, so that $\theta+\epsilon\in(-\theta_0,1+\theta_0 )$. Clearly $h\in C^1\big ((-\theta_0-\theta,1+\theta_0-\theta ) \big )$ and $h'(\epsilon)=\dfrac{1}{\alpha}\dfrac{w_1(x)-w_2(x)}{w_{\theta+\epsilon}(x)^{1-\frac{1}{\alpha}}}, \ \forall\ \epsilon\in (-\theta_0-\theta,1+\theta_0-\theta )$. Note that $0\in(-\theta_0-\theta,1+\theta_0-\theta )$.
	
	\noindent Remark that: $g_n(x)=\dfrac{w_{\theta+\epsilon_n}^{\frac{1}{\alpha}}-w_{\theta}^{\frac{1}{\alpha}}}{\epsilon_n}=\dfrac{h(\epsilon_n)-h(0)}{\epsilon_n}=h'(\tilde{\epsilon}_n)=\dfrac{1}{\alpha}\dfrac{w_1(x)-w_2(x)}{w_{\theta+\tilde{\epsilon}_n}(x)^{1-\frac{1}{\alpha}}}$, from the \textit{mean value theorem}, for some $\tilde{\epsilon}_n$ between $\epsilon_n$ and $0$. Therefore $\lim\limits_{n\to\infty} g_n(x)=\lim\limits_{n\to\infty} h'(\tilde{\epsilon}_n)=h'(0)=\dfrac{1}{\alpha}\dfrac{w_1(x)-w_2(x)}{w_{\theta}(x)^{1-\frac{1}{\alpha}}}=\gamma_{\theta}(x)$. This relation holds for a.a. $x\in\Omega$, so we have proved that $g_n\to\gamma_{\theta}$ pointwise a.e. on $\Omega$. Moreover:
	
	\begin{align*}
		|g_n|&=\dfrac{1}{\alpha}\cdot\left|\dfrac{w_1}{w_{\theta+\tilde{\epsilon}_n}}-\dfrac{w_2}{w_{\theta+\tilde{\epsilon}_n}}\right |w_{\theta+\tilde{\epsilon}_n}^{\frac{1}{\alpha}}\stackrel{\textbf{(4)}}{\leq}\dfrac{1}{\alpha}\left(2M-\dfrac{2}{3M}\right)w_{\theta+\tilde{\epsilon}_n}^{\frac{1}{\alpha}}\\
		&\stackrel{\textbf{(1)}}{\leq}\dfrac{1}{\alpha}\left(2M-\dfrac{2}{3M}\right)\left(\dfrac{3}{2}\right)^{\frac{1}{\alpha}}\max\{w_1^{\frac{1}{\alpha}},w_2^{\frac{1}{\alpha}}\}\in L^{p(x)}(\Omega).
	\end{align*}
	
	\noindent Applying \textit{Lebesgue dominated convergence theorem} we get that $\lim\limits_{n\to\infty} g_n=\gamma_{\theta}$ in $L^{p(x)}(\Omega)$. 
	
	\noindent The last thing we need to show is that $\lim\limits_{\epsilon\searrow 0} \dfrac{\partial g_n}{\partial x_i}=\dfrac{\partial\gamma_{\theta}}{\partial x_i}$ in $L^{p(x)}(\Omega)$ for each $i\in\overline{1,N}$. We have that:
	
	\begin{equation}
		\dfrac{\partial g_n}{\partial x_i}=\dfrac{1}{\epsilon_n}\left (\dfrac{\partial w_{\theta+\epsilon_n}^{\frac{1}{\alpha}}}{\partial x_i}-\dfrac{\partial w_{\theta}^{\frac{1}{\alpha}}}{\partial x_i} \right)\stackrel{\textbf{(6)}}{=}\dfrac{1}{\alpha\epsilon_n}\left (w_{\theta+\epsilon_n}^{\frac{1}{\alpha}-1}\dfrac{\partial w_{\theta+\epsilon_n}}{\partial x_i}-w_{\theta}^{\frac{1}{\alpha}-1}\dfrac{\partial w_{\theta}}{\partial x_i}\right ).
	\end{equation}
	
	\noindent Define, for a randomly fixed $x\in\Omega$ the function $\tilde{h}:(-\theta_0-\theta,1+\theta_0-\theta)\to\mathbb{R}$,

	\begin{equation} \tilde{h}(\epsilon):=\dfrac{1}{\alpha}w_{\theta+\epsilon}(x)^{\frac{1}{\alpha}-1}\dfrac{\partial w_{\theta+\epsilon}}{\partial x_i}(x)=\dfrac{1}{\alpha}\left [(\theta+\epsilon)w_1(x)+(1-\theta-\epsilon)w_2(x) \right ]^{\frac{1}{\alpha}-1}\left[(\theta+\epsilon)\dfrac{\partial w_1}{\partial x_i}+(1-\theta-\epsilon)\dfrac{\partial w_2}{\partial x_i}\right].
	\end{equation}
	
	\noindent It is easy to see that $\tilde{h}\in C^1\big((-\theta_0-\theta,1+\theta_0-\theta) \big )$ and for each $\epsilon\in(-\theta_0-\theta,1+\theta_0-\theta)$:
	
	\begin{equation}
		\tilde{h}'(\epsilon)=\dfrac{1}{\alpha}w_{\theta+\epsilon}(x)^{\frac{1}{\alpha}-1}\left(\dfrac{\partial w_1}{\partial x_i}(x)-\dfrac{\partial w_2}{\partial x_i}(x)\right )+\dfrac{1}{\alpha}\left(\dfrac{1}{\alpha}-1\right)w_{\theta+\epsilon}(x)^{\frac{1}{\alpha}-2}\big (w_1(x)-w_2(x) \big )\dfrac{\partial w_{\theta+\epsilon}}{\partial x_i}(x).
	\end{equation}
	
	\noindent Remark that:
	
	\[
	\tilde{h}'(0)=\dfrac{1}{\alpha}\dfrac{\frac{\partial w_1}{\partial x_i}(x)-\frac{\partial w_2}{\partial x_i}(x)}{w_{\theta}(x)^{1-\frac{1}{\alpha}}}+\dfrac{1}{\alpha}\left(\dfrac{1}{\alpha}-1\right)\dfrac{w_1(x)-w_2(x)}{w_{\theta}(x)^{2-\frac{1}{\alpha}}}\dfrac{\partial w_{\theta}}{\partial x_i}(x)=\dfrac{\partial\gamma_\theta}{\partial x_i}(x).
	\]
	
	\noindent Henceforth, from the \textit{mean value theorem}, there is some $\tilde{\epsilon}_n$\footnote{For the sake of simplicity we have used for $\tilde{h}$ the same notation $\tilde{\epsilon}_n$ as for the function $h$ above. They are certainly not the same.} between $\epsilon_n$ and $0$ such that:
	
	\begin{equation}
		\dfrac{\partial g_n}{\partial x_i}(x)=\dfrac{\tilde{h}(\epsilon_n)-\tilde{h}(0)}{\epsilon_n}=\tilde{h}'(\tilde{\epsilon}_n)\ \Longrightarrow\ \lim\limits_{n\to\infty}\dfrac{\partial g_n}{\partial x_i}(x)=\lim\limits_{n\to\infty} \tilde{h}'(\tilde{\epsilon}_n)=\tilde{h}'(0)=\dfrac{\partial\gamma_{\theta}}{\partial x_i}(x).
	\end{equation}
	
	\noindent This relation holds for a.a. $x\in\Omega$, so $\dfrac{\partial g_n}{\partial x_i}\longrightarrow \dfrac{\partial\gamma_\theta}{\partial x_i}$ pointwise a.e. on $\Omega$. In addition,

	\[
	\left |\dfrac{\partial g_n}{\partial x_i} \right |\leq \dfrac{1}{\alpha}\left |\dfrac{\frac{\partial w_{1}}{\partial x_i}-\frac{\partial w_{2}}{\partial x_i}}{w_{\theta+\tilde{\epsilon}_n}^{1-\frac{1}{\alpha}}} \right |+\dfrac{1}{\alpha}\left (1-\dfrac{1}{\alpha}\right )\cdot \left |\dfrac{w_{1}-w_{2}}{w_{\theta+\tilde{\epsilon}_n}^{2-\frac{1}{\alpha}}} \right |\cdot\left |\dfrac{\partial w_{\theta+\tilde{\epsilon}_n}}{\partial x_i} \right |.
	\]
	
	\noindent For the first term, using Lemma \ref{lemmatrunc} \textbf{(7)} we have that:
	
	\begin{align*}
		\dfrac{1}{\alpha}\left |\dfrac{\frac{\partial w_{1}}{\partial x_i}-\frac{\partial w_{2}}{\partial x_i}}{w_{\theta+\tilde{\epsilon}_n}^{1-\frac{1}{\alpha}}} \right |&=\left |\dfrac{w_1^{1-{\frac{1}{\alpha}}}\frac{\partial w_{1}^{\frac{1}{\alpha}}}{\partial x_i}-w_2^{1-{\frac{1}{\alpha}}}\frac{\partial w_{2}^{\frac{1}{\alpha}}}{\partial x_i}}{w_{\theta+\tilde{\epsilon}_n}^{1-\frac{1}{\alpha}}} \right |=\left |\left (\dfrac{w_{1}}{w_{\theta+\tilde{\epsilon}_n}} \right )^{1-{\frac{1}{\alpha}}}\frac{\partial w_{1}^{\frac{1}{\alpha}}}{\partial x_i}-\left (\dfrac{w_{2}}{w_{\theta+\tilde{\epsilon}_n}} \right )^{1-{\frac{1}{\alpha}}}\frac{\partial w_{2}^{\frac{1}{\alpha}}}{\partial x_i}\right |\\
		&\leq \left (\dfrac{w_{1}}{w_{\theta+\tilde{\epsilon}_n}} \right )^{1-{\frac{1}{\alpha}}}\left |\frac{\partial w_{1}^{\frac{1}{\alpha}}}{\partial x_i}\right |+\left (\dfrac{w_{2}}{w_{\theta+\tilde{\epsilon}_n}} \right )^{1-{\frac{1}{\alpha}}}\left |\frac{\partial w_{2}^{\frac{1}{\alpha}}}{\partial x_i}\right |\\
		\textbf{(4)}\ \ \ \ \ \ \ &\leq (2M)^{1-\frac{1}{\alpha}}\left (\left |\frac{\partial w_{1}^{\frac{1}{\alpha}}}{\partial x_i}\right | +\left |\frac{\partial w_{2}^{\frac{1}{\alpha}}}{\partial x_i}\right |\right )\in L^{p(x)}(\Omega).
	\end{align*}
	
	\noindent For the second term we have that:
	
	\begin{align*}
		\dfrac{1}{\alpha}\left (1-\dfrac{1}{\alpha}\right )\cdot& \left |\dfrac{w_{1}-w_{2}}{w_{\theta+\tilde{\epsilon}_n}^{2-\frac{1}{\alpha}}} \right |\cdot\left |\dfrac{\partial w_{\theta+\tilde{\epsilon}_n}}{\partial x_i} \right |=\dfrac{\alpha-1}{\alpha^2}\left |\dfrac{w_{1}}{w_{\theta+\tilde{\epsilon}_n}}-\dfrac{w_{2}}{w_{\theta+\tilde{\epsilon}_n}} \right | \dfrac{1}{w_{\theta+\tilde{\epsilon}_n}^{1-\frac{1}{\alpha}}}\left |(\theta+\tilde{\epsilon}_n)\dfrac{\partial w_{1}}{\partial x_i}+(1-\theta-\tilde{\epsilon}_n)\dfrac{\partial w_{2}}{\partial x_i} \right | \\
		\textbf{(4)}\ \ \ \  &\leq\dfrac{\alpha-1}{\alpha^2}\left (2M-\dfrac{2}{3M}\right )\dfrac{1}{w_{\theta+\tilde{\epsilon}_n}^{1-\frac{1}{\alpha}}}\left |(\theta+\tilde{\epsilon}_n)\dfrac{\partial w_{1}}{\partial x_i}+(1-\theta-\tilde{\epsilon}_n)\dfrac{\partial w_{2}}{\partial x_i} \right | \\
		\text{(Lemma \ref{lemmatrunc} \textbf{(7)})}	\ \ \ \  	&=\dfrac{\alpha-1}{\alpha}\left (2M-\dfrac{2}{3M}\right )\dfrac{1}{w_{\theta+\tilde{\epsilon}_n}^{1-\frac{1}{\alpha}}}\left |(\theta+\tilde{\epsilon}_n) w_1^{1-\frac{1}{\alpha}}\dfrac{\partial w_1^{\frac{1}{\alpha}}}{\partial x_i}+(1-\theta-\tilde{\epsilon}_n)w_2^{1-\frac{1}{\alpha}}\dfrac{\partial w_2^{\frac{1}{\alpha}}}{\partial x_i} \right |\\
		&=\dfrac{\alpha-1}{\alpha}\left (2M-\dfrac{2}{3M}\right )\left |(\theta+\tilde{\epsilon}_n)\left(\dfrac{w_{1}}{w_{\theta+\tilde{\epsilon}_n}} \right )^{1-\frac{1}{\alpha}}\dfrac{\partial w_1^{\frac{1}{\alpha}}}{\partial x_i}+(1-\theta-\tilde{\epsilon}_n)\left(\dfrac{w_{2}}{w_{\theta+\tilde{\epsilon}_n}} \right )^{1-\frac{1}{\alpha}}\dfrac{\partial w_2^{\frac{1}{\alpha}}}{\partial x_i} \right |\\
		&\leq \dfrac{\alpha-1}{\alpha}\left (2M-\dfrac{2}{3M}\right )\left [|\theta+\tilde{\epsilon}_n|\left(\dfrac{w_{1}}{w_{\theta+\tilde{\epsilon}_n}} \right )^{1-\frac{1}{\alpha}}\left |\dfrac{\partial w_1^{\frac{1}{\alpha}}}{\partial x_i} \right |+|1-\theta-\tilde{\epsilon}_n|\left(\dfrac{w_{2}}{w_{\theta+\tilde{\epsilon}_n}} \right )^{1-\frac{1}{\alpha}}\left |\dfrac{\partial w_2^{\frac{1}{\alpha}}}{\partial x_i} \right | \right ]\\
		\textbf{(4)}\ \ \ \ &\leq \dfrac{\alpha-1}{\alpha}\left (2M-\dfrac{2}{3M}\right )\cdot (2M)^{1-\frac{1}{\alpha}}\left (|\theta+\tilde{\epsilon}_n|\cdot \left |\dfrac{\partial w_1^{\frac{1}{\alpha}}}{\partial x_i} \right | +|1-\theta-\tilde{\epsilon}_n|\cdot \left |\dfrac{\partial w_2^{\frac{1}{\alpha}}}{\partial x_i} \right | \right )\\
		&\leq\dfrac{(\alpha-1)(1+\theta_0)}{\alpha}\left (2M-\dfrac{2}{3M}\right )\cdot (2M)^{1-\frac{1}{\alpha}}\left ( \left |\dfrac{\partial w_1^{\frac{1}{\alpha}}}{\partial x_i} \right | + \left |\dfrac{\partial w_2^{\frac{1}{\alpha}}}{\partial x_i} \right | \right ) \in L^{p(x)}(\Omega).
	\end{align*}

	\noindent From the last two relations we get that for every integer $n\geq 1$:
	
	\[
	\left |\dfrac{\partial g_n}{\partial x_i} \right |\leq (2M)^{1-\frac{1}{\alpha}}\left [1+\dfrac{(\alpha-1)(1+\theta_0)}{\alpha}\left (2M-\dfrac{2}{3M}\right )\right ]\cdot \left ( \left |\dfrac{\partial w_1^{\frac{1}{\alpha}}}{\partial x_i} \right | + \left |\dfrac{\partial w_2^{\frac{1}{\alpha}}}{\partial x_i} \right | \right )\in L^{p(x)}(\Omega). 
	\]
	
	\noindent We can conclude, from the \textit{Lebesgue dominated convergence theorem}, that $\lim\limits_{n\to\infty}\dfrac{\partial g_n}{\partial x_i}=\dfrac{\partial\gamma_{\theta}}{\partial x_i}$ in $L^{p(x)}(\Omega)$. This is true for each $i\in\overline{1,N}$, so we can deduce finally that $\lim\limits_{n\to\infty} g_n=\gamma_{\theta}$ in $W^{1,p(x)}(\Omega)$, as claimed.

\end{proof}

\begin{definition}
	We define the functional $J:\overset{\bullet}{W}\to\mathbb{R}$ by:
	
	\begin{equation}
		J(w)=\mathcal{J}(\sqrt[\alpha]{w})=\int_{\Omega} A(x,\nabla \sqrt[\alpha]{w})\ dx -\int_{\Omega} \overline{F}(x,\sqrt[\alpha]{w})\ dx,\ \forall\ w\in \overset{\bullet}{W}.
	\end{equation}
\end{definition}

\begin{theorem}\label{thmbeta}
	Consider that hypothesis \textnormal{\textbf{(H7)}} holds for $r=\alpha$ and let $w_1,w_2\in\overset{\bullet}{W}$ with $\dfrac{w_1}{w_2},\ \dfrac{w_2}{w_1}\in L^{\infty}(\Omega)$. For $\theta_0$ defined as in Proposition \ref{propwtheta}, the real function $\beta:(-\theta_0,1+\theta_0)\to\mathbb{R}$ given by:
	
	\begin{equation}
		\beta(\theta)=J(w_{\theta})=J\big (\theta w_1+(1-\theta) w_2 \big ),\ \forall\ \theta\in (-\theta_0,1+\theta_0)
	\end{equation}
	
	\noindent is differentiable on $(-\theta_0,\theta_0+1)$ and the following formula holds for all $\theta\in (-\theta_0,1+\theta_0)$:
	
	\begin{equation}
		\beta^\prime(\theta)=\dfrac{1}{\alpha}\int_{\Omega} \mathbf{a}\big (x,\nabla \sqrt[\alpha]{w_{\theta}}\big )\cdot\nabla\left (\dfrac{w_1-w_2}{\sqrt[\alpha]{w_{\theta}^{\alpha-1}}}\right )\ dx-\dfrac{1}{\alpha}\int_{\Omega} \overline{f}(x,\sqrt[\alpha]{w_\theta})\dfrac{w_1-w_2}{\sqrt[\alpha]{w_{\theta}^{\alpha-1}}}\ dx.
	\end{equation}
	
	\noindent Moreover, $\beta$ is a convex function and if \textnormal{\textbf{(H13')}} holds and $w_1\not\equiv w_2$ then $\beta$ is \textbf{strictly convex}.
\end{theorem}

\begin{proof} We have that $\beta=\mathcal{J}\circ\Gamma$. Since both $\mathcal{J}:W^{1,p(x)}(\Omega)\to\mathbb{R}$ and $\Gamma:(-\theta_0,1+\theta_0)\to W^{1,p(x)}(\Omega)$ are differentiable (see Proposition \ref{propc1j} and Proposition \ref{propwtheta} \textbf{(8)}), it follows that their composition $\beta$ is also differentiable\footnote{See  the differentiation of a composition of mappings (chain rule) from \cite[page 63]{zor}.}  and $ \forall\ \theta\in (-\theta_0,1+\theta_0)$:
	
	\begin{equation}
		\beta'(\theta)=\langle\mathcal{J}'(\Gamma(\theta)),\Gamma'(\theta)\rangle=\dfrac{1}{\alpha}\int_{\Omega} \mathbf{a}\big (x,\nabla \sqrt[\alpha]{w_{\theta}}\big )\cdot\nabla\left (\dfrac{w_1-w_2}{\sqrt[\alpha]{w_{\theta}^{\alpha-1}}}\right )\ dx-\dfrac{1}{\alpha}\int_{\Omega} \overline{f}(x,\sqrt[\alpha]{w_\theta})\dfrac{w_1-w_2}{\sqrt[\alpha]{w_{\theta}^{\alpha-1}}}\ dx.
	\end{equation}
	
\noindent In order to show that $\beta$ is a convex real function we will prove that $\beta':(-\theta_0,1+\theta_0)\to\mathbb{R}$ is an increasing function. Set any $1+\theta_0>\theta_1>\theta_2>-\theta_0$. For this aim we will use the proved \textbf{D\'{i}az-Saa inequality} (Theorem \ref{thmdiazsaa}) for $U_{\theta_1}:=\sqrt[\alpha]{w_{\theta_1}},U_{\theta_2}:=\sqrt[\alpha]{w_{\theta_2}}\in W^{1,p(x)}(\Omega)$ with $\dfrac{U_{\theta_1}}{U_{\theta_2}}=\left (\dfrac{w_{\theta_1}}{w_{\theta_2}}\right )^{\frac{1}{\alpha}},\ \dfrac{U_{\theta_2}}{U_{\theta_1}}=\left (\dfrac{w_{\theta_2}}{w_{\theta_1}}\right )^{\frac{1}{\alpha}}\in L^{\infty}(\Omega)$ and $U_{\theta_1},U_{\theta_2}>0$ a.e. on $\Omega$. Note that:

\begin{equation}
	U_{\theta_1}^{\alpha}-U_{\theta_2}^\alpha=w_{\theta_1}-w_{\theta_2}=\big (\theta_1w_1+(1-\theta_1)w_2\big )-\big (\theta_2w_1+(1-\theta_2)w_2 \big )=(\theta_1-\theta_2)(w_1-w_2).
\end{equation}

\noindent We have that:

\begin{equation}
	\begin{cases} \beta'(\theta_1)&=\dfrac{1}{\alpha}\displaystyle\int_{\Omega} \mathbf{a}\big (x,\nabla \sqrt[\alpha]{w_{\theta_1}}\big )\cdot\nabla\left (\dfrac{w_1-w_2}{\sqrt[\alpha]{w_{\theta_1}^{\alpha-1}}}\right )\ dx-\dfrac{1}{\alpha}\int_{\Omega} \overline{f}(x,\sqrt[\alpha]{w_{\theta_1}})\dfrac{w_1-w_2}{\sqrt[\alpha]{w_{\theta_1}^{\alpha-1}}}\ dx \\[3mm]
		&=\dfrac{1}{\alpha(\theta_1-\theta_2)}\displaystyle\int_{\Omega} \mathbf{a}\big (x,\nabla U_{\theta_1}\big )\cdot\nabla\left (U_{\theta_1}-\dfrac{U_{\theta_2}^\alpha}{U_{\theta_1}^{\alpha-1}}\right )\ dx+\dfrac{1}{\alpha}\int_{\Omega} -\dfrac{\overline{f}(x,\sqrt[\alpha]{w_{\theta_1}})}{{\sqrt[\alpha]{w_{\theta_1}^{\alpha-1}}}}(w_1-w_2)\ dx\\[3mm]
	 \beta'(\theta_2)&=\dfrac{1}{\alpha}\displaystyle\int_{\Omega} \mathbf{a}\big (x,\nabla \sqrt[\alpha]{w_{\theta_2}}\big )\cdot\nabla\left (\dfrac{w_1-w_2}{\sqrt[\alpha]{w_{\theta_2}^{\alpha-1}}}\right )\ dx-\dfrac{1}{\alpha}\int_{\Omega} \overline{f}(x,\sqrt[\alpha]{w_{\theta_2}})\dfrac{w_1-w_2}{\sqrt[\alpha]{w_{\theta_2}^{\alpha-1}}}\ dx\\[3mm]
	 &=\dfrac{1}{\alpha(\theta_1-\theta_2)}\displaystyle\int_{\Omega} \mathbf{a}\big (x,\nabla U_{\theta_2}\big )\cdot\nabla\left (\dfrac{U_{\theta_1}^\alpha}{U_{\theta_2}^{\alpha-1}}-U_{\theta_2}\right )\ dx+\dfrac{1}{\alpha}\int_{\Omega} -\dfrac{\overline{f}(x,\sqrt[\alpha]{w_{\theta_2}})}{{\sqrt[\alpha]{w_{\theta_2}^{\alpha-1}}}}(w_1-w_2)\ dx\end{cases}
\end{equation}

\noindent Henceforth:

\begin{align*}
	\beta'(\theta_1)-\beta'(\theta_2)&=\dfrac{1}{\alpha(\theta_1-\theta_2)}\left[\displaystyle\int_{\Omega} \mathbf{a}\big (x,\nabla U_{\theta_1}\big )\cdot\nabla\left (U_{\theta_1}-\dfrac{U_{\theta_2}^\alpha}{U_{\theta_1}^{\alpha-1}}\right )\ dx-\displaystyle\int_{\Omega} \mathbf{a}\big (x,\nabla U_{\theta_2}\big )\cdot\nabla\left (\dfrac{U_{\theta_1}^\alpha}{U_{\theta_2}^{\alpha-1}}-U_{\theta_2}\right )\ dx\right ]\\
	&\ \ +\dfrac{1}{\alpha(\theta_1-\theta_2)}\displaystyle\int_{\Omega}\left [-\dfrac{\overline{f}(x,\sqrt[\alpha]{w_{\theta_1}})}{{\sqrt[\alpha]{w_{\theta_1}^{\alpha-1}}}}-\left (-\dfrac{\overline{f}(x,\sqrt[\alpha]{w_{\theta_2}})}{{\sqrt[\alpha]{w_{\theta_2}^{\alpha-1}}}} \right) \right] (w_{\theta_1}-w_{\theta_2})\ dx\\
(\text{D\'{i}az-Saa, \textbf{(H13)}})\ \ 	&\geq 0.
\end{align*}

\noindent So far we have proved that $\beta'$ is increasing, which means that $\beta$ is convex. If \textbf{(H13')} holds then we will show that $\beta'$ is strictly increasing. Indeed, if we set $\Omega'=\{x\in\Omega\ |\ w_1(x)\neq w_2(x)\}=\{x\in\Omega\ |\ w_{\theta_1}(x)\neq w_{\theta_2}(x)\}$, then using the fact that $w_1\not\equiv w_2$ we deduce that $|\Omega'|>0$, and therefore from Proposition \ref{propf} \textbf{(4)}:

\begin{equation}
	\beta'(\theta_1)-\beta'(\theta_2)\geq \dfrac{1}{\alpha(\theta_1-\theta_2)}\displaystyle\int_{\Omega'}\left [-\dfrac{\overline{f}(x,\sqrt[\alpha]{w_{\theta_1}})}{{\sqrt[\alpha]{w_{\theta_1}^{\alpha-1}}}}-\left (-\dfrac{\overline{f}(x,\sqrt[\alpha]{w_{\theta_2}})}{{\sqrt[\alpha]{w_{\theta_2}^{\alpha-1}}}} \right) \right] (w_{\theta_1}-w_{\theta_2})\ dx>0.
\end{equation}

\noindent So $\beta$ is strictly convex and the theorem is proved.
\end{proof}

\begin{corollary} The following inequalities hold for the functional $J:\overset{\bullet}{W}\to\mathbb{R}$:
	
\begin{enumerate}
	\item[(i)] If \textnormal{\textbf{(H13)}} holds then for any $\theta\in [0,1]$: $J(\theta w_1+(1-\theta)w_2)\leq\theta J(w_1)+(1-\theta) J(w_2)$.
	
	\item[(ii)] If \textnormal{\textbf{(H13')}} holds and $w_1\not\equiv w_2$ then for any $\theta\in (0,1)$: $J(\theta w_1+(1-\theta)w_2)<\theta J(w_1)+(1-\theta) J(w_2)$.
\end{enumerate}  

\noindent These inequalities can be used only if $w_1,w_2\in\overset{\bullet}{W}$ and $\dfrac{w_2}{w_1},\dfrac{w_1}{w_2}\in L^{\infty}(\Omega)$.
\end{corollary}

\begin{proof} (i) The conclusion follows if we just use the convexity of $\beta$. More precisely, for any $\theta\in [0,1]$:
	
\begin{equation}
	J(\theta w_1+(1-\theta)w_2)=J(w_{\theta})=\beta(\theta)=\beta(\theta\cdot 1+(1-\theta)\cdot 0)\leq \theta\beta(1)+(1-\theta)\beta(0)=\theta J(w_1)+(1-\theta)J(w_2).
\end{equation}

\noindent (ii) Analogously, for any $\theta\in (0,1)$. from the strict convexity of $\beta$ we have that:

\begin{equation}
	J(\theta w_1+(1-\theta)w_2)=\beta(\theta)=\beta(\theta\cdot 1+(1-\theta)\cdot 0)< \theta\beta(1)+(1-\theta)\beta(0)=\theta J(w_1)+(1-\theta)J(w_2).
\end{equation}

\noindent 
	
\end{proof}

\begin{theorem}[\textbf{Uniqueness}]\label{thmuniq} Consider that hypothesis \textnormal{\textbf{(H7)}} holds for $r=\alpha$.
	
\begin{enumerate}
	\item[\textbf{(1)}] If $U_1,U_2$ are two distinct weak solutions ($U_1\not\equiv U_2$) of \eqref{eqmphase} that are \textbf{strongly positive}, i.e. $\underset{\Omega}{\operatorname{ess\ inf}}\ U_1$, $ \underset{\Omega}{\operatorname{ess\ inf}}\ U_2>0$, then there is a constant $\lambda>0$ such that:
	
	\medskip
	
\begin{enumerate}
	\item[$\bullet$] $U_2(x)=\lambda U_1(x)$, for a.a. $x\in\Omega$,
	
	\medskip
	
	\item[$\bullet$] $\Phi\big (x,\lambda|\nabla U_1(x)|\big )=\lambda^{\alpha-1}\Phi\big (x,|\nabla U_1(x)|\big ),\ \text{for a.a.}\ x\in\Omega$,
	
	\medskip
	
	\item[$\bullet$] $f\big (x,\lambda U_1(x)\big )=\lambda^{\alpha-1} f\big (x,U_1(x)\big )$, for a.a. $x\in\Omega$.
\end{enumerate} 

\bigskip
	
\item[\textbf{(2)}] If hypothesis \textnormal{\textbf{(H13')}} holds then problem \eqref{eqmphase} has at most one weak solution $U$ that is \textbf{strongly positive}.

\bigskip

\item[\textbf{(3)}] If hypothesis \textnormal{\textbf{(H7')}} holds for $r=\alpha$ then problem \eqref{eqmphase} has at most one weak solution $U$ that is \textbf{strongly positive}.
\end{enumerate}
\end{theorem}

\begin{proof}\textbf{(1)} Suppose that the problem \eqref{eqmphase} has two distinct solutions $U_1,U_2\in W^{1,p(x)}(\Omega)$, $U_1\not\equiv U_2$, with $\epsilon\leq U_1,U_2\leq 1$ for some small enough $\epsilon\in (0,1)$. Denote $w_1:=U_1^{\alpha}$ and $w_2:=U_2^{\alpha}$. Remark that $w_1,w_2>0$ a.e. on $\Omega$, $\dfrac{w_1}{w_2},\dfrac{w_2}{w_1}\in \left [\epsilon^{\alpha},\dfrac{1}{\epsilon^{\alpha}} \right ]\Longrightarrow\ \dfrac{w_2}{w_1},\dfrac{w_1}{w_2}\in L^{\infty}(\Omega)$, and $\sqrt[\alpha]{w_1}=U_1,\sqrt[\alpha]{w_2}=U_2\in W^{1,p(x)}(\Omega)$. Thus $w_1,w_2\in\overset{\bullet}{W}$.
	
\noindent Because $U_1,U_2$ are weak solutions of \eqref{eqmphase}, from Proposition \ref{propc1j} we have that for each $\phi\in W^{1,p(x)}(\Omega)$:

\begin{equation}\label{eq78}
\begin{cases}\langle \mathcal{J}'(U_1),\phi\rangle=\displaystyle\int_{\Omega} \mathbf{a}(x,\nabla U_1(x))\cdot\nabla\phi(x)\ dx-\int_{\Omega} \overline{f}(x,U_1(x))\phi(x)\ dx=0\\[3mm]\langle \mathcal{J}'(U_2),\phi\rangle=\displaystyle\int_{\Omega} \mathbf{a}(x,\nabla U_2(x))\cdot\nabla\phi(x)\ dx-\int_{\Omega} \overline{f}(x,U_2(x))\phi(x)\ dx=0 \end{cases}
\end{equation}

\noindent In what follows we will use the notations introduced in Proposition \ref{propwtheta} and Theorem \ref{thmbeta}. Choosing for $\theta=1\in (-\theta_0,1+\theta_0)$ and $\theta=0\in (-\theta_0,1+\theta_0)$, $\phi:=\gamma_{1}=\dfrac{1}{\alpha}\cdot\dfrac{w_1-w_2}{\sqrt[\alpha]{w_{1}^{\alpha-1}}}\in W^{1,p(x)}(\Omega)$ and then $\phi:=\gamma_{0}=\dfrac{1}{\alpha}\cdot\dfrac{w_1-w_2}{\sqrt[\alpha]{w_{2}^{\alpha-1}}}\in W^{1,p(x)}(\Omega)$ (from Proposition \ref{propwtheta} \textbf{(7)}), the two relation in \eqref{eq78} become:

\begin{equation}
	\begin{cases}\beta'(1)=\displaystyle\dfrac{1}{\alpha}\int_{\Omega} \mathbf{a}(x,\nabla \sqrt[\alpha]{w_1})\cdot\nabla\left ( \dfrac{w_1-w_2}{\sqrt[\alpha]{w_{1}^{\alpha-1}}}\right)\ dx-\dfrac{1}{\alpha}\int_{\Omega} \dfrac{\overline{f}(x,\sqrt[\alpha]{w_1})}{\sqrt[\alpha]{w_{1}^{\alpha-1}}}(w_1-w_2)\ dx=0\\[3mm] \beta'(0)=\displaystyle\dfrac{1}{\alpha}\int_{\Omega} \mathbf{a}(x,\nabla \sqrt[\alpha]{w_2})\cdot\nabla\left ( \dfrac{w_1-w_2}{\sqrt[\alpha]{w_{2}^{\alpha-1}}}\right)\ dx-\dfrac{1}{\alpha}\int_{\Omega} \dfrac{\overline{f}(x,\sqrt[\alpha]{w_2})}{\sqrt[\alpha]{w_{2}^{\alpha-1}}}(w_1-w_2)\ dx=0\end{cases}
\end{equation}

\noindent Up to this point we have that $\beta'(1)=\beta'(0)=0$. Knowing from Theorem \ref{thmbeta} that $\beta'$ is an increasing function we get that $\beta'(\theta)=0$ for every $\theta\in [0,1]$. Thus there is some constant $\beta_0$ such that $\beta(\theta)=\beta_0$ for any $\theta\in [0,1]$. 

\noindent Also from $\beta'(1)=\beta'(0)$ we have from the proof of Theorem \ref{thmbeta} that equality occurs in \textbf{D\'{i}az-Saa inequality} for the functions $U_1,\ U_2$ for $r=\alpha>1$. It is very important to check the condition $\dfrac{\nabla U_1}{U_1},\dfrac{\nabla U_2}{U_2}\in L^1_{\textnormal{loc}}(\Omega)$. Indeed, since $\dfrac{1}{U_1},\dfrac{1}{U_2}\in\left[1,\dfrac{1}{\epsilon}\right]$, we have that both of them are in $L^{\infty}(\Omega)$. So $\dfrac{\nabla U_1}{U_1}=\underbrace{\dfrac{1}{U_1}}_{\in L^{\infty}(\Omega)}\underbrace{\nabla U_1}_{\in L^{p(x)}(\Omega)^N}\in L^{p(x)}(\Omega)^N\subset L^1_{\text{loc}}(\Omega)^N$. In the same manner $\dfrac{\nabla U_2}{U_2}\in L^1_{\text{loc}}(\Omega)^N$. So we are in position to apply Theorem \ref{thmequalitydiaz} (ii). Therefore there is some $\lambda>0$ with $U_2=\lambda U_1$ a.e. on $\Omega$ and $\lambda\neq 1$ (because $U_1\not\equiv U_2$). Moreover the following relation holds:

\begin{equation}
	\Phi(x,\lambda|\nabla U_1|)=\lambda^{\alpha-1}\Phi(x,|\nabla U_1|),\ \text{for a.a.}\ x\in\Omega.
\end{equation}

\noindent This is not all, since we also have this equality at our disposal:

\begin{equation}\label{eqcontrad}
	\displaystyle\int_{\Omega}\left [-\dfrac{\overline{f}(x,\sqrt[\alpha]{w_{1}})}{{\sqrt[\alpha]{w_{1}^{\alpha-1}}}}-\left (-\dfrac{\overline{f}(x,\sqrt[\alpha]{w_{2}})}{{\sqrt[\alpha]{w_{2}^{\alpha-1}}}} \right) \right] (w_{1}-w_{2})\ dx=0.
\end{equation}

\noindent We have that $w_2=\lambda^{\alpha}w_1$ a.e. on $\Omega$, so $w_1\neq w_2$ a.e. on $\Omega$. The above equality and Proposition \ref{propf} \textbf{(3)} allow us to write that a.e. on $\Omega$ the following equation:

\begin{equation}
	\dfrac{\overline{f}(x,\sqrt[\alpha]{w_{1}})}{{\sqrt[\alpha]{w_{1}^{\alpha-1}}}}=\dfrac{\overline{f}(x,\sqrt[\alpha]{w_{2}})}{{\sqrt[\alpha]{w_{2}^{\alpha-1}}}}\ \Longleftrightarrow\ \overline{f}(x,\sqrt[\alpha]{w_1})=\dfrac{\overline{f}(x,\lambda\sqrt[\alpha]{w_1})}{\lambda^{\alpha-1}}\ \Longleftrightarrow\ f(x,\lambda U_1)=\lambda^{\alpha-1} f(x,U_1).
\end{equation}

\noindent Here we mention that $\lambda U_1=U_2\leq 1$, so we can use $f$ instead of $\overline{f}$.

\bigskip

\noindent Part \textbf{(2)} and Part \textbf{(3)} of this theorem both follow immediately from the last part of Theorem \ref{thmequalitydiaz} and from equation \eqref{eqcontrad} respectively. They will provide the desired contradictions. Hence in both cases we obtain $U_1\equiv U_2$ and the uniqueness of the strongly positive weak solution follows.
\end{proof}

\section{Example 1: A multiple-phase problem}

\noindent Consider the following general multiple phase problem:

\begin{equation}\tag{$MP$}\label{eqmphase2}
	\begin{cases}-\Delta_{p_1(x),w_1}U-\Delta_{p_2(x),w_2}U-\hdots-\Delta_{p_{\ell}(x),w_{\ell}}U=f\big (x,U(x)\big ), & x\in\Omega\\[3mm] \dfrac{\partial U}{\partial\nu}=0, & x\in \partial\Omega\\[3mm] 0\leq U(x)\leq 1, & x\in\Omega\end{cases}
\end{equation}

\noindent where $-\Delta_{p_{k}(x),w_k}U:=\operatorname{div}\big (w_k(x)|\nabla U|^{p_k(x)-2}\nabla U \big )$ for each $k\in\overline{1,\ell}$ and:

\begin{enumerate}
\item[$\bullet$] $p_1,p_2,\dots, p_{\ell}:\overline{\Omega}\to (1,\infty),\ \ell\geq 1$ are continuous variable exponents such that $p:\overline{\Omega}\to (1,\infty),\ p(x)=\max\{p_1(x),p_2(x),\hdots,p_{\ell}(x)\}$ satisfies the inequality $p^-:=\displaystyle\min_{x\in\overline{\Omega}} p(x)>\dfrac{2N}{N+2}$, or

\item[$\bullet$] $p_1,p_2,\dots, p_{\ell}:\overline{\Omega}\to (1,\infty),\ \ell\geq 1$ are log-H\"{o}lder continuous variable exponents such that $p:\overline{\Omega}\to (1,\infty),\ p(x)=\max\{p_1(x),p_2(x),\hdots,p_{\ell}(x)\}$ satisfies the inequality $p^-:=\displaystyle\min_{x\in\overline{\Omega}} p(x)\geq\dfrac{2N}{N+2}$.\footnote{It it important to point out that if $p_1,p_2,\hdots, p_{\ell}$ are all continuous\textbackslash log-H\"{o}lder continuous\textbackslash Lipschitz continuous, the so is $p$. To see why it is so, just note that the following elementary inequality holds: $|p(x)-p(y)|=|\max\{p_1(x),p_2(x),\hdots,p_{\ell}(x)\}-\max\{p_1(y),p_2(y),\hdots, p_{\ell}(y)\}|\leq \max\{|p_1(x)-p_1(y)|,|p_2(x)-p_2(y)|,\hdots, |p_{\ell}(x)-p_{\ell}(y)|\}$ for each $x,y\in\overline{\Omega}$.}

\item[$\bullet$] The weights $w_1,w_2,\hdots,w_{\ell}\in L^{\infty}(\Omega)$ with the property that there is some constant $\omega>0$ such that $\displaystyle\min_{j\in\overline{1,\ell}}\underset{x\in\Omega}{\operatorname{ess\ inf}}\ w_j(x)\geq\omega$.

\end{enumerate}

\noindent In this case:

\begin{equation}
\mathbf{a}:\overline{\Omega}\times\mathbb{R}^N\to\mathbb{R}^N,\ \mathbf{a}(x,\xi)=\begin{cases} w_1(x)|\xi|^{p_1(x)-2}\xi+w_2(x)|\xi|^{p_2(x)-2}\xi+\hdots+w_\ell(x)|\xi|^{p_\ell(x)-2}\xi, & \xi\neq \mathbf{0}\\ \mathbf{0}, & \xi=\mathbf{0}\end{cases}.
\end{equation}

\noindent From \cite[Proposition 3.1]{max1} we know that $\mathbf{a}$ satifies all of our hypotheses. We only need to check that \textbf{(H7')} holds for $r=\alpha\in (1,p_{-}
)$, where $p_{-}:=\displaystyle\min_{x\in\overline{\Omega}}\min_{k\in\overline{1,\ell}} p_k(x)>1$. Indeed, for a.a. $x\in\Omega$, the function:

\begin{align*}
	[0,\infty)\ni s\mapsto \dfrac{\Phi(x,s)}{s^{\alpha-1}}&=\dfrac{w_1(x)s^{p_1(x)-1}+w_2(x)s^{p_2(x)-1}+\hdots+w_\ell(x)s^{p_\ell(x)-1}}{s^{\alpha-1}}\\
	&=w_1(x)s^{p_1(x)-\alpha}+w_2(x)s^{p_2(x)-\alpha}+\hdots+w_\ell(x)s^{p_\ell(x)-\alpha},
\end{align*} 

\noindent is \textbf{strictly increasing} being a sum of strictly increasing functions. It is essential that $p_{k}(x)\geq p_{-}>\alpha$ for every $k\in\overline{1,\ell}$. Therefore applying Theorem \ref{thmuniq} in this context we obtain the following corollary:

\begin{corollary}\label{cormultiple}
	Problem \eqref{eqmphase2} has at most one strong positive solution.
\end{corollary}

\section{Example 2: An application from image processing}

\noindent This example is used in image processing tasks, like restoration, and it is inspired and adapted from \cite[Page 1]{Abo}. Fix any $\varepsilon,\delta >0$ and $p:\overline{\Omega}\to (1,\infty)$ be a variable exponent that satisfies \textbf{(H2)} and $p^->\alpha>1$. Consider the following problem:

\begin{equation}\tag{$IP$}\label{eqimage}
	\begin{cases}-\operatorname{div}\left (\begin{cases}|\nabla U|^{p(x)-2}\ln^{\delta}(1+|\nabla U|)\nabla U, & 0\leq |\nabla U|\leq\varepsilon\\[3mm] \varepsilon^{p(x)-\alpha}|\nabla U|^{\alpha-2}\ln^{\delta}(1+|\nabla U|)\nabla U, & |\nabla U|\geq\varepsilon \end{cases}\right )=f\big (x,U(x)\big ), & x\in\Omega\\[3mm] \dfrac{\partial U}{\partial\nu}=0, & x\in \partial\Omega\\[3mm] 0\leq U(x)\leq 1, & x\in\Omega\end{cases}
\end{equation}

\noindent In this example we have that:

\begin{enumerate}
	\item[$\bullet$] $\Psi:\Omega\times (0,\infty)\to\mathbb{R},\ \Psi(x,s)=\begin{cases} s^{p(x)-2}\ln^{\delta}(1+s), & 0<s\leq\varepsilon\\[3mm] \varepsilon^{p(x)-\alpha}s^{\alpha-2}\ln^{\delta}(1+s), & s\geq\varepsilon\end{cases}$.
	
	\noindent Clearly $\Psi(\cdot,s)$ is measurable for all $s\in (0,\infty)$. In order to show that for some arbitrarily fixed $x\in\Omega$ we have that $\Psi(x,\cdot)$ is locally absolutely continuous, consider some compact interval $[a,b]\subset (0,\infty)$. We have that:
	
	\begin{equation}
		\dfrac{\partial\Psi(x,s)}{\partial s}=\begin{cases}s^{p(x)-3}\ln^{\delta-1}(1+s)\left [\big (p(x)-2\big )\ln(1+s)+\dfrac{\delta s}{1+s} \right ], & 0<s<\varepsilon \\[3mm] \varepsilon^{p(x)-\alpha}s^{\alpha-3}\ln^{\delta-1}(1+s)\left [(\alpha-2)\ln(1+s)+\dfrac{\delta s}{1+s} \right ]   & s>\varepsilon\end{cases}
	\end{equation}
	
	\noindent Note that for $s=\varepsilon$ we have finite lateral derivative for $\Psi(x,\cdot)$. This proves in particular that $\dfrac{\partial\Psi(x,s)}{\partial s}$ is integrable on any $[a,b]\subset (0,\infty)$ (having just one point of discontinuity of the first kind: i.e. a jump discontinuity in $\varepsilon$). Thus, if $[a,b]\subset (0,\varepsilon]$ or $[a,b]\subset [\varepsilon,\infty)$ then $\dfrac{\partial\Psi(x,s)}{\partial s}\in C^1([a,b])$ which means that $\dfrac{\partial\Psi(x,s)}{\partial s}$ is Lipschitz on $[a,b]$ and in particular it follows that $\dfrac{\partial\Psi(x,s)}{\partial s}$ is absolutely continuous as needed. If $\varepsilon\in (a,b)$, it is easy to see that:
	
	\begin{align*}
		&\ \ \ \ \ \ \ \Psi(x,a)+\int_{a}^s\dfrac{\partial\Psi(x,\tau)}{\partial s}\ d\tau=\\
		&=\begin{cases}\Psi(x,a)+\Psi(x,s)-\Psi(x,a), & a\leq s\leq\varepsilon\\[3mm] \Psi(x,a)+\displaystyle\int_{a}^{\varepsilon} \dfrac{\partial\Psi(x,\tau)}{\partial s}\ d\tau+\displaystyle\int_{\varepsilon}^s\dfrac{\partial\Psi(x,\tau)}{\partial s}\ d\tau=\Psi(x,a)+\Psi(x,\varepsilon)-\Psi(x,a)+\Psi(x,s)-\Psi(x,\varepsilon),& \varepsilon<s\leq b \end{cases}\\[3mm]
		&=\Psi(x,s),\ \forall\ s\in [a,b].
	\end{align*}
	
	\noindent Using Theorem 13 from Chapter 6 in \cite{royden} we deduce that $\dfrac{\partial\Psi(x,s)}{\partial s}$ is absolutely continuous on $[a,b]$ being the indefinite integral of a Lebesgue integrable function. We have proved that $\Psi$ satisfies \textbf{(H3)}. Since $p(x),\alpha>1$ we have that $\lim\limits_{s\to 0^+} \Psi(x,s)s=0$ for a.a. $x\in\Omega$, i.e. \textbf{(H4)} holds. Also \textbf{(H5)} is satisfied from the same reason, since it is obvious to see that for a.a. $x\in\Omega$ the function $(0,\infty)\ni s\mapsto \Psi(x,s)s=\begin{cases}s^{p(x)-1}\ln^{\delta}(1+s), & 0\leq s\leq\varepsilon\\[3mm] \varepsilon^{p(x)-\alpha}s^{\alpha-1}\ln^{\delta}(1+s), & |s|\geq\varepsilon \end{cases}$ is strictly increasing, being continuous at $s=\varepsilon$ and strictly increasing on both branches. 
	
	\item[$\bullet$] $\Phi:\overline{\Omega}\times\mathbb{R}\to [0,\infty),\ \Phi(x,s)=\begin{cases}|s|^{p(x)-1}\ln^{\delta}(1+|s|), & 0\leq|s|\leq\varepsilon\\[3mm] \varepsilon^{p(x)-\alpha}|s|^{\alpha-1}\ln^{\delta}(1+|s|), & |s|\geq\varepsilon \end{cases}$.

	\noindent It is time to check \textbf{(H6)}:
	
	\begin{equation}
		|\Phi(x,s)|=|s|^{p(x)-1}\cdot\begin{cases}\ln^{\delta}(1+|s|), & 0\leq|s|\leq\varepsilon\\[3mm] \left [\dfrac{\ln(1+|s|)}{\left (\frac{|s|}{\varepsilon}\right )^{\frac{p(x)-\alpha}{\delta}}}\right ]^{\delta}, & |s|\geq\varepsilon\end{cases}\leq |s|^{p(x)-1}\max\{C^{\delta},\ln^{\delta}(1+\varepsilon)\},\ \forall\ s\in \mathbb{R}.
	\end{equation}
	
	\noindent We explain now who is this $C$. Remark that $\lim\limits_{|s|\to\infty} \dfrac{\ln(1+|s|)}{\left (\frac{|s|}{\varepsilon}\right )^{\frac{p(x)-\alpha}{\delta}}}=0$ because $p(x)\geq p^{-}>\alpha$. So there is some $\tilde{\varepsilon}\geq\varepsilon$ such that $\dfrac{\ln(1+|s|)}{\left (\frac{|s|}{\varepsilon}\right )^{\frac{p(x)-\alpha}{\delta}}}\leq 1$ for any $s\in\mathbb{R}$ with $|s|\geq\tilde{\varepsilon}$. But the function $[\varepsilon,\tilde{\varepsilon}]\ni\tau\mapsto \dfrac{\ln(1+\tau)}{\left (\frac{\tau}{\varepsilon}\right )^{\frac{p(x)-\alpha}{\delta}}}$ is bounded (being continuous on a compact interval) by some $C$ that will be assumed, without loss of generality, to be greater than or equal to $1$. This allows us to assert that $\dfrac{\ln(1+|s|)}{\left (\frac{|s|}{\varepsilon}\right )^{\frac{p(x)-\alpha}{\delta}}}\leq C,\ \forall\ s\in\mathbb{R}$ with $|s|\geq\varepsilon$. So \textbf{(H6)} is verified.
	
	\noindent Next we show that \textbf{(H7')} holds for $r=\alpha$. Indeed, for any $s\in [0,\infty)$ we have that:
	
	\begin{equation}
		\dfrac{\Phi(x,s)}{s^{\alpha-1}}=\begin{cases} s^{p(x)-\alpha}\ln^{\delta}(1+s), & 0\leq s\leq\varepsilon\\[3mm] \varepsilon^{p(x)-\alpha}\ln^{\delta}(1+s), & s\geq\varepsilon \end{cases}.
	\end{equation}
	
	\noindent Since this function is strictly increasing on both branches ($p(x)\geq p^->\alpha$) and it is continuous at $s=\varepsilon$ we infer that the function is \textbf{strictly increasing} on $[0,\infty)$ as needed.

	\item[$\bullet$] $\textbf{a}:\overline{\Omega}\times\mathbb{R}^N\to\mathbb{R}^N$, $\mathbf{a}(x,\xi)=\begin{cases} \mathbf{0}, & \xi=\mathbf{0}\\[3mm] |\xi|^{p(x)-2}\ln^{\delta}(1+|\xi|)\xi, & 0<|\xi|\leq\varepsilon\\[3mm] \varepsilon^{p(x)-\alpha}|\xi|^{\alpha-2}\ln^{\delta}(1+|\xi|)\xi, & |\xi|\geq\varepsilon \end{cases}$.
	
	\item[$\bullet$] $A:\overline{\Omega}\times\mathbb{R}^N\to [0,\infty)$,
	
	\begin{equation}
	A(x,\xi)=\displaystyle\int_{0}^{|\xi|} \Phi(x,s)\ ds=\begin{cases}\displaystyle\int_{0}^{|\xi|}s^{p(x)-1}\ln^{\delta}(1+s) \ ds, & 0\leq |\xi|\leq\varepsilon \\[3mm]\displaystyle\int_{0}^{\varepsilon}s^{p(x)-1}\ln^{\delta}(1+s)  \ ds+\varepsilon^{p(x)-\alpha}\displaystyle\int_{\varepsilon}^{|\xi|}s^{\alpha-1}\ln^{\delta}(1+s) \ ds, & |\xi|\geq\varepsilon \end{cases}.                                                                                                        
	\end{equation}
	
	\item[$\bullet$] $\mathcal{A}:W^{1,p(x)}(\Omega)\to [0,\infty),\ \mathcal{A}(U)=\displaystyle\int_{\Omega} A(x,\nabla U(x))\ dx=\int_{\Omega}\displaystyle\int_{0}^{|\xi|} \Phi(x,s)\ ds$.
	
	\begin{remark}
		It is clear that $\Phi$ is not a $p(x)-1$--homogeneous function, and \textbf{D\'{i}az-Saa inequality} as stated in Theorem \ref{thmdiazsaa} applies to $\mathbf{a}$. The older versions of \textbf{D\'{i}az-Saa inequality} given in \cite{tak}, \cite{giaco1}, \cite{diazsaa}, \cite{brezosw} fail to apply here. Take also a look at Lemma \ref{lemmasimplehomo}.
	\end{remark}
	
	\bigskip
	
	\noindent The last thing we need to check is hypotheses \textbf{(H8)}. We can prove, using the inequality $\ln(1+s)\geq\dfrac{s}{1+s}$ for any $s\geq 0$, after some manipulations that there are two constants $c_1,c_2>0$ such that:
	
	\begin{equation}
		\mathcal{A}(U)\geq c_1\int_{\Omega} |\nabla U|^\alpha \ dx- c_2,\ \forall\ U\in W^{1,p(x)}(\Omega).
	\end{equation}
	
	\noindent So hypothesis \textbf{(H8)} fails to be satisfied. However the results proved here (in particular uniqueness of the strong positive weak solution) are applicable to this example, because hypothesis \textbf{(H8)} was required only for the existence of the solution, as it was proved in \cite{max1}. The existence of the solution for this particular example may follow via other techniques. Using again Theorem \ref{thmuniq} we obtain:
	
	\begin{corollary}\label{corimage}
		Problem \eqref{eqimage} has at most one strong positive solution.
	\end{corollary}
	
\end{enumerate}

\section{Concluding remarks}

\noindent In a future research I will try to state and prove a strong maximum principle for problem \eqref{eqmphase} similar to the one proved in \cite[Theorem 1.1]{qihu}. The strong maximum principle combined with a regularity result (at least continuity up to the boundary or H\"{o}lder continuity up to the boundary)\footnote{Some results of this kind are stated here \cite[Remark 2.4]{Fan6}.} for the same problem will imply that any positive weak solution of \eqref{eqmphase} is zero a.e. on $\Omega$ or it is a strong positive solution, so that the results proven in the present paper will play a decisive role in establishing the number of positive solutions of problem \eqref{eqmphase}.

\section{Appendix}

\begin{lemma}\label{lemrunu}
	For any $a,b\geq 0$ and any $0\leq r\leq 1$ the following inequalities hold:\footnote{See \cite[page 251]{Jones}.}
	
	\begin{enumerate}
		\item[\textnormal{\textbf{(1)}}] $|a^r-b^r|\leq |a-b|^r$
		
		\medskip
		
		\item[\textnormal{\textbf{(2)}}] $(a+b)^r\leq a^r+b^r$.
	\end{enumerate}
	
\end{lemma}

\begin{proof} Without loss of generality we may assume that $b\geq a$. Consider the function $h:[0,\infty)\to\mathbb{R},\ h(x)=(a+x)^r-a^r-x^r$. We see that $h(0)=0$, $h$ is continuous on $[0,\infty)$ and differentiable on $(0,\infty)$ with $h'(x)=r\big [(a+x)^{r-1}-x^{r-1}\big ]\leq 0$. Here we have used that $0\leq r\leq 1$. So $h$ is decreasing on $[0,\infty)$. Thus $h(b-a)\leq h(0)$ from which $b^r-a^r\leq (b-a)^r$, and $h(b)\leq h(0)$ from which $(a+b)^r\leq a^r+b^r$, as desired. 
\end{proof}

	\bibliographystyle{apalike}
	\bibliography{diaz}
\end{document}